\documentclass[12pt]{article}

\usepackage{amssymb,amsmath,amsthm,amsfonts,mathrsfs}
\usepackage{subfigure,epsfig,psfrag,epstopdf}
\usepackage{fullpage,graphicx,color}
\usepackage{algorithm,algorithmic}
\usepackage{tabularx,array,booktabs,bm,multirow}
\usepackage{longtable}
\usepackage{cite,url}
\usepackage[bookmarks=false]{hyperref}
\usepackage[labelfont=bf]{caption}
\usepackage{enumitem}

\newcommand{\mB}{\mathcal{B}}
\newcommand{\mD}{\mathcal{D}}
\newcommand{\mJ}{\mathcal{J}}

\newcommand{\mC}{\mathcal{C}}

\newcommand{\mO}{\mathcal{O}}
\newcommand{\mP}{\mathcal{P}}
\newcommand{\RR}{\mathbb{R}}

\def\bnabla{\boldsymbol{\nabla}}

\def\Dpartial#1#2{\dfrac{\partial #1}{\partial #2}}
\newcommand{\argmin}{\operatorname{argmin}}

\newcommand{\mean}{\operatorname{mean}}
\newcommand*{\mybox}[1]{\framebox{#1}}

\begin{document}

\title{Efficient Gradient-based Optimization for Reconstructing
       Binary Images in Applications to Electrical Impedance Tomography}
\date{}
\author{{\bf \normalsize Paul R.~Arbic II} \\
{\it \small Department of Mathematical Sciences, Florida Institute of Technology} \and
{\bf \normalsize Vladislav Bukshtynov}\footnote{Corresponding author: \url{vbukshtynov@fit.edu}} \\
{\it \small Department of Mathematics and Systems Engineering, Florida Institute of Technology,}\\
{\it \small Melbourne, FL 32901, USA}}
\maketitle

\begin{abstract}

A novel and highly efficient computational framework for reconstructing
binary-type images suitable for models of various complexity seen in
diverse biomedical applications is developed and validated. Efficiency
in computational speed and accuracy is achieved by combining the advantages
of recently developed optimization methods that use sample solutions with
customized geometry and multiscale control space reduction, all paired with
gradient-based techniques. The control space is effectively reduced based
on the geometry of the samples and their individual contributions. The entire
3-step computational procedure has an easy-to-follow design due to a nominal
number of tuning parameters making the approach simple for practical
implementation in various settings. Fairly straightforward methods for
computing gradients make the framework compatible with any optimization
software, including black-box ones. The performance of the complete computational
framework is tested in applications to 2D inverse problems of cancer
detection by electrical impedance tomography~(EIT) using data from
models generated synthetically and obtained from medical images showing
the natural development of cancerous regions of various sizes and shapes.
The results demonstrate the superior performance of the new method and
its high potential for improving the overall quality of the EIT-based
procedures.


\end{abstract}

\section{Introduction}
\label{sec:introduction}

In this paper, we propose a novel computational approach for the
optimal reconstruction of biomedical images to assist with the
recognition of cancer-affected regions while solving an inverse
problem of cancer detection~(IPCD). The main focus is on applying
our new approach to be aligned with the electrical impedance
tomography~(EIT) technique. However, one could easily extend this
methodology to a broad range of problems in biomedical sciences, physics,
geology, chemistry, etc. EIT is a rapidly developing non-invasive
imaging technique gaining popularity during the last few decades by
enabling various medical applications to perform screening for cancer
detection \cite{Zou2003Review,Brown2003,Adler2008,AdlerHolder2022,Uhlmann2009,
Lionheart2004,Abascal2011}. This technique uses a well-known fact that
the electrical properties of different tissues, e.g., electrical
conductivity or permittivity, are different if they are healthy
or affected by cancer. This phenomenon allows EIT to produce images
of biological tissues by interpreting their response to applied
electrical voltages (potentials) or injected currents
\cite{Brown2003,Cheney1999,AdlerHolder2022}. More specifically, the inverse
EIT problem reconstructs the electrical conductivity by measuring
voltages or currents at electrodes placed on the surface of the
tested medium. This so-called Calderon-type inverse problem
\cite{Calderon1980} is highly ill-posed; we refer the reader to
Borcea's topical review paper \cite{Borcea2002}. Since the 1980s, various
computational techniques have been suggested to solve this highly
nonlinear inverse problem computationally; see the recent review
papers \cite{Adler2015,Bera2018,Wang2020} on the current state of
the art and existing open problems associated with EIT and its
applications.

Our particular interest is in creating a robust and computationally
efficient EIT-based optimization framework that is useful in various
applications for medical practices dealing with models characterized
by parameters close to binary-type distributions such as electrical
conductivity. Recent papers
\cite{AbdullaBukshtynovSeif2021,KoolmanBukshtynov2021,ChunEdwardsBukshtynov2024}
propose to convert the inverse EIT statement into a PDE-constrained optimal
control framework and apply multilevel control space reduction at various
scales to improve the quality of the obtained binary images. In
\cite{ArbicBukshtynov2022,ArbicMS2020}, the authors proposed a novel
(although fairly simple) computing algorithm built around a derivative-free
optimization supported by a set of sample solutions. These samples are generated
synthetically with a geometry based on some prior knowledge of the
simulated phenomena and the expected structure of obtained images.
Although the ease of parallelization allows operations on very large
sample sets, enabling the best approximations for the initial guess,
practical applications may be limited to reconstructions of cancerous
spots bearing simple geometry, e.g., with circular shapes. This
algorithm also requires some prior knowledge of the simulated
phenomena in the form of approximated properties (electrical
conductivity) for healthy tissues and regions affected by cancer.

In the current work, we propose a generalized 3-step optimization
procedure that removes these limitations, making the computational
framework applicable to models of various levels of complexity.
The superior performance of this generalized methodology is achieved
by adding the following new computational components.
\begin{enumerate}
  \item[(1)] We designed new methods for computing gradients and validating
    their correctness to further enhance the computational performance at
    Step~2, as discussed in Section~\ref{sec:step_2}. These gradients are
    computed with respect to all control variables used in the main optimization
    cycle (all three steps).
  \item[(2)] We also discussed the applicability of various gradient-based methods
    (optimizers) to further enhance the optimization performance.
  \item[(3)] The limitations of our simplified 2-step approach presented in
    \cite{ArbicBukshtynov2022,ArbicMS2020} are overcome by the addition of a new
    (third) step in the current procedure, namely a binary tuning optimization,
    as discussed in Section~\ref{sec:step_3}. We applied the multiscale-based
    post-processing filtering \cite{ChunMS2022,ChunEdwardsBukshtynov2024} at
    the coarse scale to further enhance the quality of images obtained as the
    final output of our new optimization procedure.
\end{enumerate}

The proposed computational framework has an easy-to-follow design
split into three phases and tuned by a nominal (meaning, minimal)
number of computational
parameters making the approach simple for practical implementation
for various applications far beyond biomedical imaging.
We tested the performance of our new algorithm computationally
in applications to 2D inverse problems of cancer detection using noisy
data generated synthetically and obtained from medical images showing
cancerous regions of various sizes and shapes developed naturally in
the human body.

This paper proceeds as follows. Section~\ref{sec:math} presents
a general mathematical description of the inverse EIT problem
formulated as an optimal control problem. The entire computational
procedure for solving this optimization problem is discussed in Section~\ref{sec:solution}. Model descriptions and detailed computational
results, including a discussion of chosen methods, are presented in Section~\ref{sec:results}. Concluding remarks are provided in
Section~\ref{sec:remarks}.

\section{Inverse EIT Problem}
\label{sec:math}

\subsection{Mathematical Model for Optimization}
\label{sec:math_opt}

As discussed at length in
\cite{AbdullaBukshtynovSeif2021,KoolmanBukshtynov2021,ArbicBukshtynov2022,ChunEdwardsBukshtynov2024},
the inverse EIT problem is formulated as a PDE-constrained optimal
control problem. Here, we note that while we refer to the control theory
throughout the entire paper, in a very general sense, we can treat this
formulation as an example of a PDE-constrained optimization problem
defined for an open and bounded set (domain) $\Omega \subset \RR^n,
\ n = 2, 3$, representing a body with electrical conductivity at point
$x \in \Omega$ given by function $\sigma(x): \, \Omega \rightarrow
\RR_+$. In this paper, we use the so-called ``voltage--to--current"
model where voltages (electrical potentials) $U = (U_{\ell})_{\ell=1}^m
\in \RR^m$ are applied to $m$ electrodes $(E_{\ell})_{\ell=1}^m$ with
contact impedances $(Z_{\ell})_{\ell=1}^m \in \RR^m_+$ subject to the
ground (zero potential) condition
\begin{equation}
  \sum_{\ell=1}^m U_{\ell}=0.
  \label{eq:ground_conds}
\end{equation}
These voltages initiate electrical currents $(I_{\ell})_{\ell=1}^m \in
\RR^m$ through the same electrodes $E_{\ell}$ placed at the periphery
$\partial \Omega$ of domain $\Omega$. We compute the electrical currents
\begin{equation}
  I_{\ell} = \int_{E_{\ell}}  \sigma(x) \Dpartial{u(x)}{n} \, ds,
  \quad \ell = 1, \ldots, m
  \label{eq:el_current}
\end{equation}
based on conductivity field $\sigma(x)$ and a distribution of electrical
potential $u(x): \, \Omega \rightarrow \RR$ obtained by solving the
following elliptic problem
\begin{subequations}
  \begin{alignat}{3}
    \bnabla \cdot \left[ \sigma(x) \bnabla u(x) \right] &= 0,
    && \quad x \in \Omega \label{eq:forward_1}\\
    \Dpartial{u(x)}{n} &= 0, && \quad x \in \partial \Omega -
    \bigcup\limits_{\ell=1}^{m} E_{\ell} \label{eq:forward_2}\\
    u(x) +  Z_{\ell} \sigma(x) \Dpartial{u(x)}{n} &= U_{\ell},
    && \quad x \in E_{\ell}, \ \ell= 1, \ldots, m
    \label{eq:forward_3}
  \end{alignat}
  \label{eq:forward}
\end{subequations}
in which $n$ is an external unit normal vector on $\partial \Omega$.
A complete description and analysis of electrode models used
in electric current computed tomography may be found, e.g., in
\cite{Somersalo1992}.

We set conductivity $\sigma(x)$ in \eqref{eq:forward} as a control variable
(or optimization variable; see our remark at the beginning of this section)
and formulate the inverse EIT (conductivity) problem
\cite{Calderon1980} as a PDE-constrained optimal control problem
\cite{AbdullaBukshtynovSeif2021} by considering least-square
minimization of mismatches $\left( I_{\ell} - I_{\ell}^* \right)^2$,
where $(I_{\ell}^*)_{\ell=1}^m \in \RR^m$ are measurements of electrical
currents $I_{\ell}$. In addition, we have to mention a well-known
fact that this inverse EIT problem to be solved in a discretized
domain $\Omega$ is highly ill-posed. Therefore, we enlarge the data
up to a size of $m^2$ by adding new measurements following the ``rotation
scheme'' described in detail in \cite{AbdullaBukshtynovSeif2021}
while keeping the size of the unknown parameters, i.e., elements in the
discretized description for $\sigma(x)$, fixed. Having a new set of
data $\mD^* = \left\{ I_{\ell}^{k*} \right\} _{\ell,k=1}^m$ and in
light of the Robin condition~\eqref{eq:forward_3} used together with
\eqref{eq:el_current}, we define a complete form of the cost functional
\begin{equation}
  \mJ (\sigma) = \sum_{k=1}^{m} \sum_{\ell=1}^m
  \left[ \int_{E_{\ell}} \dfrac{U^k_{\ell}-u^k(x;\sigma)}{Z_{\ell}}
  \, ds - I^{k*}_{\ell} \right]^2
  \label{eq:cost_functional}
\end{equation}
for the optimal control problem
\begin{equation}
  \hat \sigma(x) = \underset{\sigma}{\argmin} \ \mJ(\sigma)
  \label{eq:minJ_sigma}
\end{equation}
subject to PDE constraint as each function $u^k(\cdot; \sigma), \ k = 1,
\ldots, m$, solves elliptic PDE (forward EIT) problem \eqref{eq:forward}.
Here, we assume that optimal solution $\hat \sigma(x)$ in \eqref{eq:minJ_sigma}
is unique. We also note that after applying \eqref{eq:el_current} and adding some
noise, these solutions may be used for generating various model examples
(synthetic data) for inverse EIT problems to adequately mimic the
presence of regions affected by cancer we mentioned above.

\subsection{Solution by Sample-based Parameterization}
\label{sec:math_sol}

In the view of the binary type of a solution we seek to reconstruct,
we assume that the actual (true) electrical conductivity
$\sigma_{true}(x)$ is represented by
\begin{equation}
  \sigma_{true}(x) = \left\{
  \begin{aligned}
    \sigma_c, & \quad x \in \Omega_c,\\
    \sigma_h, & \quad x \in \Omega_h,
  \end{aligned}
  \right.
  \quad \Omega_c \cap \Omega_h = \emptyset,
  \label{eq:sigma_true}
\end{equation}
where $\sigma_c$ and $\sigma_h$ are some constants for the respective
cancer-affected region $\Omega_c$ and the healthy tissue part $\Omega_h$.
As such, we seek the solution of \eqref{eq:minJ_sigma} in a form
of the linear (convex) combination
\begin{equation}
  \sigma(x) = \sum_{i=1}^{N_s} \alpha_i \bar \sigma_i(x),
  \qquad 0 \leq \alpha_i \leq 1, \qquad  \sum_{i=1}^{N_s} \alpha_i = 1,
  \label{eq:sigma_main}
\end{equation}
where $\bar \sigma_i(x)$, $i = 1, \ldots, N_s$, are sample solutions
generated synthetically based on a general assumption for the solution
structure provided in \eqref{eq:sigma_true}; see Section~\ref{sec:step_0}
for details. The proposed computational algorithm for solving optimization
problem \eqref{eq:minJ_sigma} could be executed in three steps and requires
a preprocessing phase (Step~0), as shown briefly below and detailed in
Section~\ref{sec:solution}.
\begin{itemize}
  \item[{\bf Step~0:}] Preprocessing includes choosing (approximating) constants
    $\sigma_c$ and $\sigma_h$ in \eqref{eq:sigma_true}, deciding on the geometry
    of sample solutions $\bar \sigma_i$, and generating a collection $\mC(N)$ of
    $N$ samples with corresponding ``measurements'' obtained by using
    \eqref{eq:el_current}--\eqref{eq:forward}; see Section~\ref{sec:step_0} for
    details.
  \item[{\bf Step~1:}] The initial basis of samples
    \begin{equation}
      \mB^0 = \left\{ \bar \sigma_i(x) \right\}_{i=1}^{N_s} \subset \mC(N)
      \label{eq:smpl_basis}
    \end{equation}
    is defined by choosing the $N_s$ ``best'' samples out of $\mC(N)$ that
    provide the best measurement fit in terms of cost functional
    \eqref{eq:cost_functional} (we elaborate more on this initialization
    procedure in Section~\ref{sec:step_1}). This basis will serve as the
    initial guess for the fine-scale optimization performed in Step~2.
  \item[{\bf Step~2:}] All parameters in the description of basis $\mB^0$
    and weights $\alpha_i$ in \eqref{eq:sigma_main}--\eqref{eq:smpl_basis}
    are set as controls to perform gradient-based
    optimization for solving problem~\eqref{eq:minJ_sigma} numerically to
    find the optimal solution $\hat \sigma(x)$ at a {\it fine scale} (i.e., using
    a fine mesh on which we reconstruct images for electrical conductivity $\sigma(x)$)
    via optimal control set $(\hat \mB, \hat \alpha)$ as shown in Section~\ref{sec:step_2}.
  \item[{\bf Step~3:}] The fine-scale optimal solution $\hat \sigma(x)$ is
    then finally tuned to create an optimal binary image $\hat \sigma_b(x)$
    by moving $\hat \sigma(x)$ from the {\it fine} to the {\it coarse scale}
    (where $\hat \sigma_b(x)$ leaves) using a multiscale optimization algorithm.
    Here, we use the ``coarse scale'' notation to refer to an upscaled (reduced-dimensional)
    control space containing just a few controls representing reconstructed values for high
    and low electrical conductivities inside each cancerous region and outside, respectively.
    In short, this step will tune the fine-scale solution obtained in Step~2
    by recreating it as a piecewise constant reconstruction; refer to Section~\ref{sec:step_3}
    for more details.
\end{itemize}
To comment more on the ``philosophy'' behind the entire approach, we would reiterate that
the fine-scale optimization of Step~2 is used to approximate the location of regions
with high and low values of electrical conductivity $\sigma(x)$, and projecting
obtained solutions onto the coarse scale in Step~3 provides a dynamical (sharp-edge)
filtering to the fine-scale images optimized to better represent the structure
of the cancer-affected regions (i.e., their location and boundaries).

\section{Gradient-based Optimization Framework}
\label{sec:solution}

\subsection{Step~0: Sample Preprocessing}
\label{sec:step_0}

Without loss of generality, in this paper, we discuss the application
of the new algorithm to solving optimization problem \eqref{eq:minJ_sigma}
in the 2D ($n = 2$) domain, e.g.,
\begin{equation}
  \Omega = \left\{ x \in \RR^2 : \ \| x \| ^2 < R^2 \right\},
  \label{eq:domain}
\end{equation}
which is a disc of radius $R$. However, the same concept could be
easily extended to 3D ($n = 3$) regions of various complexity.

The entire collection of $N$ samples
\begin{equation}
  \mC(N) = \left\{ \bar \sigma_i(x) \right\}_{i=1}^{N}, \qquad N \gg N_s
  \label{eq:smpl_collect}
\end{equation}
could be generated based on various assumptions made for the
(geometrical) structure of the reconstructed images with the binary
type restriction. Here, we assume that flexibility in reconstructing
images of various complexity and also convenient simplicity could be
achieved by combining simple convex geometric shapes (elements) in 2D
such as triangles, squares, circles, etc. For example, in this paper,
the $i$th sample in $\mC(N)$ consists of $N_c^i$ circles of various
radii $r \in \RR_+$ and centers $x^0 = (x^{01}, x^{02}) \in \RR^2$
located inside domain $\Omega$, i.e.,
\begin{equation}
  \bar \sigma_i(x) = \left\{
  \begin{aligned}
    \tilde \sigma_c, &\quad \mid x - x^0_j  \mid ^2 \leq r^2_j,
    \quad j = 1, \ldots, N_c^i\\
    \tilde \sigma_h, &\quad {\rm otherwise}
  \end{aligned}
  \right.
  \label{eq:smpl_par}
\end{equation}
where some approximations $\tilde \sigma_c$ and $\tilde \sigma_h$
for respective $\sigma_c$ and $\sigma_h$ in \eqref{eq:sigma_true}
are required and considered as {\it a priori} knowledge needed for
applying the approach in practice. In \eqref{eq:smpl_par}, all $N_c^i$
circles (i.e., cancer-affected regions) are parameterized
by the set of triplets
\begin{equation}
  \mP_i = \left\{ (x^{01}_j, x^{02}_j, r_j) \right\}_{j=1}^{N_c^i},
  \qquad i = 1, \ldots, N
  \label{eq:smpl_triplet}
\end{equation}
generated randomly subject to the following restrictions
\begin{equation}
  \begin{aligned}
    &&\mid x^0_j \mid < R+r_j, \qquad j &= 1, \ldots, N_c^i,\\
    &&1 \leq N_c^i \leq N_{c,\max}, \qquad i &= 1, \ldots, N.
  \end{aligned}
  \label{eq:smpl_restr}
\end{equation}
Parameter $N_{c,\max}$ in \eqref{eq:smpl_restr} defines the
maximum number of circles in the samples and, in fact, sets the
highest level of complexity (resolution) for the reconstructed
images (refer to \cite{ArbicMS2020,ArbicBukshtynov2022} for more
details on technicalities of this preprocessing step).

Finalizing the preprocessing procedure requires solving forward
problem~\eqref{eq:forward} and evaluating cost
functional~\eqref{eq:cost_functional} $N$ times for all samples in
$\mC(N)$. Using a fixed scheme of potentials $U$, the entire
``measurement'' data $\mD = \left\{ \mD_i \right\}_{i=1}^N$, where
$\mD_i = I^{k}_{\ell}(\bar \sigma_i) \in \RR^{m^2}$, are precomputed by
\eqref{eq:el_current}--\eqref{eq:forward} and then stored for
multiple uses with different models. In addition, this task may
run in parallel with minimal computational overhead, which
allows easy switching between various schemes for electrical
potentials. Easy parallelization enables taking $N$ to be quite large,
which helps better approximate the solution by the initial state
of basis $\mB$ during Step~1 before proceeding to Step~2.

\subsection{Step~1: Initializing Sample Basis}
\label{sec:step_1}

We set the number of samples $N_s$ in the initial basis $\mB^0$
as a hyper parameter of the algorithm and define it heuristically after
making assumptions on the model complexity or after experimentation. We suggest
$N_s$ be sufficiently large to properly support a local/global
search for optimal solution $\hat \sigma(x)$ during Step~2. At the
same time, while solving problem \eqref{eq:minJ_sigma}, this number
should allow the total number of controls, to be comparable with
the size of the data, namely $m^2$, for satisfying the well-posedness
requirement for the solution of the optimization problem in the sense
of Hadamard \cite{Hadamard1923}.

Considering models with highly complicated structures may require
increasing the number of elements (in our case, circles) in every
sample within the chosen basis $\mB^0$. In this case,
one could re-set parameters $N_c^i, \ i = 1, \ldots, N_s$, to higher
values and add missing elements, for example, by generating randomly
new circles. It will project the initial basis $\mB^0$ onto a new
control space of a higher dimension to minimize the loss in the quality
of the initial solution $\sigma(\mB^0,\alpha^0)$.

Step~1 will be completed after ranking all samples in $\mC(N)$ in
ascending order using computed cost
functionals~\eqref{eq:cost_functional} while comparing the obtained
data $\mD$ with true data $\mD^*$ available from the actual
measurements. After ranking, the first $N_s$ samples create
the initial basis $\mB^0$ to construct solution $\sigma(\mB^0,\alpha^0)$
used as the initial guess for optimization in Step~2.

\subsection{Step~2: Fine-Scale Gradient-based Optimization}
\label{sec:step_2}

As discussed in Section~\ref{sec:step_0}, all elements (circles) in
all samples of basis $\mB^0$ obtained during the Step~1 ranking procedure
are represented by a finite number of ``sample-based'' parameters
associated with set $\left\{ \mP_i \right\}_{i=1}^{N_s}$. In general,
solution $\sigma(x) = \sigma(\mP, \alpha)$ could be uniquely represented
as a function of set $\mP = \left\{ \mP_i \right\}_{i=1}^{N_s}$
and vector of weights $\alpha = (\alpha_i)_{i=1}^{N_s}$. We substitute
the continuous form of optimal control problem \eqref{eq:minJ_sigma}
with its new equivalent form
\begin{equation}
  (\hat \mP, \hat \alpha) = \underset{\mP,\alpha}{\argmin}
  \ \mJ(\mP,\alpha)
  \label{eq:minJ_ext}
\end{equation}
to be solved numerically
subject to PDE constraint \eqref{eq:forward}, linear constraints for
$\alpha_i$ in \eqref{eq:sigma_main}, and suitably established bounds
for all components of the combined control set $(\mP, \alpha)$.
As naturally followed from the structure of this new control, a
dimension of the parameterized solution space is bounded by
\begin{equation}
  \dim(\mP, \alpha) \leq N_s \cdot [ N_{c,\max} (n+1) + 1].
\end{equation}
To solve \eqref{eq:minJ_ext} iteratively, one may choose various
criteria to terminate the optimization run at the $k$th iteration,
e.g., comparing the relative decrease in the cost functional $\mJ^k$
evaluated after completing $k$ iterations
\begin{equation}
    \dfrac{\mid \mJ^k - \mJ^{k-1} \mid}{\mJ^k}
     < \epsilon, \qquad k > 0
  \label{eq:termination}
\end{equation}
with preset tolerance $\epsilon \in \RR_+$.

The prior works on using sample-based parameterization
\cite{ArbicBukshtynov2022,ArbicMS2020} employed the coordinate descent~(CD)
method to solve \eqref{eq:minJ_ext}. While achieving a good performance
when applied to simple models, CD obviously exhibits certain limitations
as it optimizes over a single control at a time. Without proper
parallelization, it requires enormous cost functional evaluations and,
as such, a large amount of computational time to complete optimization
with the reasonably small $\epsilon$. Instead, our new optimization
framework operates with fairly straightforward methods for computing
gradients derived with respect to all controls in the control set
$(\mP, \alpha)$.

We start our derivation of gradients relative to both sample-based
parameters in set $\mP$ and weights $\alpha$ by referring to the
known structure
\cite{ChunEdwardsBukshtynov2024,KoolmanBukshtynov2021,AbdullaBukshtynovSeif2021}
of gradients $\mJ'_{\sigma}$ obtained with respect to control $\sigma$

\begin{equation}
  \mJ'_{\sigma} = - \sum_{k=1}^{m} \bnabla \psi^k(x)
  \cdot \bnabla u^k(x)
  \label{eq:grad_sigma_spat}
\end{equation}
computed based on solutions $\psi^k(\cdot; \sigma): \, \Omega \rightarrow
\RR, \ k = 1, \ldots, m$, of the following adjoint PDE problem
\begin{equation}
  \begin{aligned}
    \bnabla \cdot \left[ \sigma(x) \bnabla \psi(x) \right] &= 0,
    && \quad x \in \Omega\\
    \Dpartial{\psi(x)}{n} &= 0, && \quad x \in \partial \Omega -
    \bigcup\limits_{\ell=1}^{m} E_{\ell}\\
    \psi(x) + Z_{\ell} \Dpartial{\psi(x)}{n} &=
    2 \beta_{\ell} \left[ \int_{E_{\ell}}
    \dfrac{u(x)-U_{\ell}}{Z_{\ell}} \, ds
    + I^*_{\ell} \right], && \quad
    x \in E_{\ell}, \\
    &&& \quad\ell= 1, \ldots, m
  \end{aligned}
  \label{eq:adjoint}
\end{equation}

First, we derive the gradient $\bnabla_{\alpha} \mJ$ of cost
functional $\mJ$ with respect to control $\alpha$. Using connectivity
of solution $\sigma(x)$ with individual sample $\bar \sigma_i$ weights
$\alpha_i$ provided explicitly by \eqref{eq:sigma_main}, partial
derivatives
\begin{equation}
  \Dpartial{\sigma(x)}{\alpha_i} = \bar \sigma_i, \quad i = 1, \ldots, N_s
  \label{eq:grad_sigma_i}
\end{equation}
are then used to construct the gradient
\begin{equation}
  \bnabla_{\alpha} \sigma = \left[ \bar \sigma_1 \ \bar \sigma_2 \
  \ldots \bar \sigma_{N_s} \right]^T \triangleq \bar \sigma.
  \label{eq:grad_sigma}
\end{equation}
Then using the chain rule gives the sought gradient
\begin{equation}
  \bnabla_{\alpha} \mJ = \bnabla_{\alpha} \sigma \cdot
  \mJ'_{\sigma} = \bar \sigma \cdot \mJ'_{\sigma}.
  \label{eq:grad_alpha}
\end{equation}
We also note that in the discretized settings (when domain $\Omega$
undergoes $N_{\Omega}$-component discretization), $\bar \sigma$ and
$\mJ'_{\sigma}$ are represented by $N_s \times N_{\Omega}$
matrix and $N_{\Omega}$-component column-vector, respectively.

Finally, the gradient $\bnabla_{\mP} \mJ$ of cost functional $\mJ$
with respect to sample-based control $\mP$ is derived in the same
manner by using the chain rule and precomputed gradient $\mJ'_{\sigma}$
\begin{equation}
  \bnabla_{\mP} \mJ = \bnabla_{\mP} \sigma \cdot \mJ'_{\sigma}.
  \label{eq:grad_P}
\end{equation}
Deriving $\bnabla_{\mP} \sigma$ appears complicated due to the involved
geometry of $\sigma(\mP, \alpha)$ that may neither have a straightforward
structure nor be known due to the randomness of the process described in
\eqref{eq:smpl_par}--\eqref{eq:smpl_restr}. However, recent studies
\cite{Volkov2018,Krogstad2022} suggest a flexible approach when some
partial derivatives used as a part of the adjoint-based analysis are
approximated by numerical perturbations. We could conveniently adapt this
approach as, in the current framework, it will not require reevaluating
cost functionals: only changes in the sample solution $\bar \sigma_i$
associated with its own parameter $\mP_i$ should be assessed. Thus, we
perturb every parameter $\mP_i$ in control set $\mP$
\begin{equation}
  \Dpartial{\sigma}{\mP_i} \approx \dfrac{\Delta \sigma}{\Delta \mP_i}
  = \alpha_i \dfrac{\Delta \bar \sigma_i}{\Delta \mP_i}, \quad i = 1,
  \ldots, N_s
  \label{eq:grad_P_i}
\end{equation}
by setting all perturbations $\Delta \mP_i$ to numerical values pursuing
a trade-off between being reasonably small to ensure the accuracy of
finite-difference~(FD) estimations of ${\partial \sigma}/{\partial \mP_i}$
and large enough to protect the numerator from being zero.

\subsection{Step~3: Coarse-Scale Binary Tuning Optimization}
\label{sec:step_3}

For the last step, we employ a recently designed gradient-based approach to
support multiscale optimization with multilevel control space reduction
using principal component analysis~(PCA) coupled with dynamical control
space upscaling \cite{KoolmanBukshtynov2021,ChunEdwardsBukshtynov2024}.
As pointed out before, this step introduces projecting the fine-scale solutions
$\hat \sigma(x)$ onto the coarse scale (i.e., new upscaled control space with
a significantly reduced number of controls) to perform a dynamical (sharp-edge)
filtering to the fine-scale images. New images represented by the coarse-scale
solutions $\hat \sigma_b(x)$ are optimized to better represent the structure
of the cancer-affected regions (i.e., their location and boundaries).
From the entire approach presented in \cite{KoolmanBukshtynov2021,ChunEdwardsBukshtynov2024},
we adopt only the coarse-scale phase to be applied to the fine-scale solution
$\hat \sigma(x)$ obtained during Step~2.

First, we must specify the maximum (expected) number of cancer-affected
(high conductivity) regions $N_{\max}$. By considering the healthy part
(low conductivity region $\Omega_h$) of domain $\Omega$ as a single region,
partitioning the fine mesh elements representing $\Omega$ will create
$N_{\zeta} = N_{\max} + 1$ spatial subsets. To proceed with Step~3
optimization at the coarse scale, we define a new control vector $\zeta =
(\zeta_j)_{j=1}^{2N_{\max}+1}$ of significantly reduced dimensionality,
in which the first entry is the low value of
(binary) electrical conductivity $\sigma (x)$ associated with a healthy
region $\Omega_h$. The next $N_{\max}$ controls are the high values
of $\sigma (x)$ related to areas in $\Omega_c$ affected by cancer, i.e.,
\begin{equation}
  \zeta_1 = \sigma_{low}, \quad \zeta_2 = \sigma_{high,1},
  \quad \zeta_3 = \sigma_{high,2}, \quad \ldots, \quad
  \zeta_{N_{\max}+1} = \sigma_{high,N_{\max}}.
  \label{eq:ctrls_low_high}
\end{equation}
The rest $N_{\max}$ components
\begin{equation}
  \zeta_{N_{\max}+2} = \sigma_{th,1}, \quad
  \zeta_{N_{\max}+3} = \sigma_{th,2}, \quad \ldots, \quad
  \zeta_{2N_{\max}+1} = \sigma_{th,N_{\max}}
  \label{eq:ctrls_thresh}
\end{equation}
take responsibility for the shape of those $N_{\max}$ cancerous regions.
They are set as separation thresholds to define boundaries between the low
and high-conductivity areas. The structure of control $\zeta$ allows the
creation of the systematic representation of the coarse-scale solution
$\zeta^k$ for control $\sigma^k$ at the $k$th iteration based on the current
fine-scale parameterization $\sigma(\zeta^k) = (\sigma_i(\zeta^k))_{i=1}^{N_{\Omega}}$,
i.e,
\begin{equation}
  \sigma^k_i = \left\{
  \begin{aligned}
    \sigma^k_{low},  & \quad \sigma_i(\zeta^k) < \sigma^k_{th,n},\\
    \sigma^k_{high,n}, & \quad \sigma_i(\zeta^k) \geq \sigma^k_{th,n},
  \end{aligned}
  \right. \quad i = 1, \ldots, N_{\Omega}, \quad 1 \leq n \leq N_{\max}.
  \label{eq:sigma_coarse}
\end{equation}
Here, $n = n(i)$ denotes the number of a particular cancer-affected region
defined subject to the partitioning map currently established and
used for the $k$th iteration \cite{ChunEdwardsBukshtynov2024}.
We also note that
\begin{equation}
  \begin{aligned}
    0 < \sigma^k_{low} &< \underset{1 \leq n \leq N_{\max}}{\min} \ \sigma^k_{high,n}, \\
    \underset{1 \leq i \leq N_{\Omega}}{\min} \ \sigma_i(\zeta^k) < \sigma^k_{th,n}
    &< \underset{1 \leq i \leq N_{\Omega}}{\max} \ \sigma_i(\zeta^k), \quad
    n = 1, \ldots, N_{\max}.
  \end{aligned}
  \label{eq:sigma_coarse_bounds}
\end{equation}
Simply, \eqref{eq:sigma_coarse} provides a rule for creating a
fine--to--coarse partition of discretized fine-scale solution $\hat \sigma$
with all spatial elements $\hat \sigma_i$ belonging either to cancer
affected region $\Omega_c$ ($\hat \sigma_i = \sigma_{high,n}$,
$n = 1, \ldots N_{\max}$) or the healthy tissue part $\Omega_h$
($\hat \sigma_i = \sigma_{low}$) based on the current state of control
$\zeta^k$ (at the $k$th iteration).

During the Step~3 (coarse-scale optimization) phase, control $\sigma^k$ is
updated by solving the following ($2N_{\max}+1$)-dimensional optimization
problem in the $\zeta$-space
\begin{equation}
  \hat \zeta =
  \underset{\zeta}{\argmin} \ \mJ(\zeta)
  \label{eq:minJ_zeta}
\end{equation}
subject to constraints (bounds) provided in \eqref{eq:sigma_coarse_bounds}
and optimal (binary) solution $\hat \sigma_b (x) = \sigma(\hat \zeta)$.
To solve \eqref{eq:minJ_zeta} by any approaches that
require computing gradients, their first $N_{\max} + 1$ components
could be easily obtained by using a gradient summation formula
derived in \cite{KoolmanBukshtynov2021} and upgraded in
\cite{ChunEdwardsBukshtynov2024} for the case of multiple cancer-affected
regions
\begin{equation}
  \Dpartial{\mJ(\zeta)}{\zeta_j} =
  \sum_{i=1}^{N_{\Omega}} P_{i,j} \mJ'_{\sigma,i} \Delta_i,
  \quad j = 1, \ldots, N_{\max}+1.
  \label{eq:grad_zeta_1_2}
\end{equation}
Here, $\mJ'_{\sigma,i}$ is the $i$th component of the discretized gradient
$\mJ'_{\sigma}$, and $P_{i,j}$ is the partitioning (indicator) function
defined by
\begin{equation}
  P_{i,j} = \left\{
  \begin{aligned}
    1, \quad \sigma_i &\in C_j,\\
    0, \quad \sigma_i &\notin C_j,
  \end{aligned}
  \right.
  \label{eq:part_ind}
\end{equation}
after completing the partitioning map by employing \eqref{eq:sigma_coarse}.
In \eqref{eq:part_ind}, $C_1$ is the current representation of the
healthy region $\Omega_h$, and $C_{n+1}$, $n = 1, \ldots, N_{\max}$, are
parts of the cancerous area $\Omega_c$.

The rest $N_{\max}$ components may be approximated by a finite difference
scheme, e.g., of the first order:
\begin{equation}
  \begin{aligned}
    \Dpartial{\mJ(\zeta)}{\zeta_j} &=
    \dfrac{1}{\delta_{\zeta}} \left[ \mJ
    \left(\sigma^k(\ldots,\zeta_j +
    \delta_{\zeta}, \ldots) \right) - \mJ
    \left( \sigma^k(\ldots,\zeta_j,\ldots)
    \right)\right] + \mO(\delta_{\zeta}),\\
    n &= 1, \ldots, N_{\max}, \quad j = n + N_{\max} +1
  \end{aligned}
  \label{eq:grad_zeta_b}
\end{equation}
which requires at most $N_{\max}$ extra cost functional evaluations per optimization
iteration. Following the same discussion as in Section~\ref{sec:step_2},
parameter $\delta_{\zeta}$ in \eqref{eq:grad_zeta_b} may be defined
experimentally, pursuing a trade-off between being reasonably small to
ensure accuracy and large enough to protect the numerator from being zero.

In fact, formulas \eqref{eq:sigma_coarse}--\eqref{eq:minJ_zeta}
provide a complete description of Step~3 fine--to--coarse projection
(or, as we call it before, coarse-scale binary tuning) for fine-scale optimal
control $\hat \sigma(x)$ to obtain optimal (in terms of fitting to data $\mD^*$)
binary distribution $\hat \sigma_b(x)$. An import remark should be made here for
the performance of the binary tuning optimization of Step~3. As shown in multiple
examples in \cite{KoolmanBukshtynov2021,ChunEdwardsBukshtynov2024}, the ``quality''
of the coarse-scale solution $\hat \sigma_b$ is usually worse compared to the quality
of the fine-scale images $\hat \sigma$ if evaluated in terms of fitting to the
same data $\mD^*$. This effect is expected as the coarse scale operates with
significantly reduced number of parameters (components of Step~3 control variable $\zeta$),
reconstructing in fact piecewise constant representation of electrical conductivity $\sigma(x)$
after applying sharp-edge filtering to better represent the structure of the cancer-affected
regions in terms of their location and boundaries. Such solutions could lose even more
``quality'' if the used data contains noise.

The efficiency of the entire optimization framework is confirmed
by extensive computational results for multiple models of different
complexity presented in Section~\ref{sec:results}. A summary of the
complete computational framework to perform our new optimization with
sample-based parameterization and multiscale control-space reduction
for binary tuning is provided in Algorithm~\ref{alg:main_opt}.

\begin{algorithm}
\begin{algorithmic}
  \STATE set parameters: $N$, $N_{c,\max}$, $N_s$
  \STATE \begin{center} \mybox{Step~0}
         \it{Sample Preprocessing} \end{center}
  \FOR{$i \leftarrow 1$ to $N$}
    \STATE generate $\bar \sigma_i(x)$ by
    \eqref{eq:smpl_par}--\eqref{eq:smpl_restr}
    \STATE obtain data $\mD_i = I^{k}_{\ell}(\bar \sigma_i)$ from
      sample $\bar \sigma_i$ by \eqref{eq:el_current}--\eqref{eq:forward}
  \ENDFOR
  \STATE \begin{center} \mybox{Step~1}
         \it{Initializing Sample Basis} \end{center}
  \STATE select model and obtain true data $\mD^*$
  \FOR{$i \leftarrow 1$ to $N$}
    \STATE compute $\mJ(\bar \sigma_i)$ by \eqref{eq:cost_functional}
  \ENDFOR
  \STATE choose $N_s$ best samples from $\mC(N)$ by values
    $\mJ(\bar \sigma_i)$
  \STATE form initial basis $\mB^0$
  \STATE set initial weights $\alpha^0$
  \STATE compute $\sigma^0(x)$ using $\mB^0$ and $\alpha^0$ by
    \eqref{eq:sigma_main}
  \STATE \begin{center} \mybox{Step~2}
         \it{Fine-Scale Optimization} \end{center}
  \STATE $k \leftarrow 0$
  \REPEAT
  \STATE compute state $u^k$ by solving \eqref{eq:forward}
  \STATE compute adjoint state $\psi^k$ by solving \eqref{eq:adjoint}
  \STATE compute gradient $\mJ'_{\sigma}$ using \eqref{eq:grad_sigma_spat}
  \STATE obtain gradients $\bnabla_{\mP} \mJ$ and $\bnabla_{\alpha} \mJ$
    respectively by \eqref{eq:grad_P}--\eqref{eq:grad_P_i} and
    \eqref{eq:grad_sigma_i}--\eqref{eq:grad_alpha}
  \STATE update control set $(\mP, \alpha)$ using gradients
    $\bnabla_{\mP} \mJ$ and $\bnabla_{\alpha} \mJ$
  \STATE $k \leftarrow k + 1$
  \STATE update $\sigma^k(x)$ using new basis $\mB^k = \mB(\mP^k)$ and
    weights $\alpha^k$ by \eqref{eq:sigma_main}
  \UNTIL termination criterion \eqref{eq:termination} is satisfied
    to given tolerance $\epsilon$
  \STATE obtain optimal solution $\hat \sigma(x) =
    \sigma (\hat \mB, \hat \alpha)$
  \STATE \begin{center} \mybox{Step~3}
         \it{Coarse-Scale Binary Tuning} \end{center}
  \STATE $k \leftarrow 0$
  \STATE $\sigma^0_i \leftarrow (\hat \sigma_i)_{i=1}^{N_{\Omega}}$
  \STATE define $\zeta^0$ by
    \eqref{eq:sigma_coarse}--\eqref{eq:sigma_coarse_bounds}
  \REPEAT
  \STATE compute gradient $\bnabla_{\zeta} \mJ$ using
    \eqref{eq:grad_zeta_1_2}--\eqref{eq:grad_zeta_b}
  \STATE update control $\zeta$ using gradient $\bnabla_{\zeta} \mJ$
  \STATE $k \leftarrow k + 1$
  \STATE update $\sigma_i^k$ using \eqref{eq:sigma_coarse}
  \UNTIL termination criterion \eqref{eq:termination} is satisfied
    to given tolerance $\epsilon$
  \STATE obtain optimal solution $\hat \sigma_b = \sigma_i (\hat \zeta)$
\end{algorithmic}
\caption{Computational workflow for gradient-based optimization in
  binary images reconstruction}
\label{alg:main_opt}
\end{algorithm}

\section{Computational Results}
\label{sec:results}

\subsection{2D Model Setup}
\label{sec:comp_model}

Our optimization framework integrates computational facilities for
solving forward PDE problem \eqref{eq:forward}, adjoint PDE problem
\eqref{eq:adjoint}, evaluating cost functionals by
\eqref{eq:cost_functional}, and constructing gradients according
to \eqref{eq:grad_sigma_spat}, \eqref{eq:grad_alpha}, \eqref{eq:grad_P},
and \eqref{eq:grad_zeta_1_2}--\eqref{eq:grad_zeta_b}. These facilities
are incorporated using {\tt FreeFEM} \cite{FreeFem2012}, an open-source,
high-level integrated development
environment for obtaining numerical solutions of PDEs based on the
finite element method~(FEM). For solving numerically forward PDE
problem \eqref{eq:forward}, spatial discretization is carried out
by implementing FEM triangular finite elements: P2 piecewise
quadratic (continuous) and P0 piecewise constant representations
for electrical potential $u(x)$ and conductivity field $\sigma(x)$,
respectively. Systems of algebraic equations obtained after such
discretization are solved with {\tt UMFPACK}, a solver for
nonsymmetric sparse linear systems \cite{UMFPACK}. The same technique
is used for numerical solutions of adjoint problems \eqref{eq:adjoint}.

\begin{figure}[htb!]
  \begin{center}
  \mbox{
  \subfigure[model~\#1]{\includegraphics[width=0.5\textwidth]{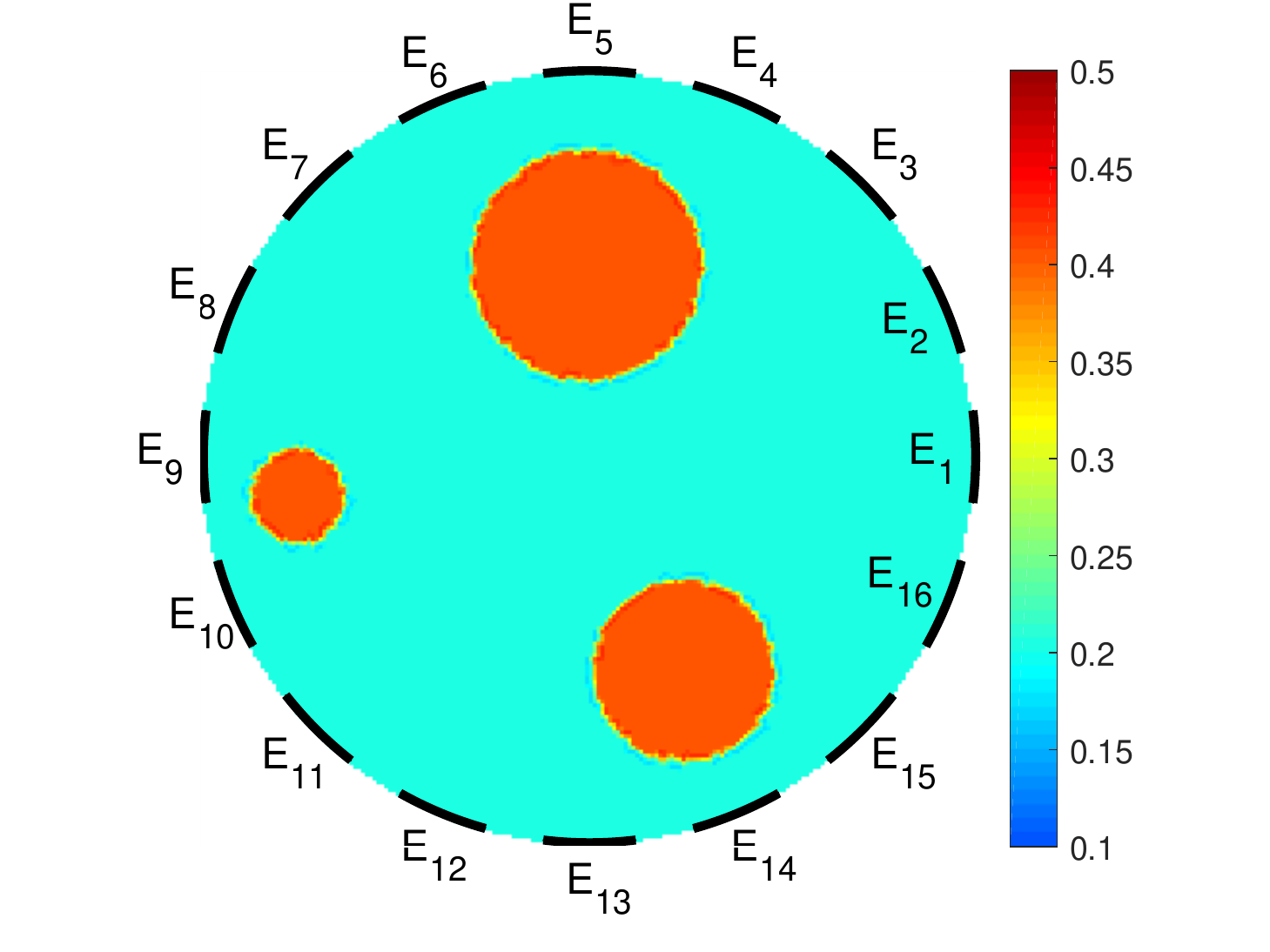}}
  \subfigure[]{\includegraphics[width=0.5\textwidth]{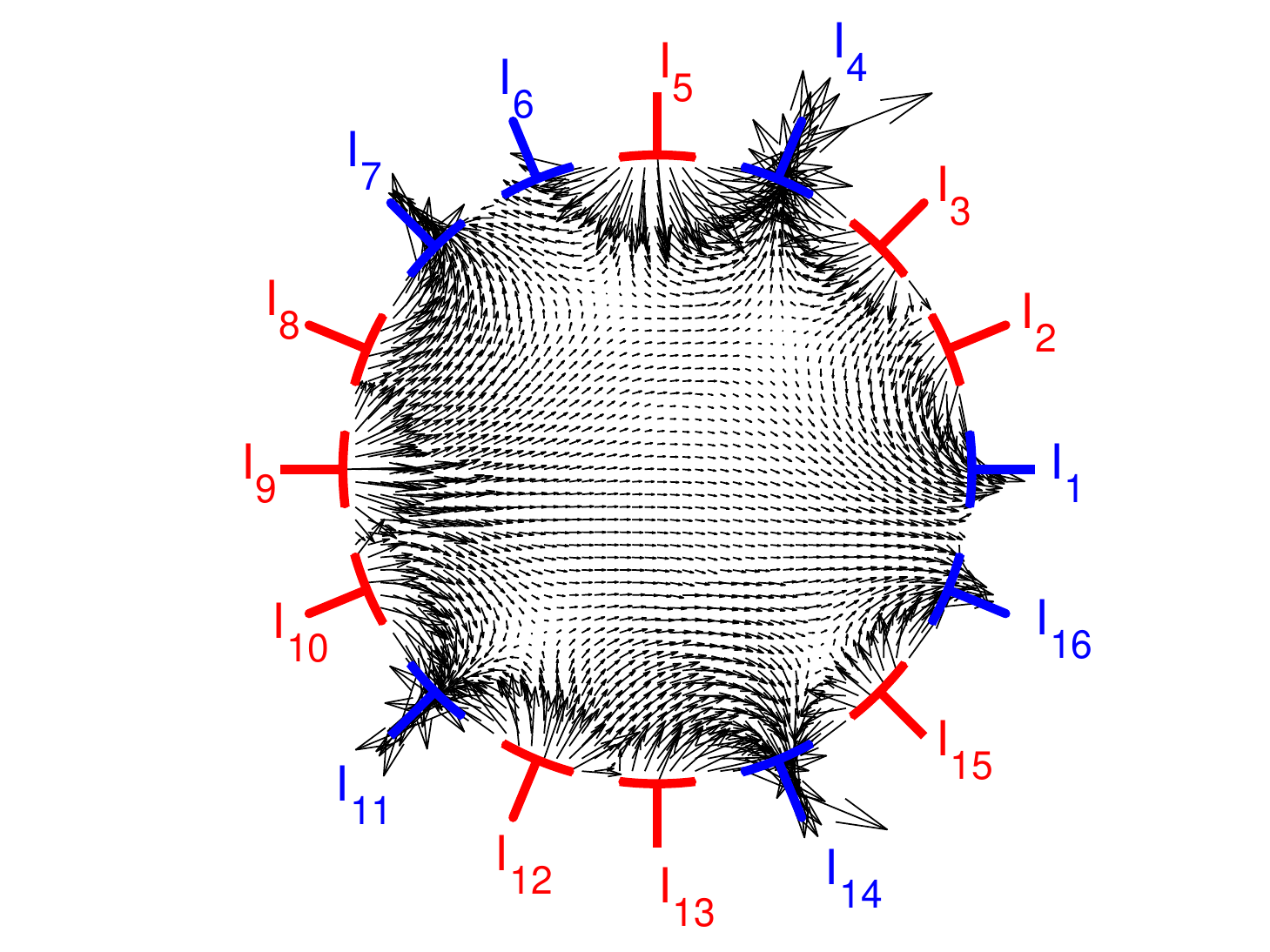}}}
  \end{center}
  \caption{(a)~EIT model~\#1: true electrical conductivity $\sigma_{true}(x)$
    and equispaced geometry of electrodes $E_{\ell}$ placed over boundary
    $\partial \Omega$. (b)~Electrical currents $I_l$ (positive in red, negative
    in blue) induced at electrodes $E_{\ell}$. Black arrows show the distribution
    of flux $\sigma(x) \bnabla u(x)$ of electrical potential $u$ in the interior
    of domain $\Omega$.}
  \label{fig:model}
\end{figure}

All computations are performed using a 2D domain \eqref{eq:domain}
which is a disc of radius $R = 0.1$ with $m = 16$ equidistant electrodes
$E_{\ell}$ with half-width $w = 0.12$~rad covering approximately 61\%
of boundary $\partial \Omega$ as shown in Figure~\ref{fig:model}(a).
Electrical potentials
\begin{equation*}
  (U_{\ell})_{\ell=1}^{16} = \left\{ -3, +1, +2, -5, +4, -1, -3, +2,
  +4, +3, -3, +3, +2, -4, +1, -3 \right\}
  \label{eq:meas_a}
\end{equation*}
are applied to electrodes $E_{\ell}$ following the ``rotation scheme''
discussed in Section~\ref{sec:math_opt} and chosen to be consistent with
the ground potential condition \eqref{eq:ground_conds}. Figure~\ref{fig:model}(b)
shows an example of the distribution of flux $\sigma(x) \bnabla u(x)$ of
electrical potential $u$ in the interior of domain $\Omega$ and measured
currents $(I_{\ell}^*)_{\ell=1}^{16}$ during the EIT procedure. Determining
the Robin part of the boundary conditions in \eqref{eq:forward_3}, we
equally set the electrode contact impedance $Z_{\ell} = 0.1$.

The actual (true) electrical conductivity $\sigma_{true}(x)$ we seek to
reconstruct is defined analytically for each model in \eqref{eq:sigma_true}
by setting $\sigma_c = 0.4$ and $\sigma_h = 0.2$ unless stated otherwise.
The initial guess for control $\mP$ at Step~2 is provided by the
parameterization of initial basis $\mB^0$ obtained after completing Step~1.
For control $\alpha$, the initial values are set to be equal,
i.e., $\alpha^0_i = 1/N_s$. The initial state of control $\zeta$ in
Step~3 is approximated by
\begin{equation}
  \begin{aligned}
    \sigma^0_{th,n} &= \sigma_{ini} = \frac{1}{2}
    \left[ \underset{1 \leq i \leq N_{\Omega}}{\max} \ \sigma_i(\zeta^0) +
    \underset{1 \leq i \leq N_{\Omega}}{\min} \ \sigma_i(\zeta^0) \right],\\
    \sigma^0_{low} &= \underset{1 \leq i \leq N_{\Omega}}{\mean}
    \left\{ \sigma_i(\zeta^0) : \ \sigma_i(\zeta^0) < \sigma_{ini} \right\},\\
    \sigma^0_{high,n} &= \underset{1 \leq i \leq N_{\Omega}}{\mean}
    \left\{ \sigma_i(\zeta^0) : P_{i,n+1} = 1, \ \sigma_i(\zeta^0) \geq
    \sigma_{ini} \right\}, \\
    n &= 1, \ldots, N_{\max}.
  \end{aligned}
  \label{eq:minJ_zeta_ini_1}
\end{equation}
Termination criteria are set by tolerance $\epsilon = 10^{-9}$ in
\eqref{eq:termination} and the total number of cost functional
evaluations of 50,000, whichever is reached first.

\begin{figure}[!htb]
  \begin{center}
  \mbox{
  \subfigure{\includegraphics[width=0.25\textwidth]{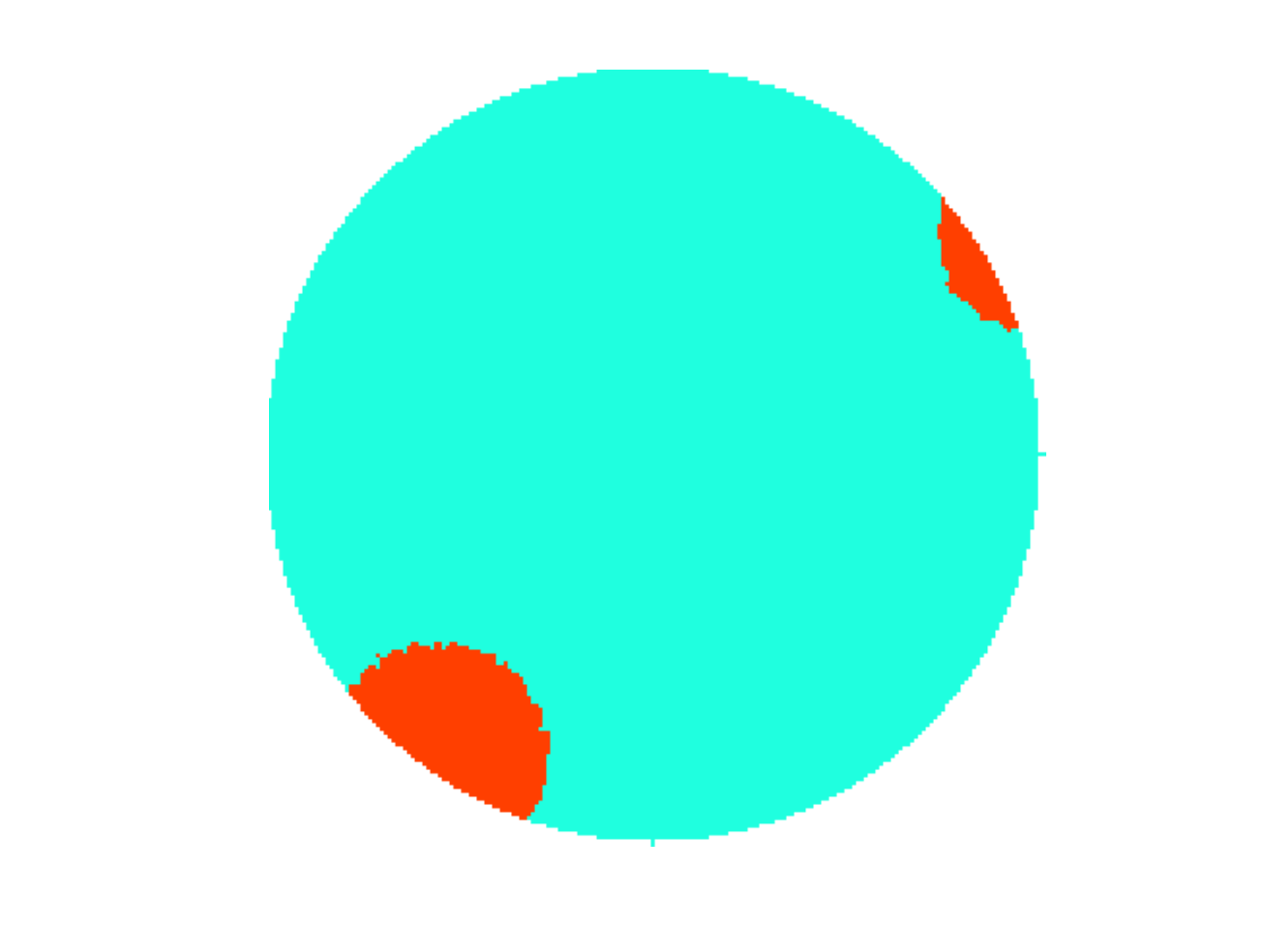}}
  \subfigure{\includegraphics[width=0.25\textwidth]{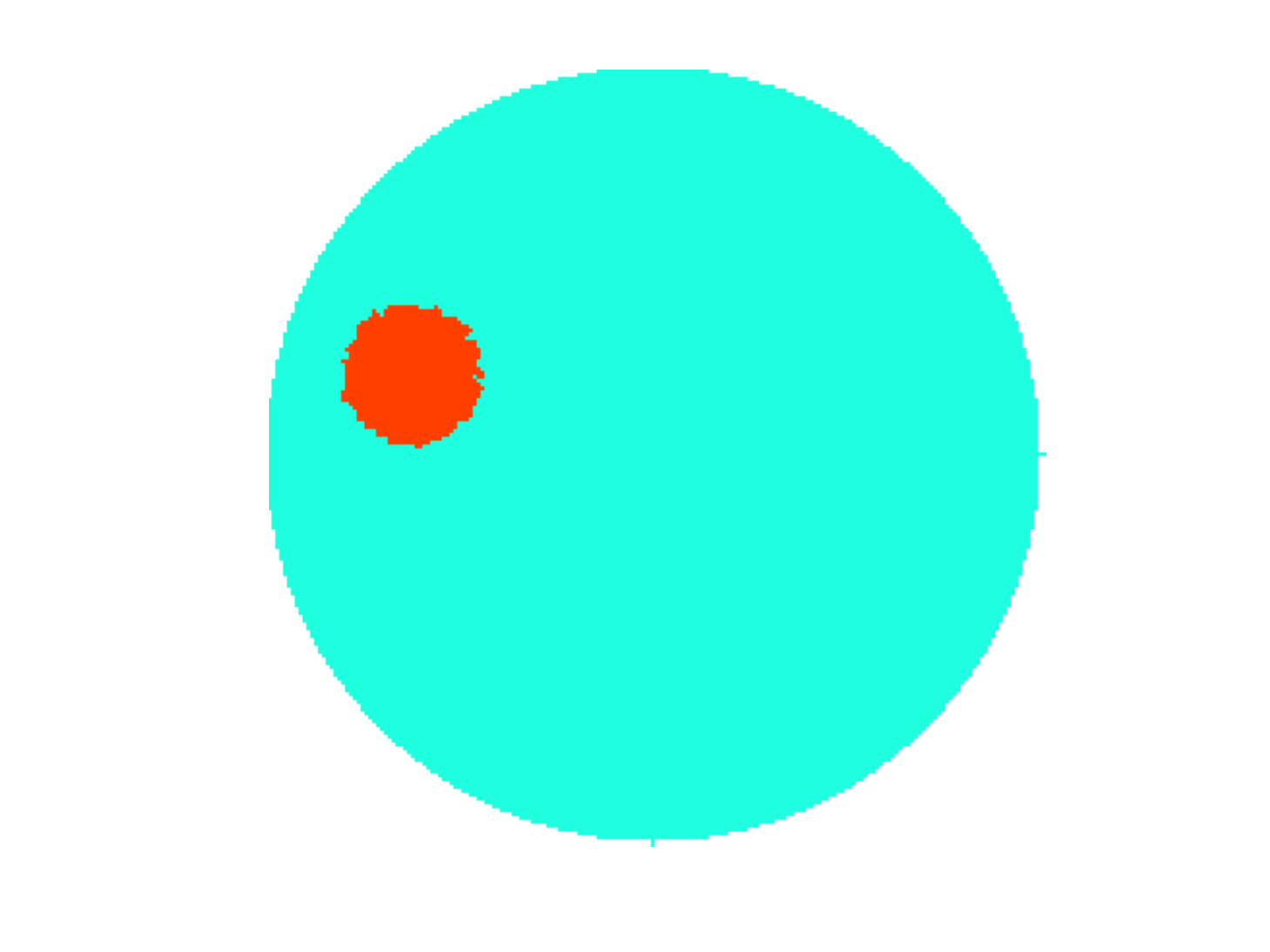}}
  \subfigure{\includegraphics[width=0.25\textwidth]{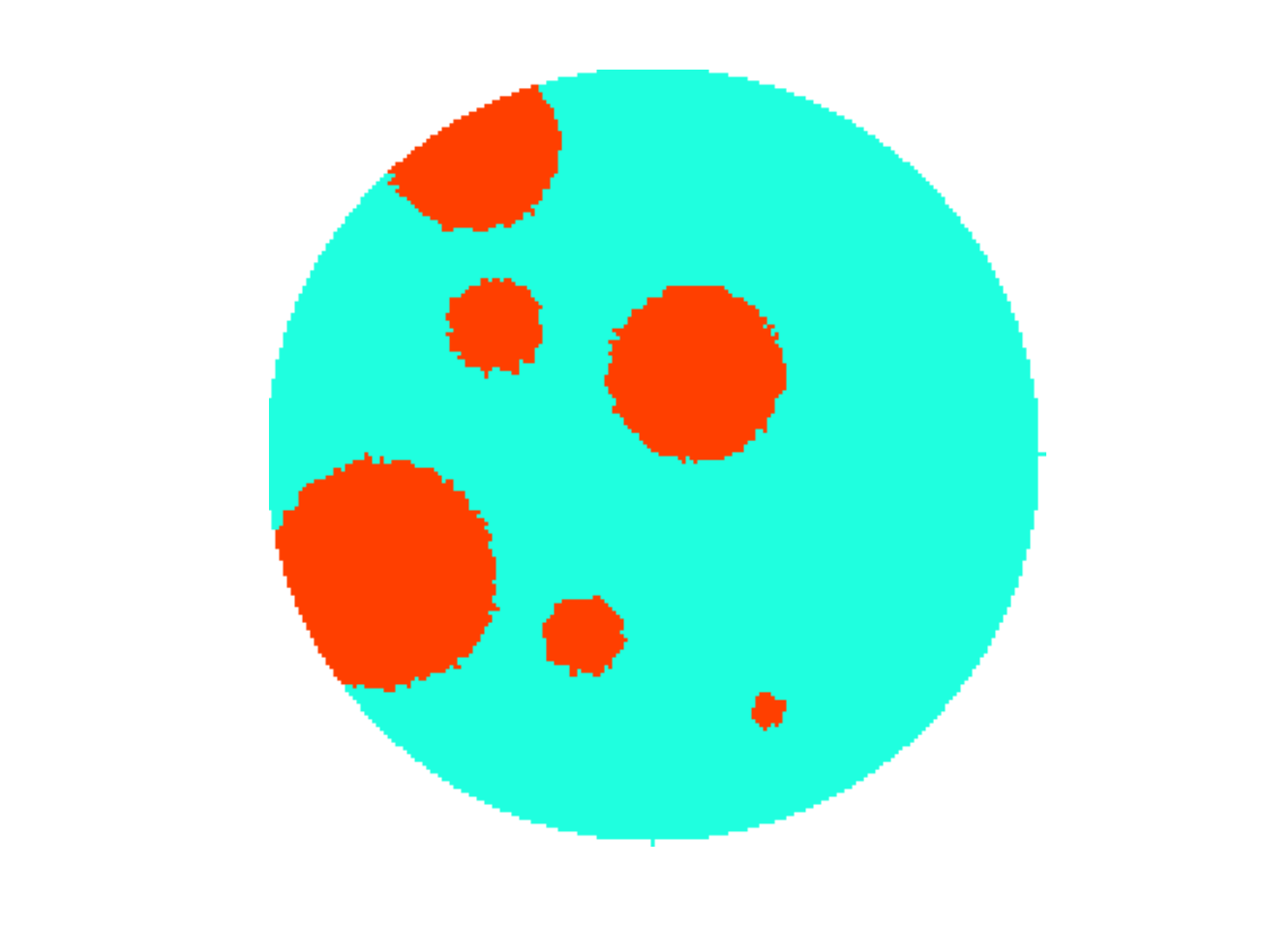}}
  \subfigure{\includegraphics[width=0.25\textwidth]{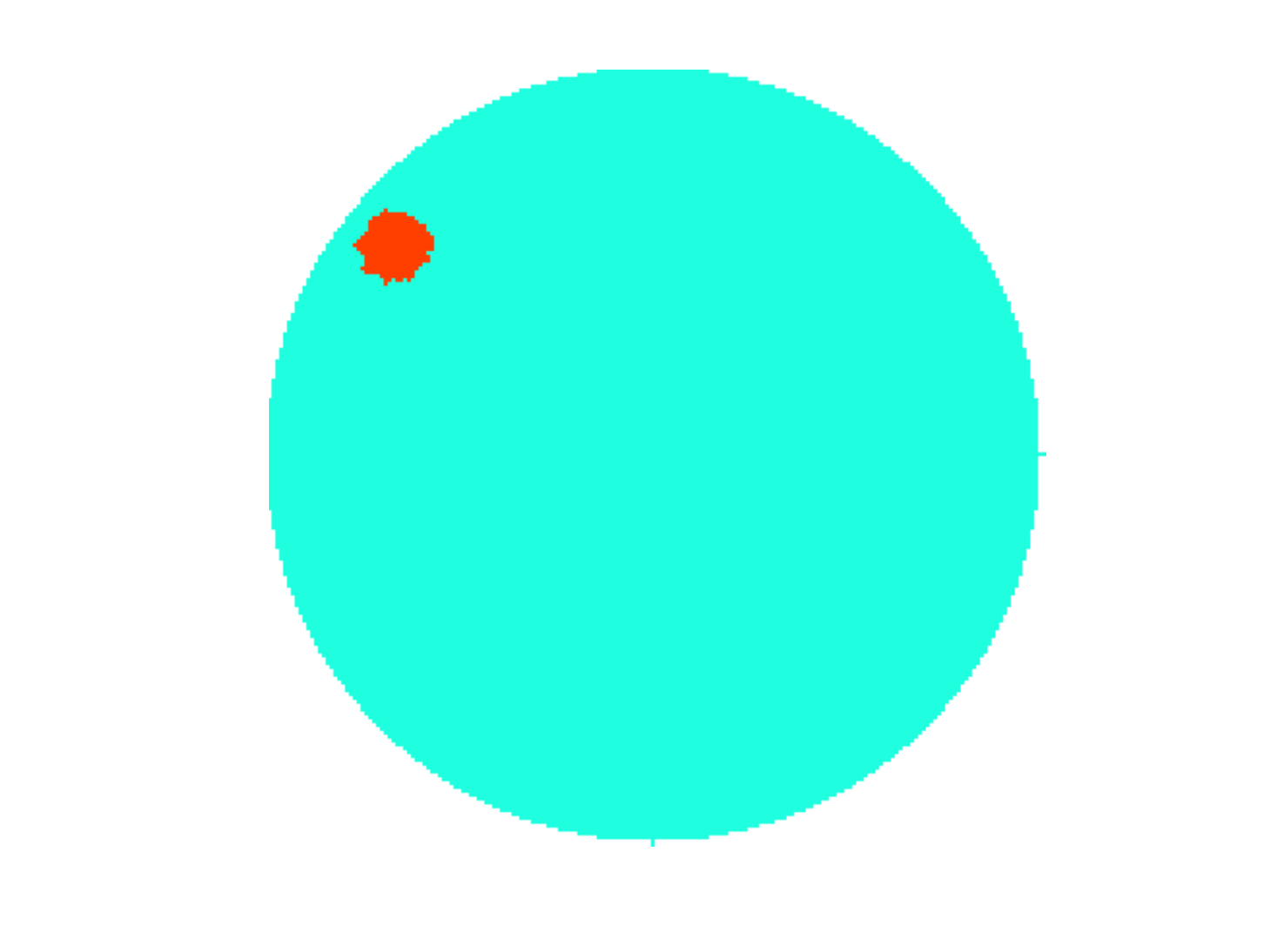}}}
  \mbox{
  \subfigure{\includegraphics[width=0.25\textwidth]{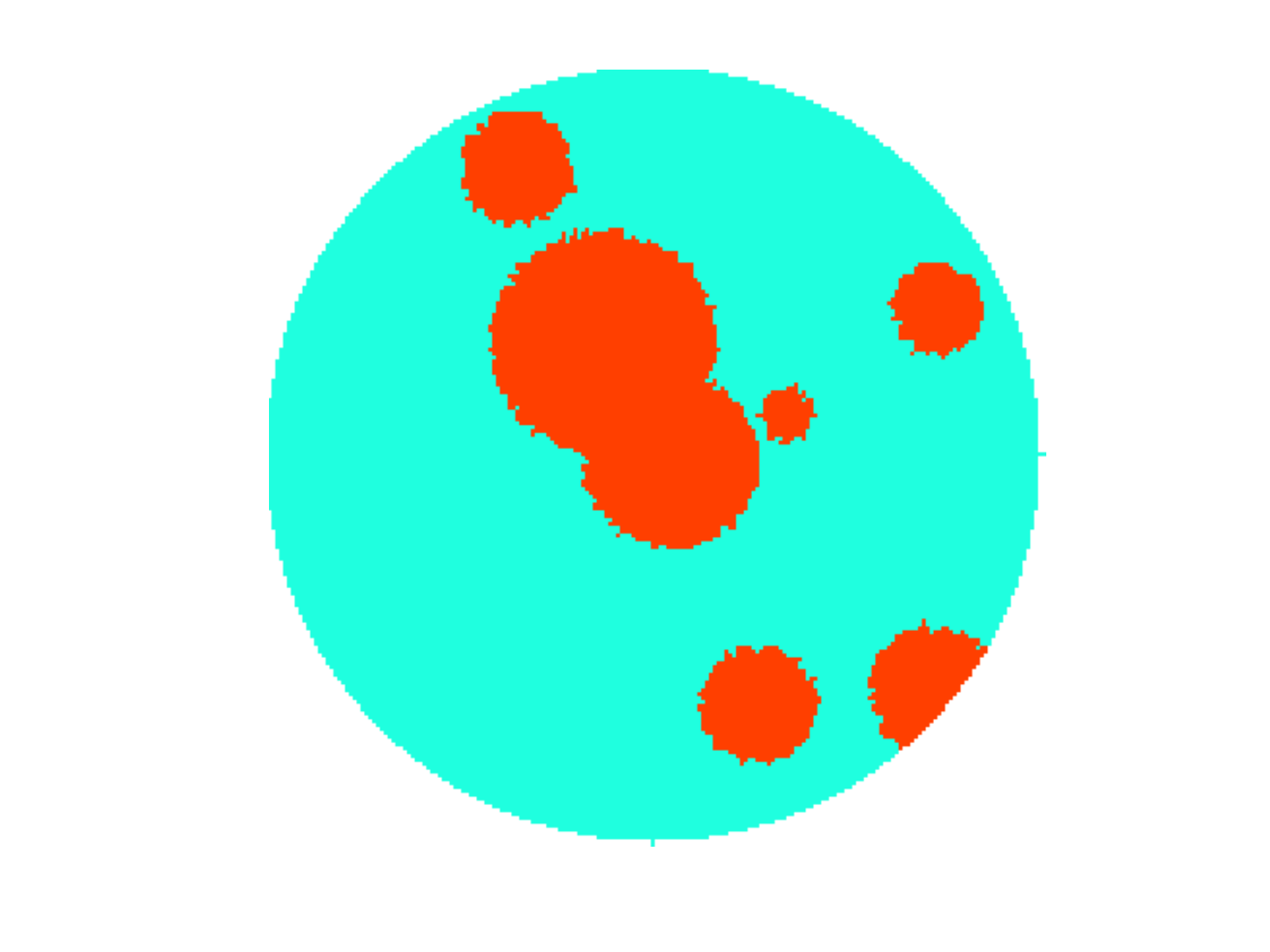}}
  \subfigure{\includegraphics[width=0.25\textwidth]{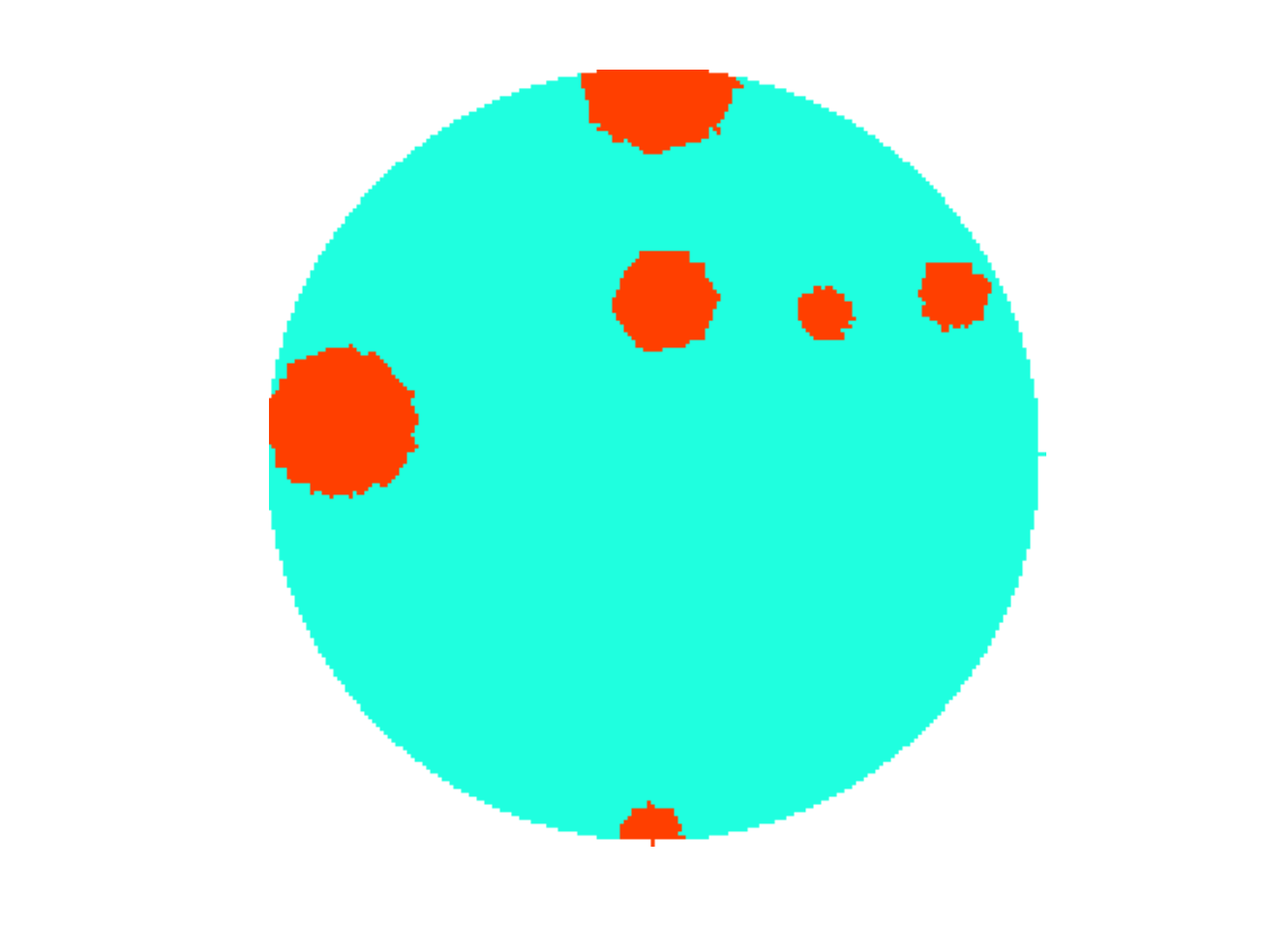}}
  \subfigure{\includegraphics[width=0.25\textwidth]{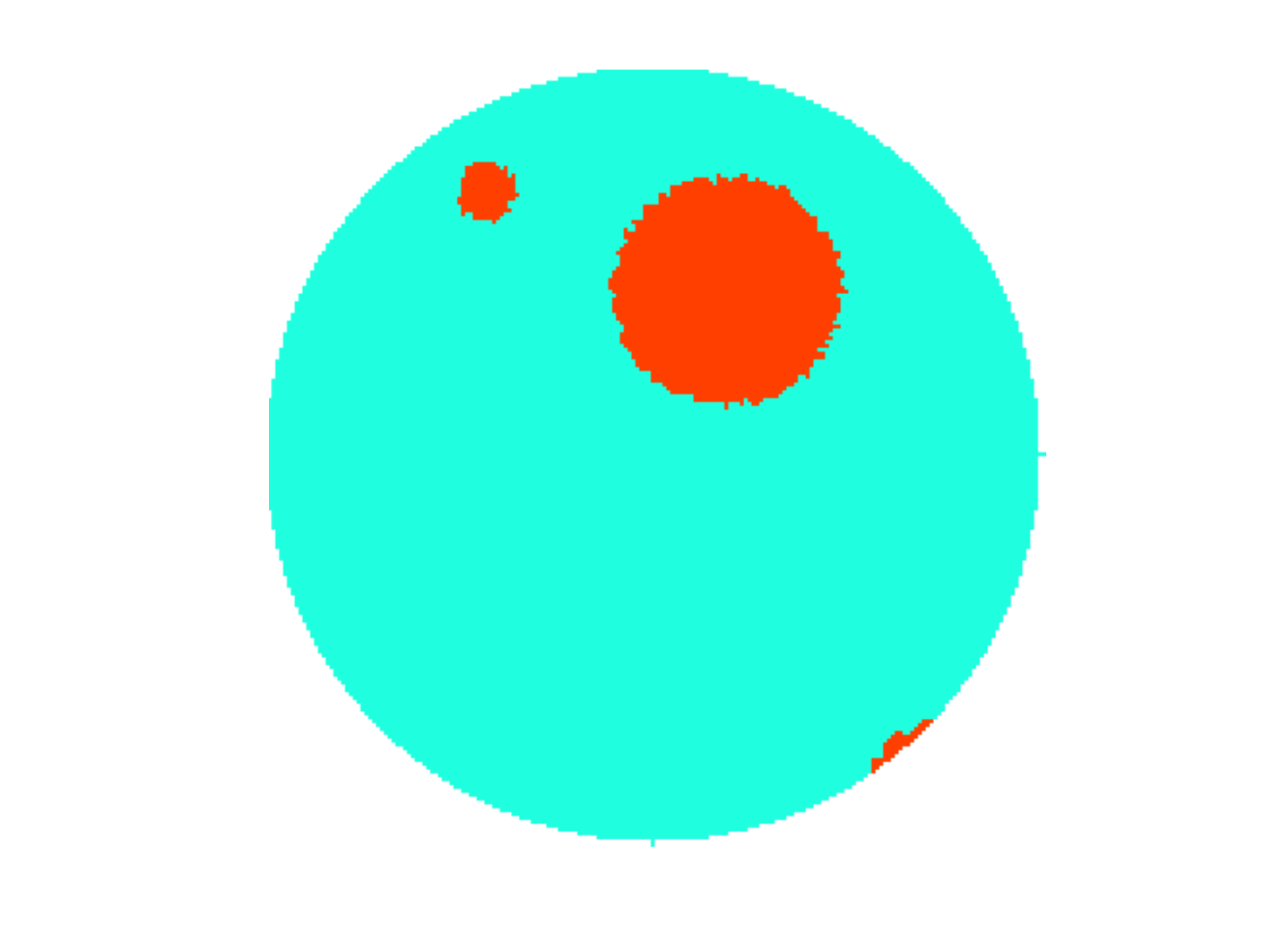}}
  \subfigure{\includegraphics[width=0.25\textwidth]{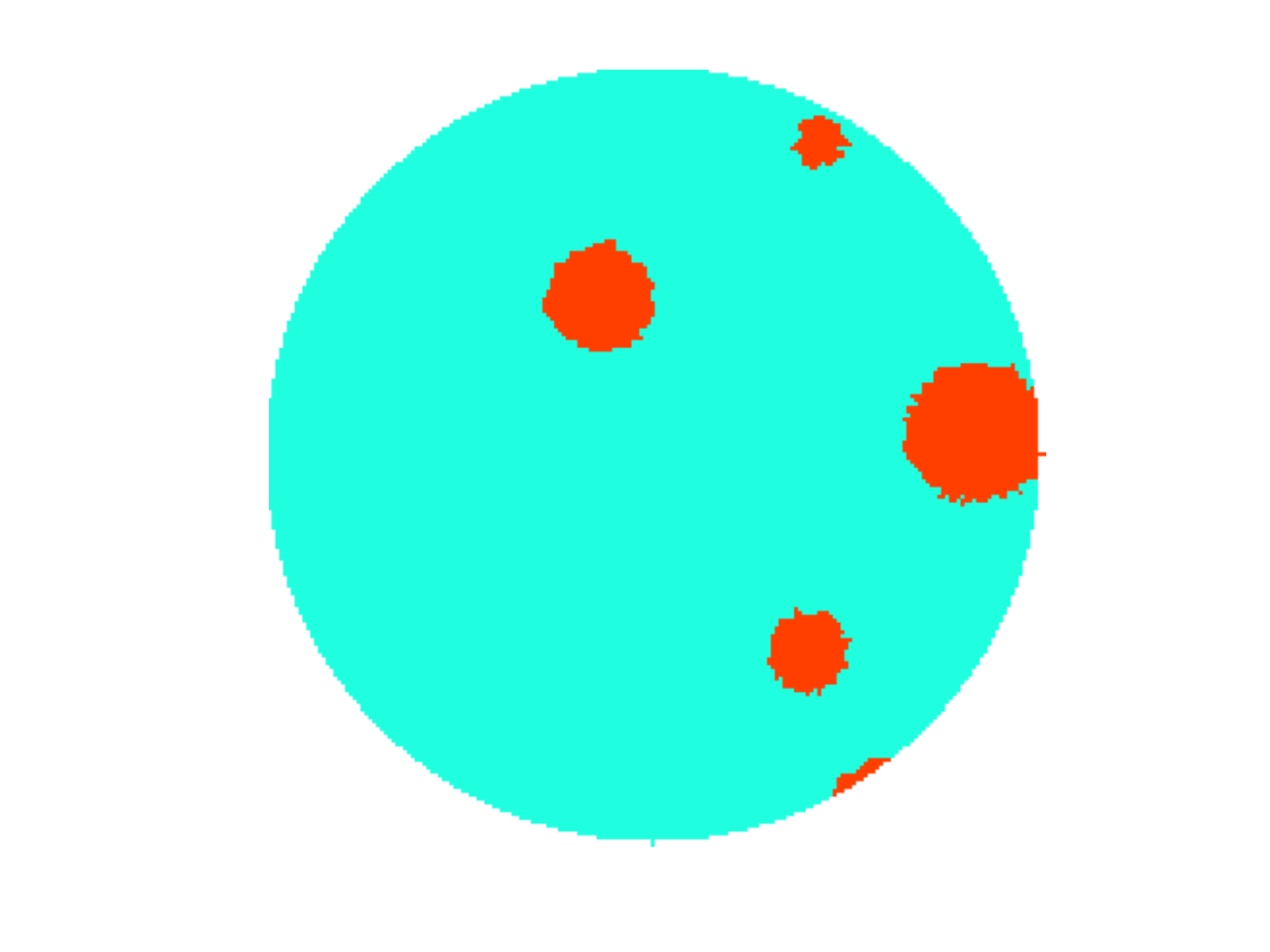}}}
  \end{center}
  \caption{8 first sample solutions $\bar \sigma_i(x), \, i = 1, \ldots, 8$,
     from the set $\mC(10,000)$.}
  \label{fig:samples}
\end{figure}

For generating samples in $\mC(N)$ collections discussed in
Section~\ref{sec:step_0}, we use $N = 10,000$ and $N_{c,max} = 8$. This set
is precomputed using a generator of uniformly distributed random numbers.
Therefore, each sample $\bar \sigma_i(x)$ ``contains'' from one to eight
``cancer-affected'' areas with $\sigma_c = 0.4$. Each area is located randomly
within domain $\Omega$ and represented by a circle of randomly chosen radius
$0 < r \leq 0.3 R$ as exemplified in Figure~\ref{fig:samples}. Also, we fix
the number of samples $N_s$ to 10 for all
numerical experiments shown in this paper. Finally, the results of our previous
research \cite{ArbicBukshtynov2022,ArbicMS2020} confirm that our sample-based
parameterization enables high-level stability of the obtained results towards
the noise present in measurements. Thus, for all numerical experiments shown
in this paper, we use measurement data contaminated with $0.5\%$ normally
distributed (Gaussian) noise.

\subsection{Optimization Framework Validation}
\label{sec:frame_valid}

To demonstrate the applicability of the proposed computational framework discussed
in Section~\ref{sec:solution} and check its overall performance through a sequence
of steps while solving the inverse EIT problem, we use model~\#1 featuring three
circular-shaped cancerous regions of various sizes, as shown in
Figure~\ref{fig:model}(a). We use this model to mimic a typical situation for
a cancer-affected biological tissue containing several spots suspected of being
tumorous and, as such, having elevated electrical conductivity.

\subsubsection{Validating Gradients}
\label{sec:grads}

\begin{figure}[htb!]
  \begin{center}
  \mbox{
  \subfigure[]{\includegraphics[width=0.5\textwidth]{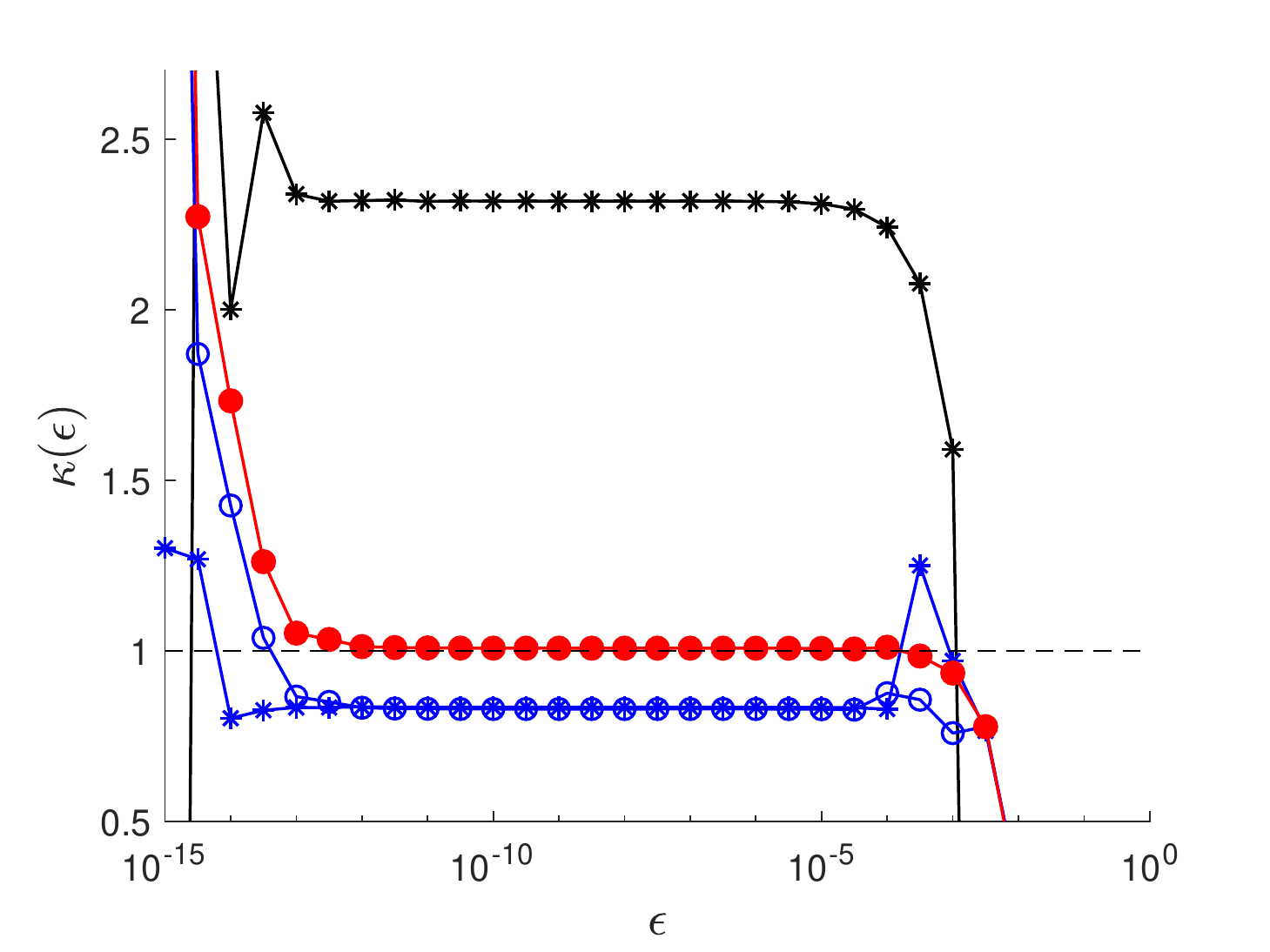}}
  \subfigure[]{\includegraphics[width=0.5\textwidth]{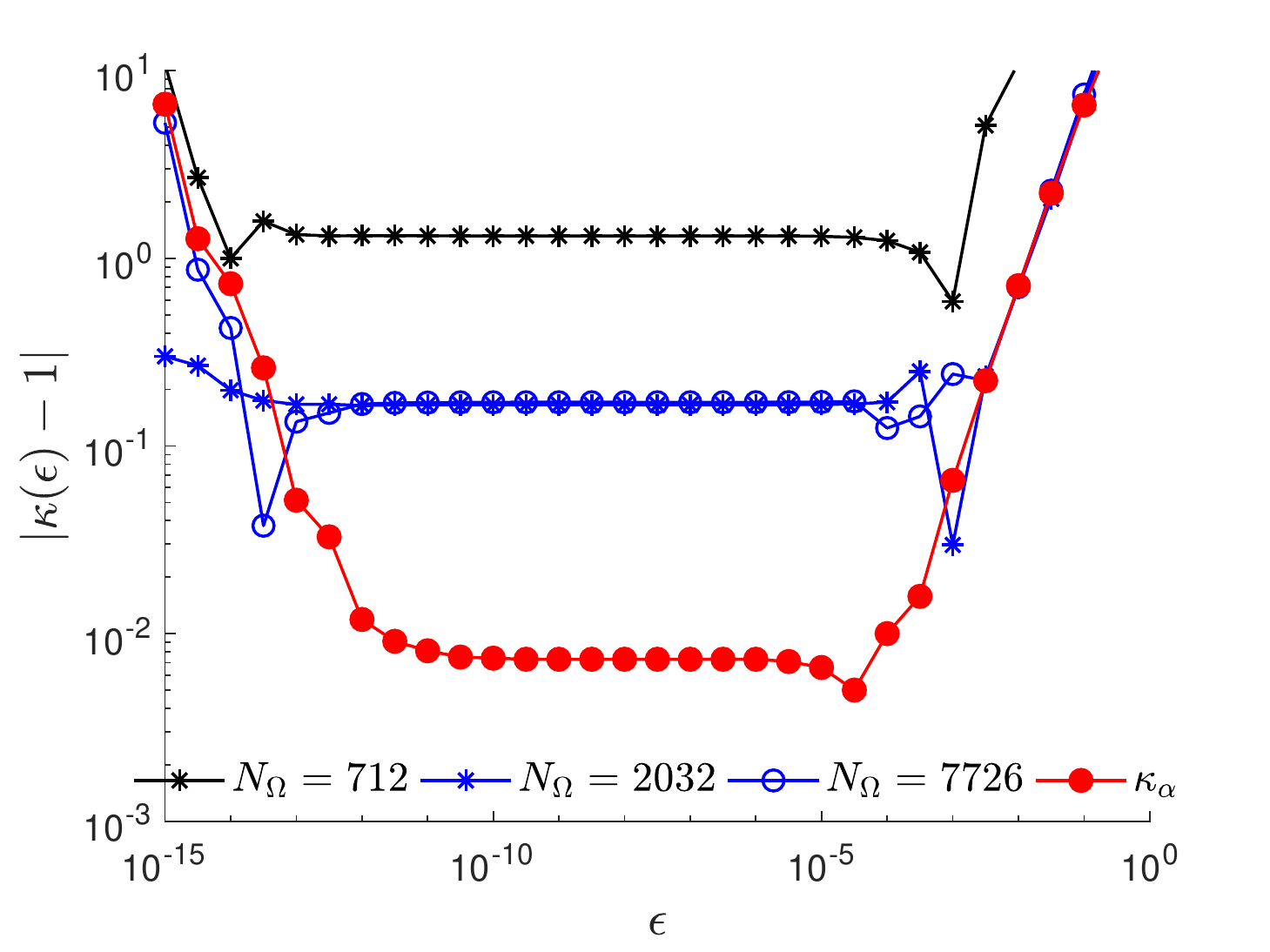}}}
  \end{center}
  \caption{The behavior of (a)~$\kappa (\epsilon)$ and (b)~$\log_{10} \lvert
    \kappa(\epsilon) - 1 \rvert$ as functions of $\epsilon$ while checking the
    consistency of used gradients computed for model~\#1.
    In both graphs, red circles show results for $\kappa_{\alpha}(\epsilon)$
    and $N_{\Omega} = 7,726$, while the rest relates to the application of the
    $\kappa$-test to both controls $\alpha$ and $\mP$ in domain $\Omega$ with
    different spatial discretizations.}
  \label{fig:kappa_all}
\end{figure}

\begin{figure}[htb!]
  \begin{center}
  \mbox{
  \subfigure[]{\includegraphics[width=0.5\textwidth]{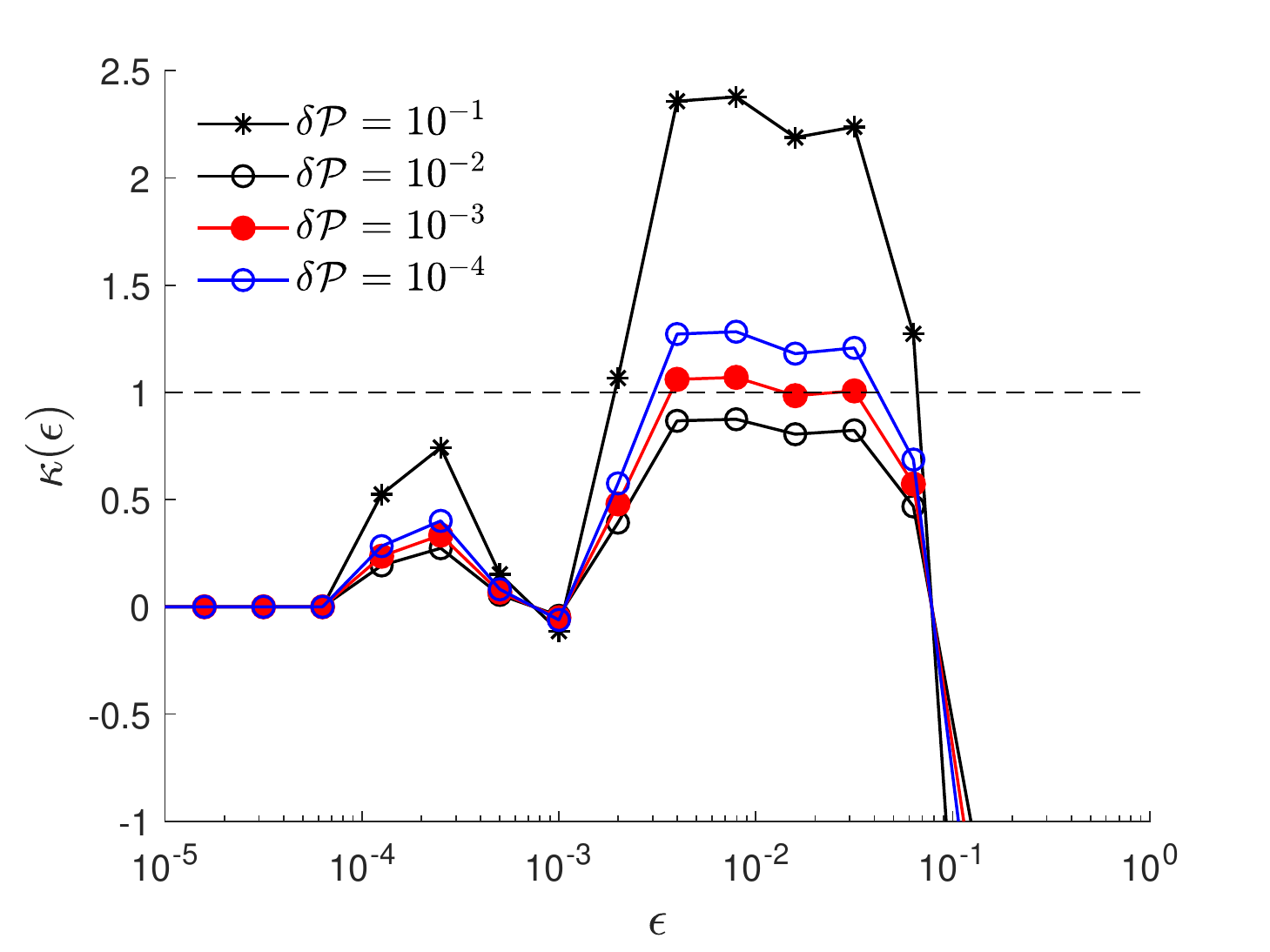}}
  \subfigure[]{\includegraphics[width=0.5\textwidth]{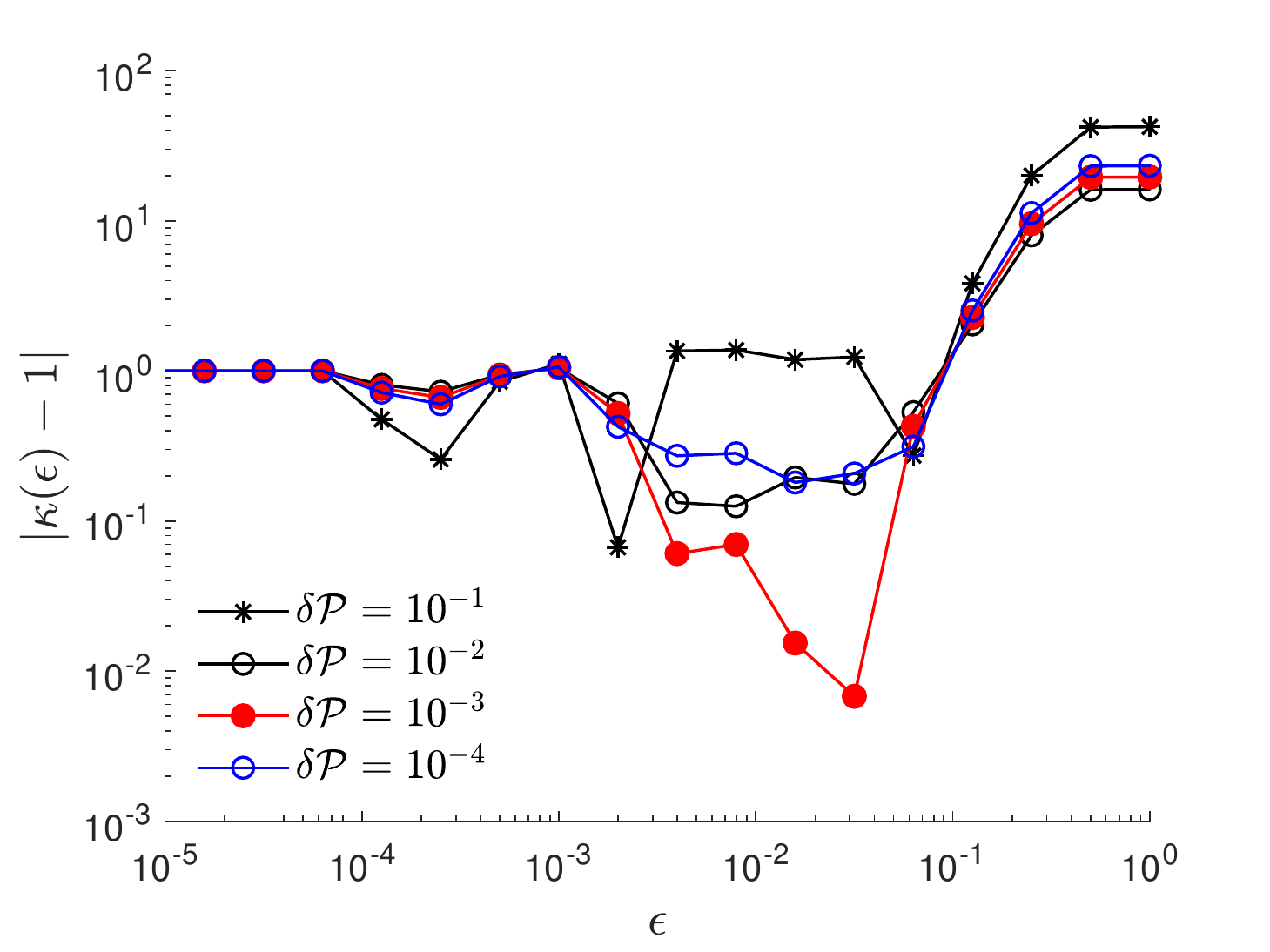}}}
  \end{center}
  \caption{The behavior of (a)~$\kappa_{\mP} = \kappa (\epsilon)$ and (b)~$\log_{10}
    \lvert \kappa(\epsilon) - 1 \rvert$ as functions of $\epsilon$ while checking the
    consistency of gradients $\bnabla_{\mP} \mJ$ computed for model~\#1 with different
    values of parameter $\delta \mP$.}
  \label{fig:kappa_P}
\end{figure}

The gradient-based concept is central to the proposed methodology; therefore, our
first results demonstrate the consistency of the gradients computed for various
controls (optimization variables), as discussed in Section~\ref{sec:step_2}.
Figures~\ref{fig:kappa_all} and \ref{fig:kappa_P} show the results of a diagnostic
test ($\kappa$-test) commonly employed to verify the correctness of the discretized
gradients; see, e.g.,
\cite{Bukshtynov2011,Bukshtynov2013,ChunEdwardsBukshtynov2024,BukshtynovBook2023}.
It consists in computing the directional differential, e.g., $\mJ'(\alpha;
\delta \alpha) = \left\langle \bnabla_{\alpha} \mJ, \delta \alpha \right\rangle_{L_2}$,
for some selected variations (perturbations) $\delta \alpha$ in two different ways:
namely, using a finite--difference approximation versus using \eqref{eq:grad_alpha}
with \eqref{eq:grad_sigma} and \eqref{eq:grad_sigma_spat}, and then examining the ratio
of the two quantities, i.e.,
\begin{equation}
  \kappa_{\alpha} (\epsilon) = \dfrac{\frac{1}{\epsilon} \left[ \mJ(\alpha + \epsilon \,
  \delta \alpha) - \mJ(\alpha) \right]}
  {\int_{\Omega} \bnabla_{\alpha} \mJ \delta \alpha \, dx}
  \label{eq:kappa}
\end{equation}
for a range of values of $\epsilon$. If these gradients are computed correctly
then for intermediate values of $\epsilon$, $\kappa_{\alpha}(\epsilon)$ will be
close to the unity.

Figure~\ref{fig:kappa_all}(a) demonstrates such behavior over a range of $\epsilon$
spanning about 10 orders of magnitude for $\bnabla_{\alpha} \mJ$ (in red) discretized
in the same manner as domain $\Omega$, e.g., with $N_{\Omega} = 7,726$. As can be
expected, the quantity $\kappa(\epsilon)$ deviates from the unity for very small values
of $\epsilon$ due to the subtractive cancelation (round--off) errors and also for large
values of $\epsilon$ due to the truncation errors (both of which are well--known effects).
In addition, the quantity $\log_{10} \lvert \kappa(\epsilon) - 1 \rvert$ plotted in
Figure~\ref{fig:kappa_all}(b) shows how many significant digits of accuracy are captured
in a given gradient evaluation.

Similarly, we apply the same $\kappa$-test to check the correctness of all steps involved
in gradient computations associated with control $\mP$, as given by \eqref{eq:grad_P} with
\eqref{eq:grad_P_i} and \eqref{eq:grad_sigma_spat}. However, as discussed in
Section~\ref{sec:step_2}, computing $\bnabla_{\mP} \mJ$ involves FD estimations of
$\bnabla_{\mP} \sigma$ whose accuracy depends on the choice of perturbations $\Delta \mP_i$.
Figure~\ref{fig:kappa_P}(a) shows the obtained results after applying the $\kappa$-test
when $\Delta \mP_i$ are constant vectors with all components equal to $\delta \mP$ ranging
between $10^{-4}$ and $10^{-1}$. The expected ``unity plateau'' forms, confirming the
correctness of all computations; however, its length is limited to about two orders of
$\epsilon$ ($10^{-3} \div 10^{-1}$) values. In the same fashion, Figure~\ref{fig:kappa_P}(b)
depicts significant digits of accuracy. We use the results of this test to ``calibrate'' our
computational framework by setting $\delta \mP$ to $10^{-3}$ for all experiments.

Finally, we explore the results obtained by applying the $\kappa$-test to both controls
(simultaneously to $\alpha$ and $\mP$) to determine a proper discretization of domain
$\Omega$ as an obviously better approximation of continuous gradients by finer spatial
discretizations here is affected by the added FD approximations in $\bnabla_{\mP} \mJ$
parts. Although the results provided in Figures~\ref{fig:kappa_all}(a) and
\ref{fig:kappa_all}(b) do not show a significant difference in using $N_{\Omega} = 2,032$
and $N_{\Omega} = 7,726$, we will use the latter as the number of FEM elements for all
numerical experiments and all models in this paper to ensure better resolution for the
reconstructed images.

\subsubsection{Model~\#1: Validating Performance}
\label{sec:model_1}

\begin{figure}[htb!]
  \begin{center}
  \mbox{
  \subfigure[]{\includegraphics[width=0.5\textwidth]{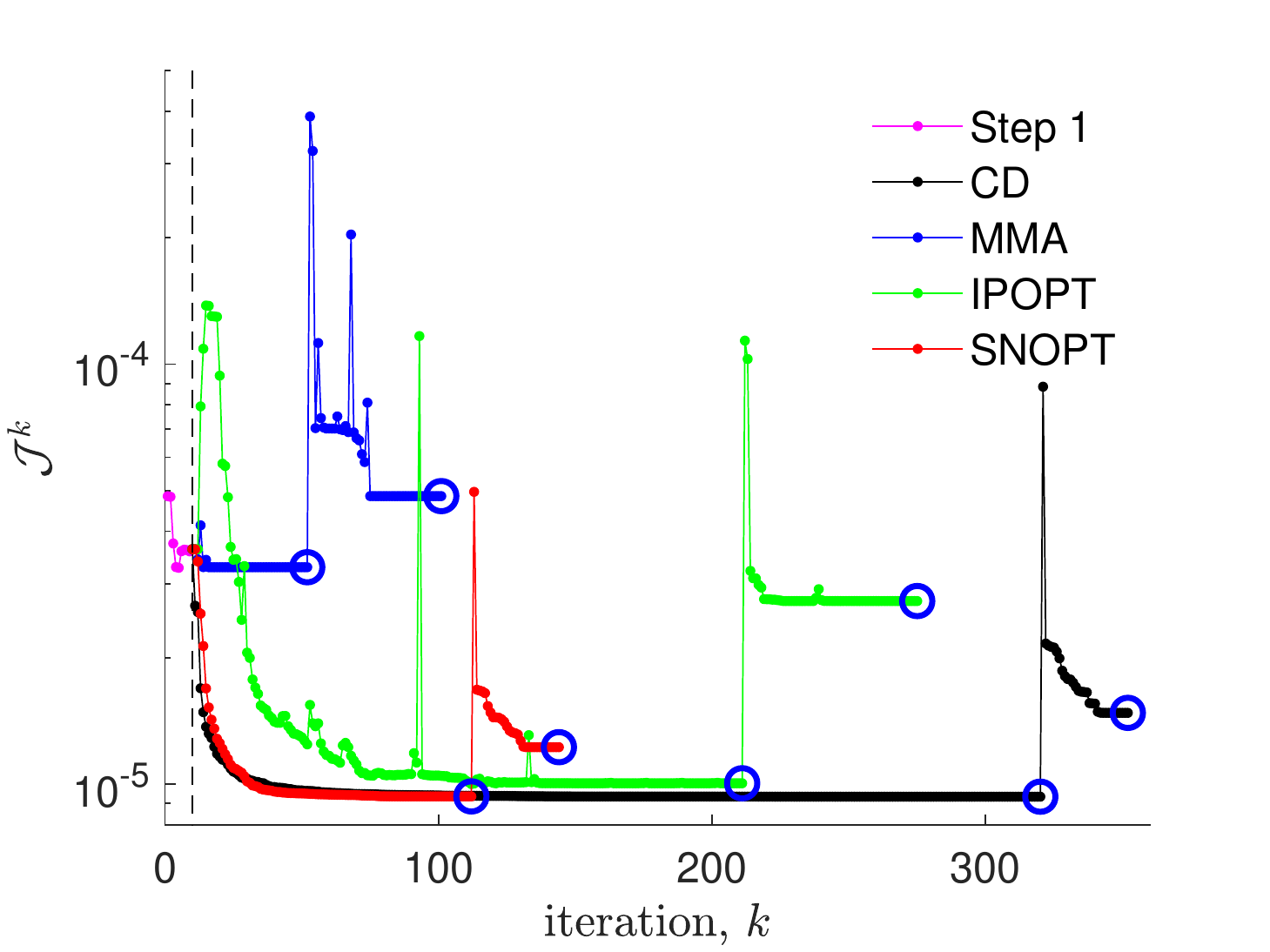}}
  \subfigure[]{\includegraphics[width=0.5\textwidth]{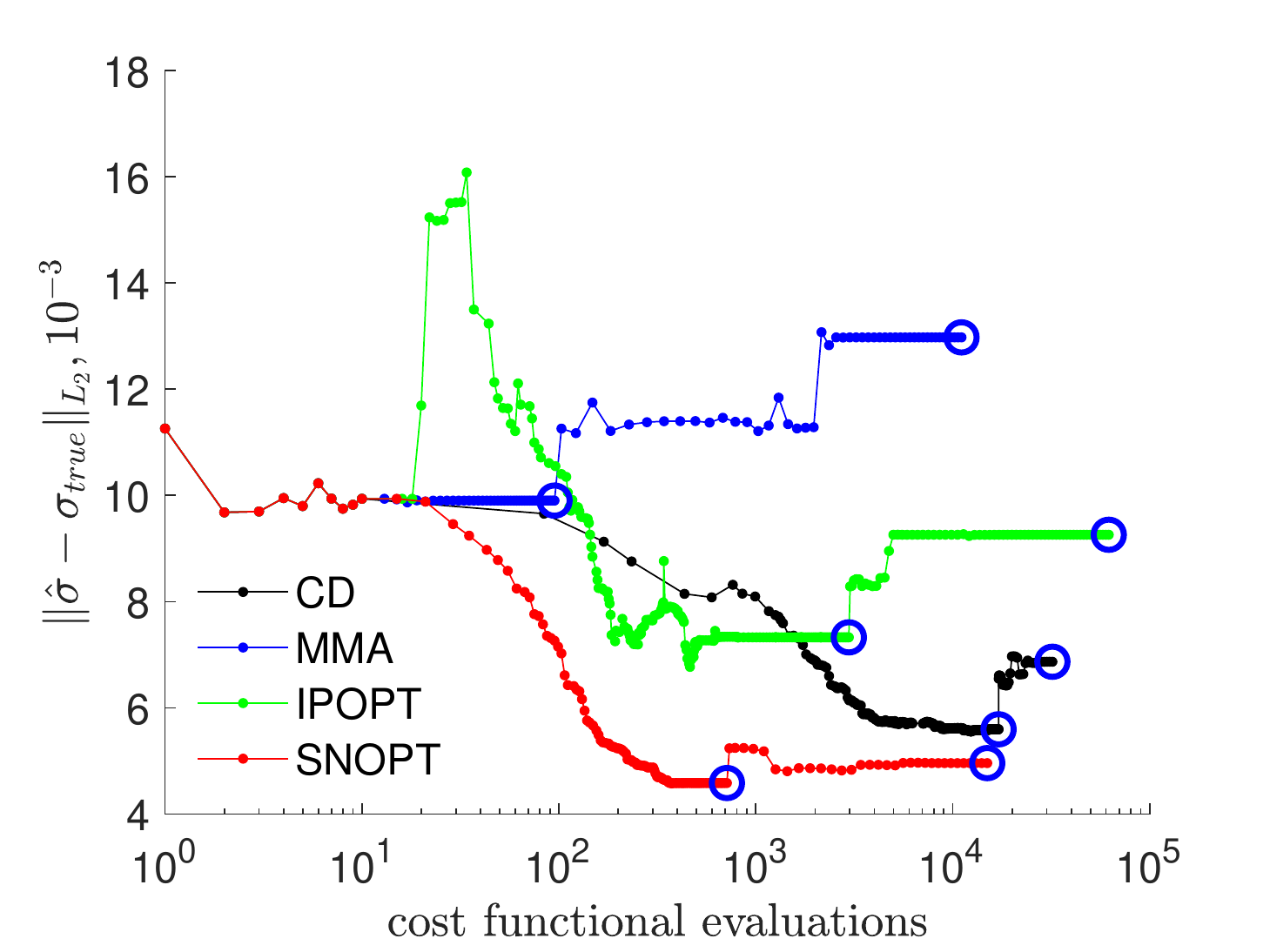}}}
  \end{center}
  \caption{Optimization results for model~\#1: (a)~cost functionals $\mJ^k$
    as functions of iteration count $k$ and (b)~solution errors $\| \sigma^k -
    \sigma_{true} \|_{L_2}$ as functions of a number of cost functional evaluations
    evaluated while employing different optimizers (CD, {\tt MMA}, {\tt IPOPT},
    and {\tt SNOPT}). Pink dots in (a) represent Step~1 solutions ($k = 1, 2, \ldots, 10$),
    and empty blue circles in (a) and (b) show solutions obtained after Step~2 and 3 phases
    are complete.}
  \label{fig:model_8_a}
\end{figure}

In this section, we will prove the superior performance of the proposed computational
framework supplied with the fine-scale gradient-based and coarse-scale binary tuning
optimization algorithms, as detailed in Section~\ref{sec:solution}. First, we refer
to Figure~\ref{fig:model_8_a}; it shows the results of applying this methodology to
our model~\#1 and a comparison of the computational performance observed while employing
different optimizers, namely
\begin{itemize}
  \item gradient-based sequential quadratic, nonlinear interior point, and sequential
    convex programming algorithms by respectively {\tt SNOPT} \cite{SNOPTManual},
    {\tt IPOPT} \cite{IPOPTManual}, and {\tt MMA} \cite{Svanberg1987},
  \item previously used in \cite{ArbicBukshtynov2022,ArbicMS2020} derivative-free
    CD method customized to a predefined order of controls.
\end{itemize}

Figure~\ref{fig:model_8_a}(a) clearly distinguishes gradient-based {\tt SNOPT} from other
methods by the best results provided by the cost functional $\mJ^k$ evaluated
after all optimization stages (including fine-scale optimization and binary tuning
at the coarse scale) are complete. We expected a comparable quality of the reconstructed
images obtained by {\tt SNOPT} and CD methods as they use the same core concept of
utilizing sample solutions with customized geometry. On top of that, our new approach
uses multiscale control space reduction for binary tuning, also paired with
gradient-based techniques. It makes the gradient-based implementation (by {\tt SNOPT})
more advantageous in both the accuracy of the final solutions and the computational speed,
as demonstrated in Figure~\ref{fig:model_8_a}(b), by comparing solution errors $\| \hat
\sigma - \sigma_{true} \|_{L_2}$ as functions of a total number of cost functional
evaluations (including cases of evaluating cost functionals $\mJ$ to complete computations
for constructing gradients, choosing optimal step size in the gradient-based methods, etc.).
We see this measure to examine the overall performance of various approaches to be
reasonable, as cost functional evaluations paired with numerical solutions for forward
EIT problem \eqref{eq:forward} contribute to the major part of the computational load
(all other parts, in fact, take fractions of a second of the CPU time to be completed).

\begin{figure}[htb!]
  \begin{center}
  \mbox{
  \subfigure[MMA]{\includegraphics[width=0.25\textwidth]{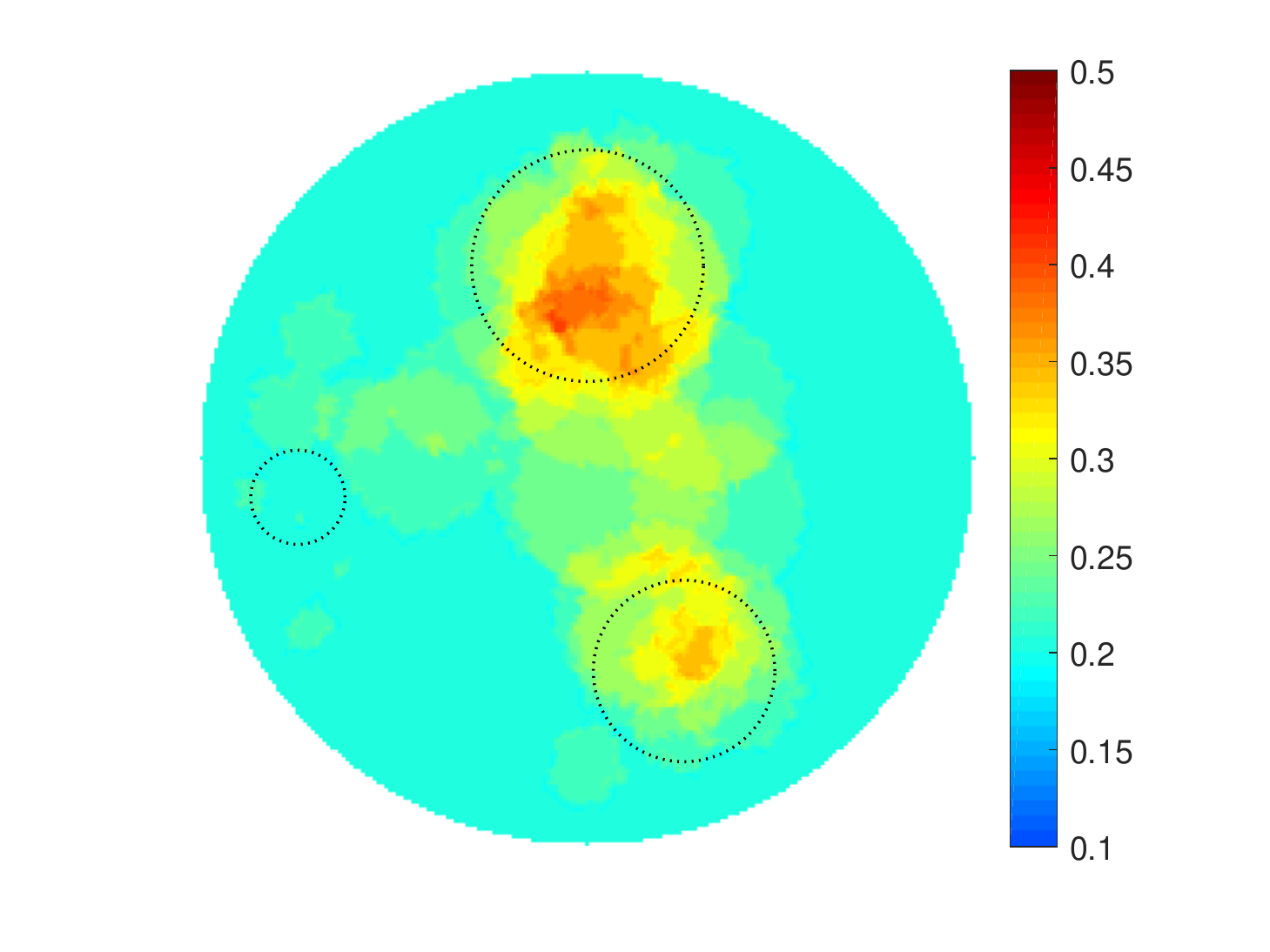}}
  \subfigure[IPOPT]{\includegraphics[width=0.25\textwidth]{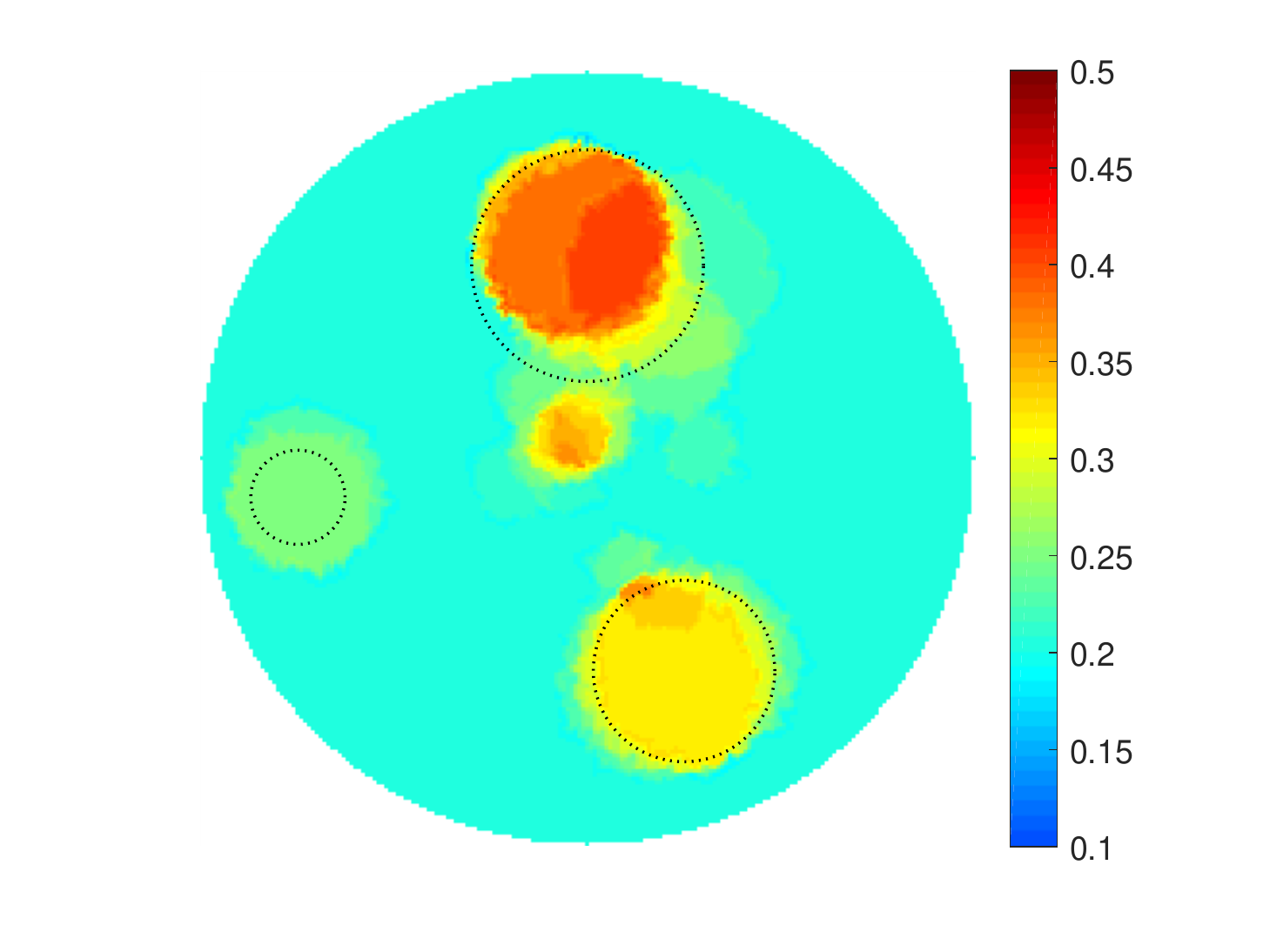}}
  \subfigure[SNOPT]{\includegraphics[width=0.25\textwidth]{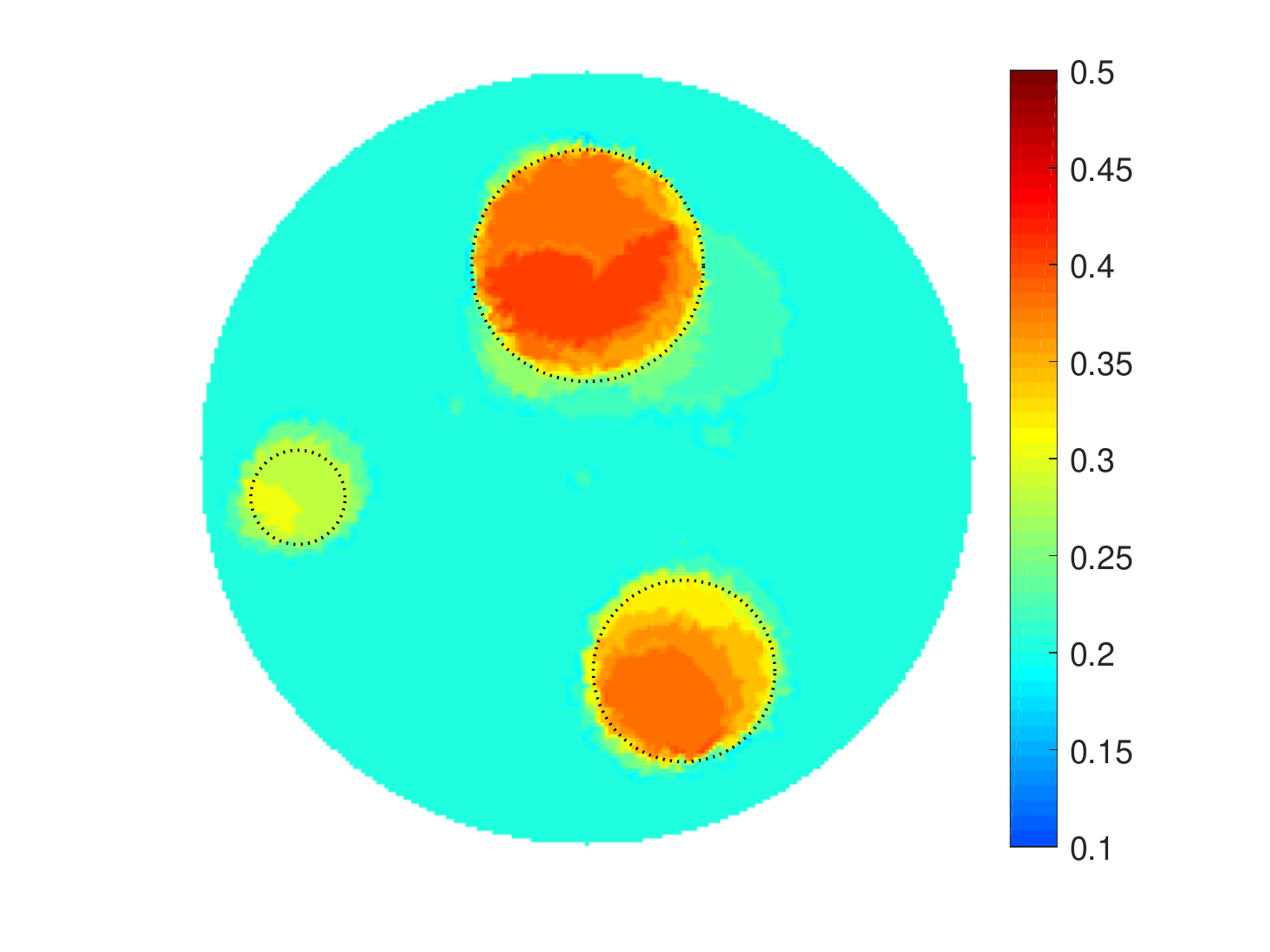}}
  \subfigure[CD]{\includegraphics[width=0.25\textwidth]{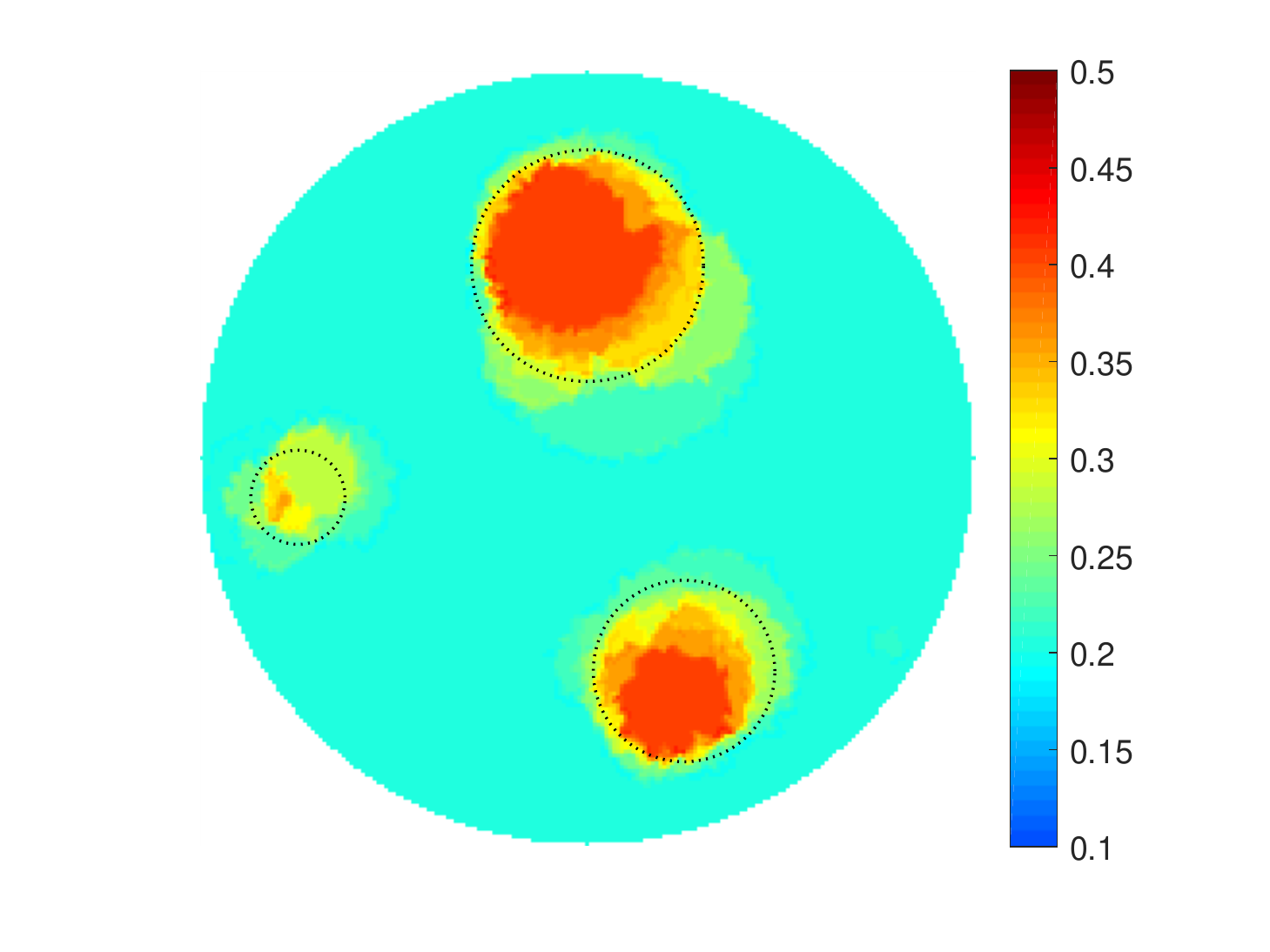}}}
  \mbox{
  \subfigure[MMA]{\includegraphics[width=0.25\textwidth]{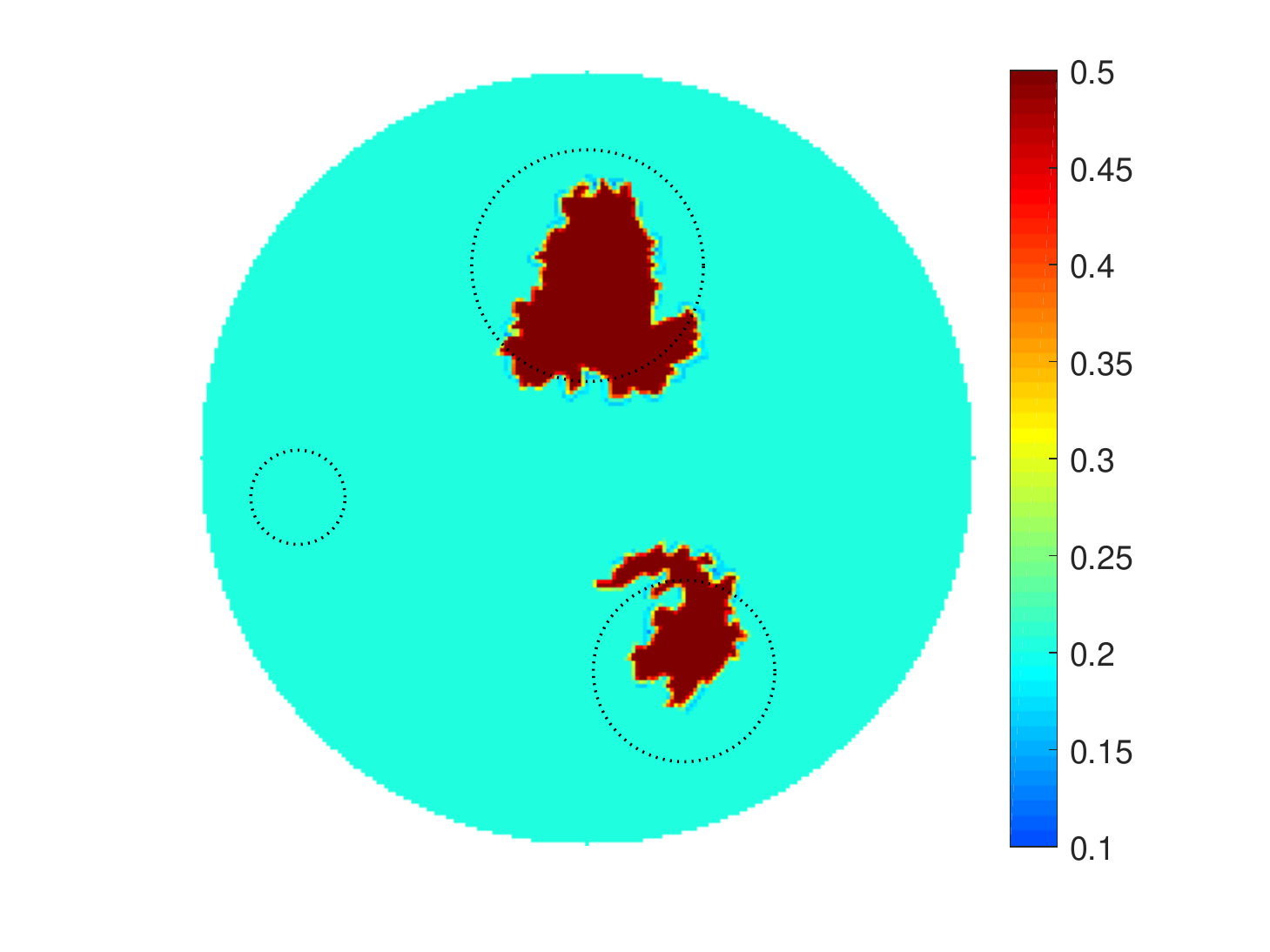}}
  \subfigure[IPOPT]{\includegraphics[width=0.25\textwidth]{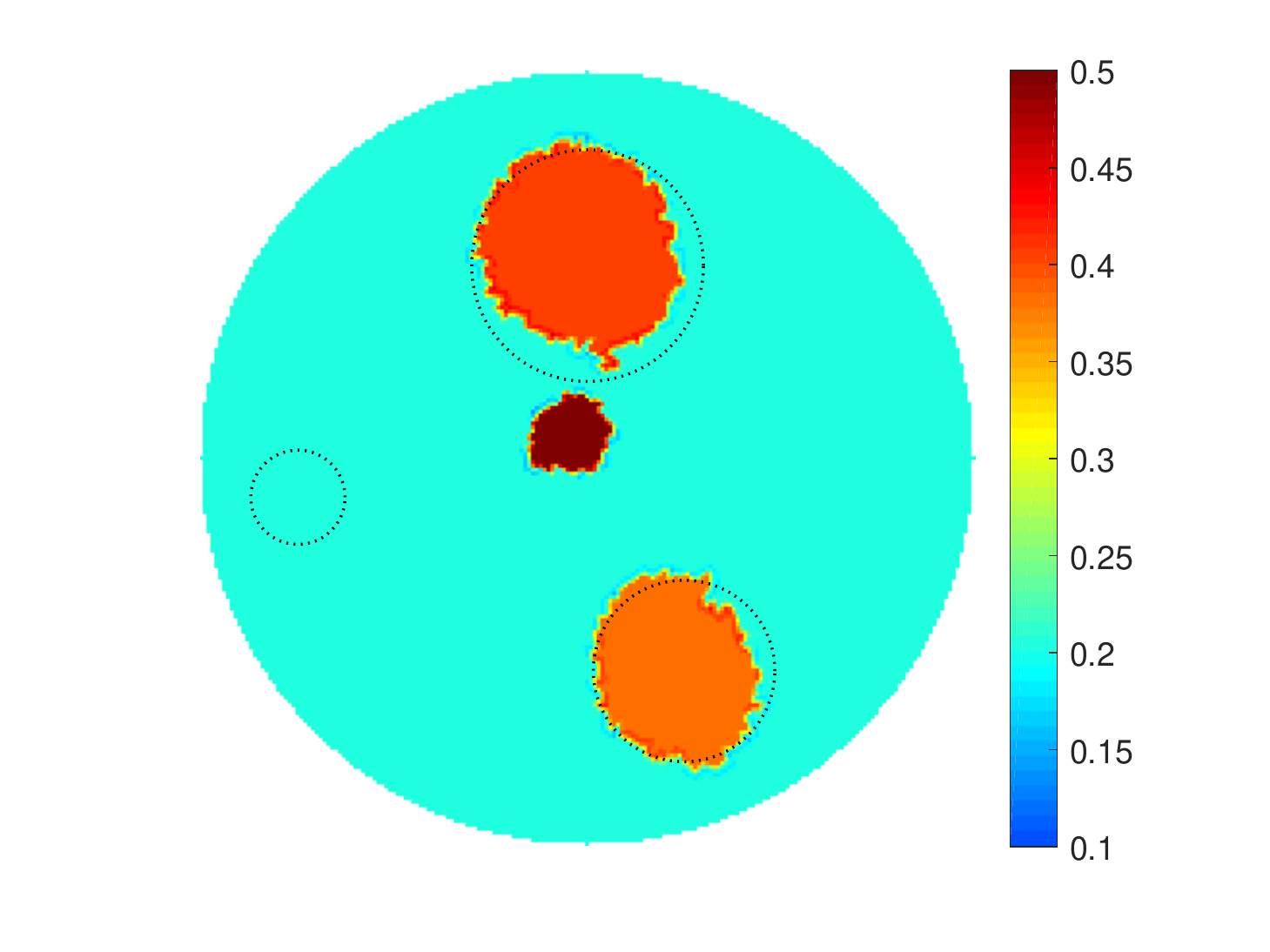}}
  \subfigure[SNOPT]{\includegraphics[width=0.25\textwidth]{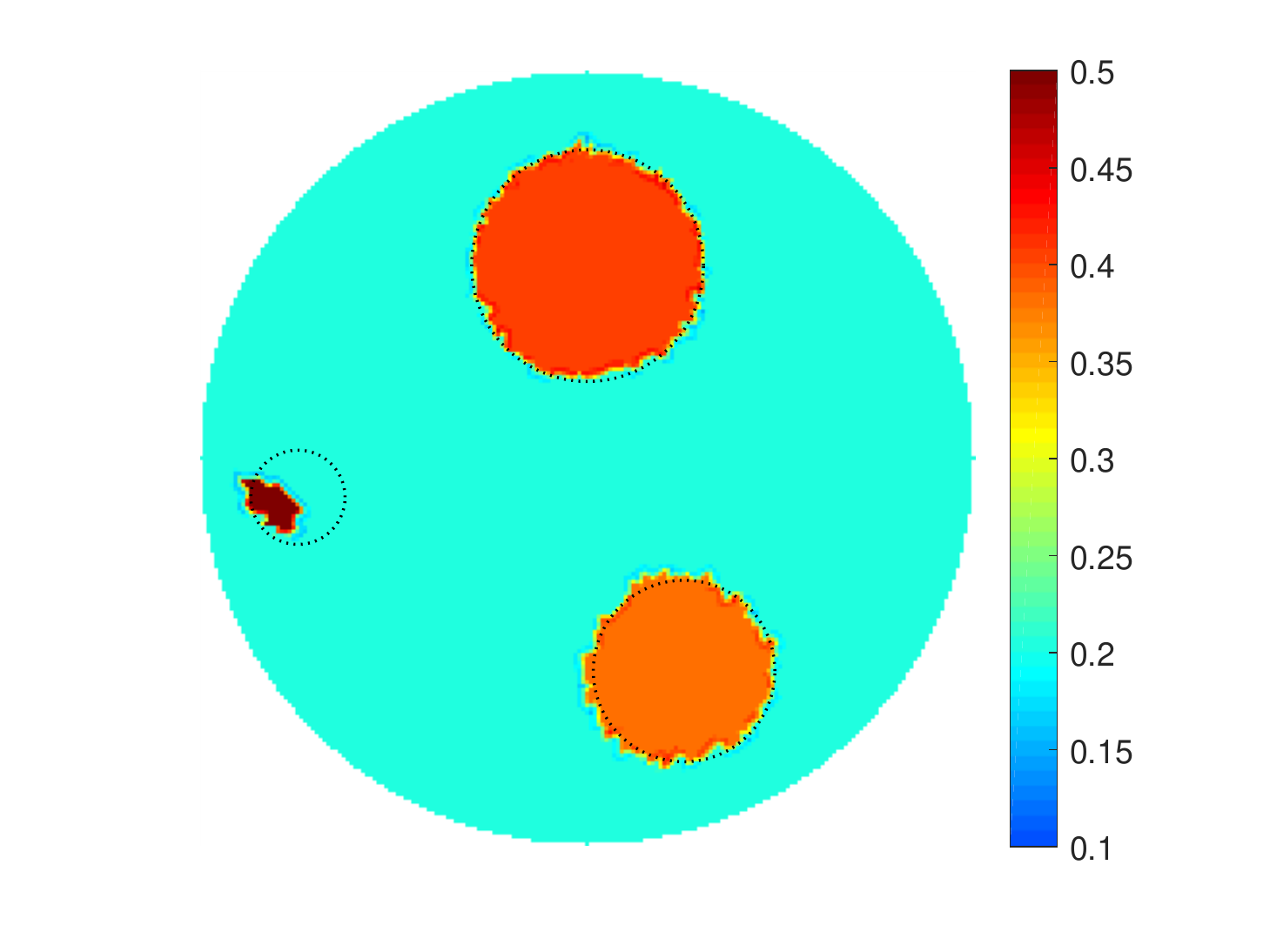}}
  \subfigure[CD]{\includegraphics[width=0.25\textwidth]{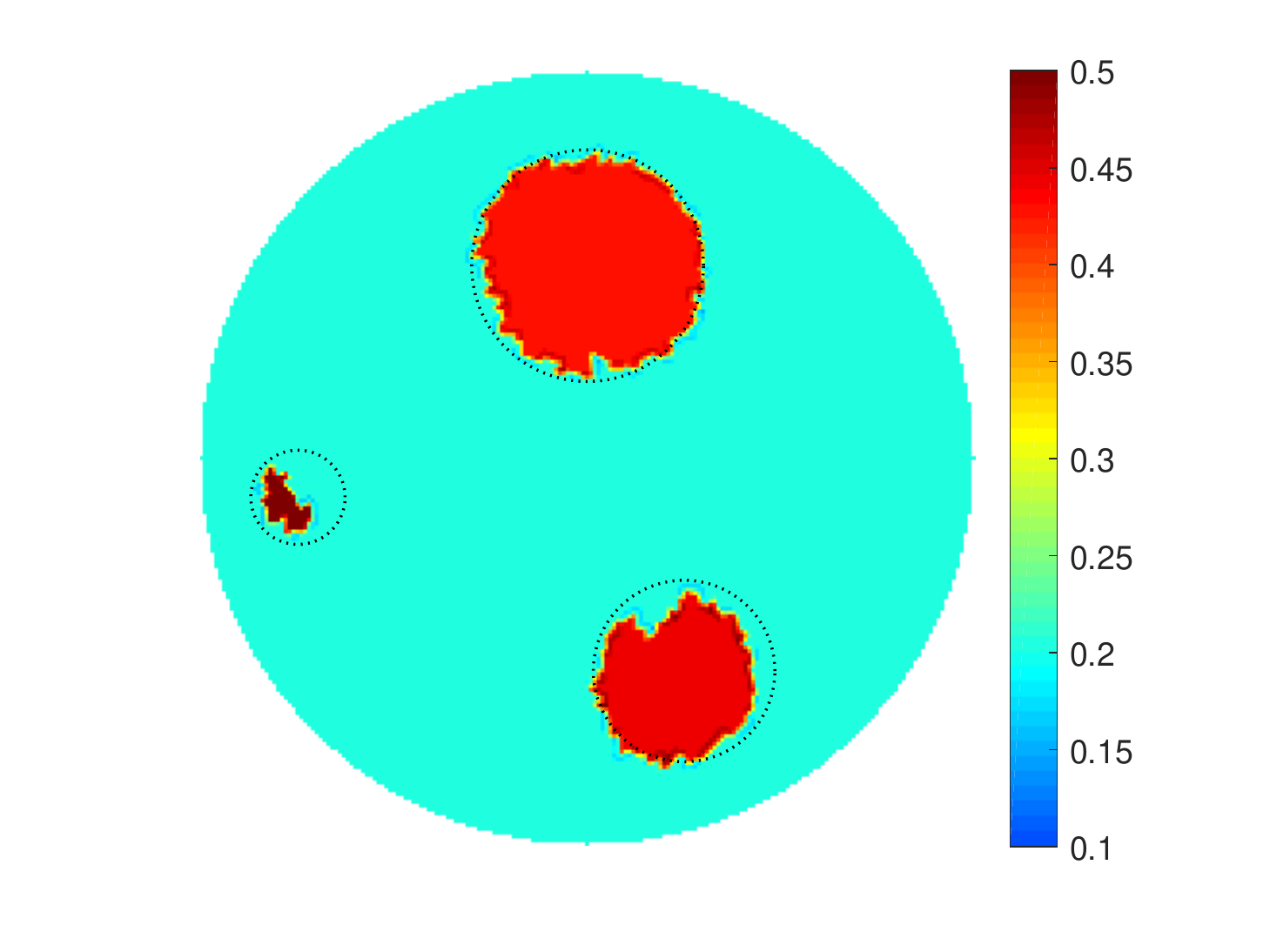}}}
  \mbox{
  \subfigure[]{\includegraphics[width=0.33\textwidth]{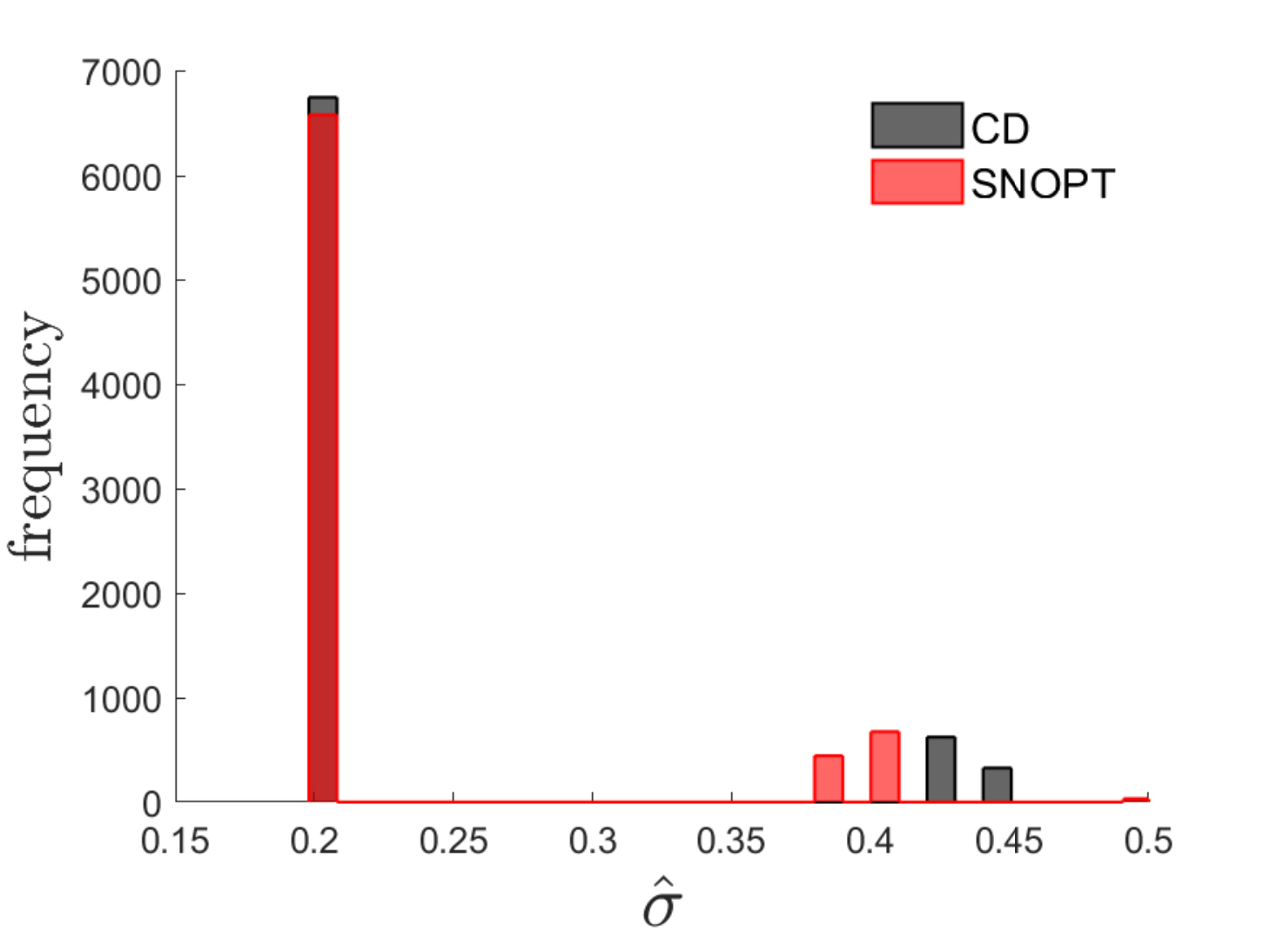}}
  \subfigure[SNOPT]{\includegraphics[width=0.33\textwidth]{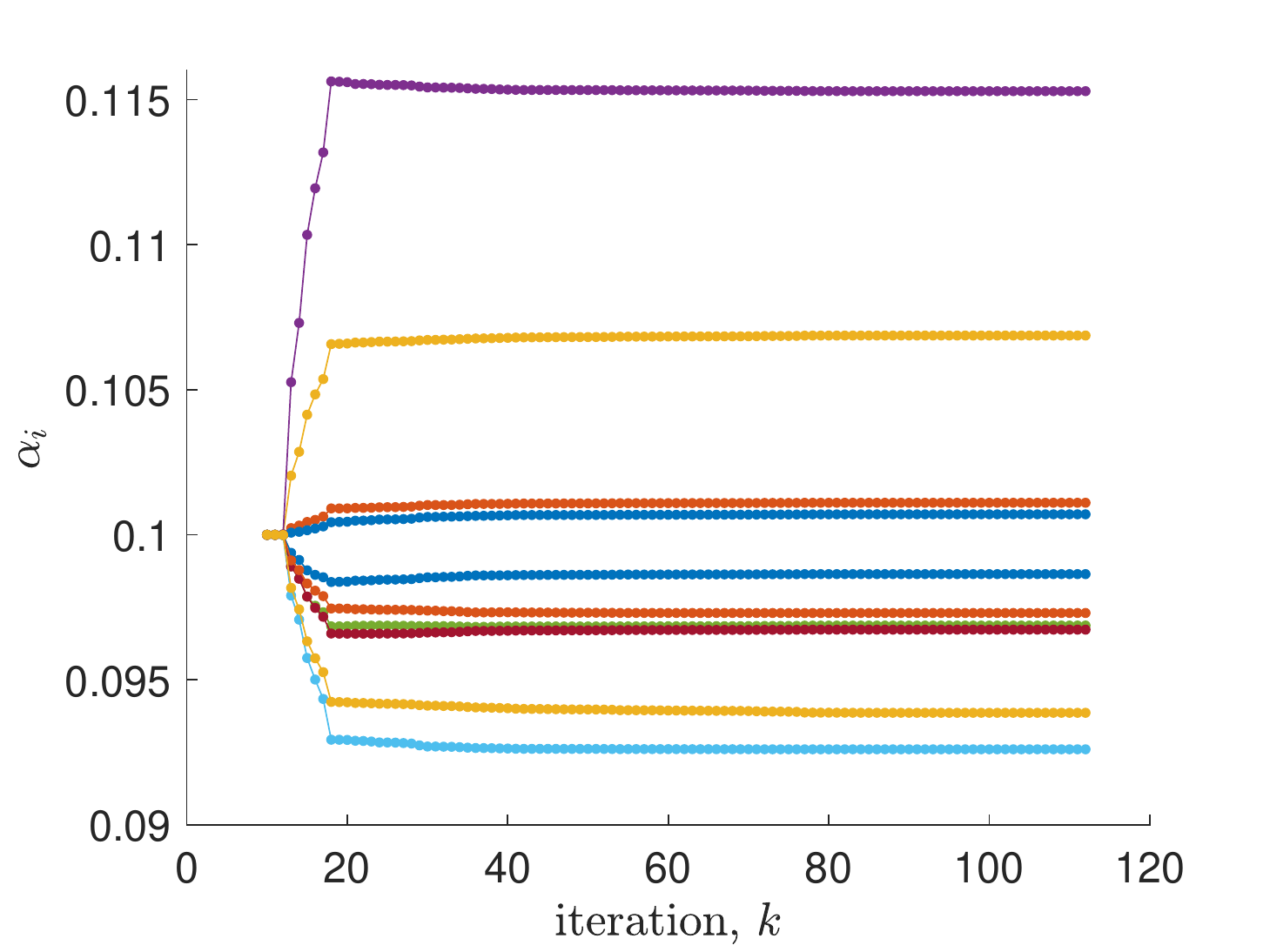}}
  \subfigure[CD]{\includegraphics[width=0.33\textwidth]{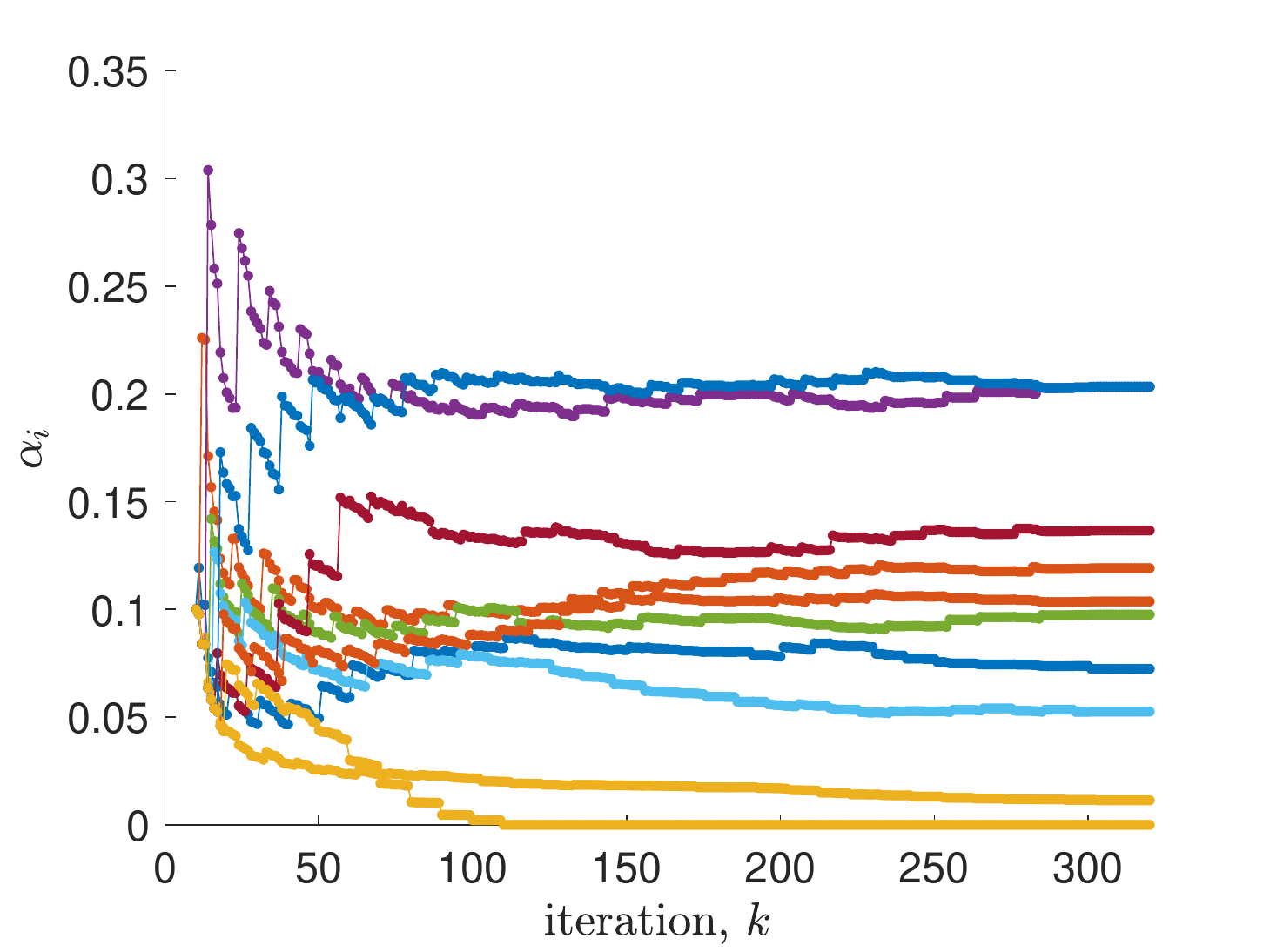}}}
  \end{center}
  \caption{Model~\#1. (a-h)~Solution images obtained by (a,e)~MMA, (b,f)~{\tt IPOPT},
    (c,g)~{\tt SNOPT}, and (d,h)~CD after completing (a-d)~Step~2 and (e-h)~Step~3.
    (i)~Histograms for solutions obtained by (black) CD and (red) {\tt SNOPT}.
    The dashed circles are added to represent the location of cancer-affected regions
    taken from known $\sigma_{true}(x)$ in Figure~\ref{fig:model}(a).
    (j,k)~History of changes in weights $\alpha_i$ when applying (j)~{\tt SNOPT} and
    (k)~CD methods.}
  \label{fig:model_8_b}
\end{figure}

We are particularly interested in the approach that allows fast and accurate
fine-scale images at Step~2 and also accurate coarse-scale binary images to be
of comparable quality after applying the Step~3 tuning procedure. The sequential
quadratic optimizer {\tt SNOPT} proves its superior suitability for this task
while competing with its predecessor (CD) and other methods that use the same
gradients. Figures~\ref{fig:model_8_b}(a-h) provide the results obtained by all
used methods on both fine and coarse scales. Figure~\ref{fig:model_8_b}(g)
confirms an almost ideal reconstruction of big spots (both for color and shape)
and rather satisfactory (due to the size comparable with the size of the applied
boundary electrodes) quality of the smallest spot image. Figure~\ref{fig:model_8_b}(i)
contributes to this conclusion by comparing the reconstructed values of
$\hat \sigma_{high}$ that are much closer to the ``known'' value $\sigma_c = 0.4$
in the case of gradient-based {\tt SNOPT}. Finally, we reiterate and conclude on the
reasons for the improved computational speed (1,666 vs.~18,112 cost functional
evaluations for {\tt SNOPT} and CD, respectively). Evidently, gradient-based methods
are faster as they change all (or, at least, many) controls while CD works only with
one control at a time. In addition to this, as seen in Figures~\ref{fig:model_8_b}(j)
and \ref{fig:model_8_b}(k), CD spends some time ``re-ranking'' all $N_s$ samples
used to update the reconstructed image by changing their weights $\alpha_i$. On
the other hand, based on the sensitivity provided in $\bnabla_{\alpha} \mJ$,
{\tt SNOPT} sets the ``ranks" at the beginning of Step~2 and focuses on their
qualities afterward.

\subsection{Validation on Complicated Models}
\label{sec:frame_compl}

Based on the results obtained for our model~\#1 and described in the previous
section, we conclude on the superior performance brought to the proposed
computational framework by employing a multiscale gradient-based search and
control space reduction for binary tuning. As shown by practical experiments,
the efficacy of the former is subject to a particular optimizer, and the
suitability of the latter appears questionable for more complicated models.
Therefore, in this section, we discuss the results obtained using this new
framework now applied to models with a significantly increased level of
complexity to explore any bounds for its applicability. So far, the new
algorithm confirms its ability to accurately reconstruct circular-shaped
cancerous regions of various sizes at multiple locations. Here, the added
complications are non-circular shapes for those regions and their varied
conductivities $\sigma_c$. The rest of our computational results will compare
the performance of the {\tt SNOPT} and CD optimizers based on data with
0.5\% noise.

\begin{figure}[htb!]
  \begin{center}
  \mbox{
  \subfigure[model~\#2]{\includegraphics[width=0.33\textwidth]{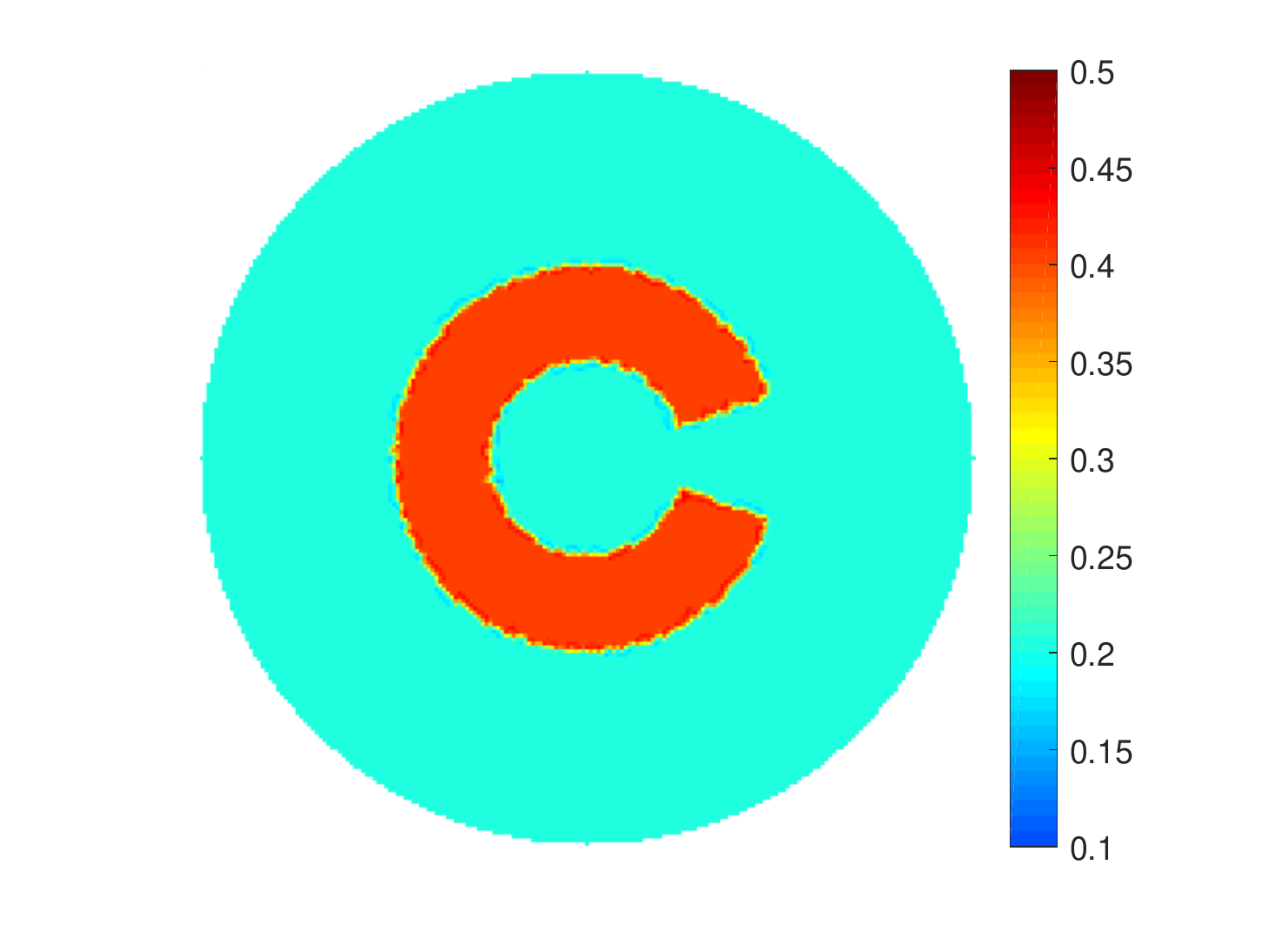}}
  \subfigure[SNOPT: Step~2]{\includegraphics[width=0.33\textwidth]{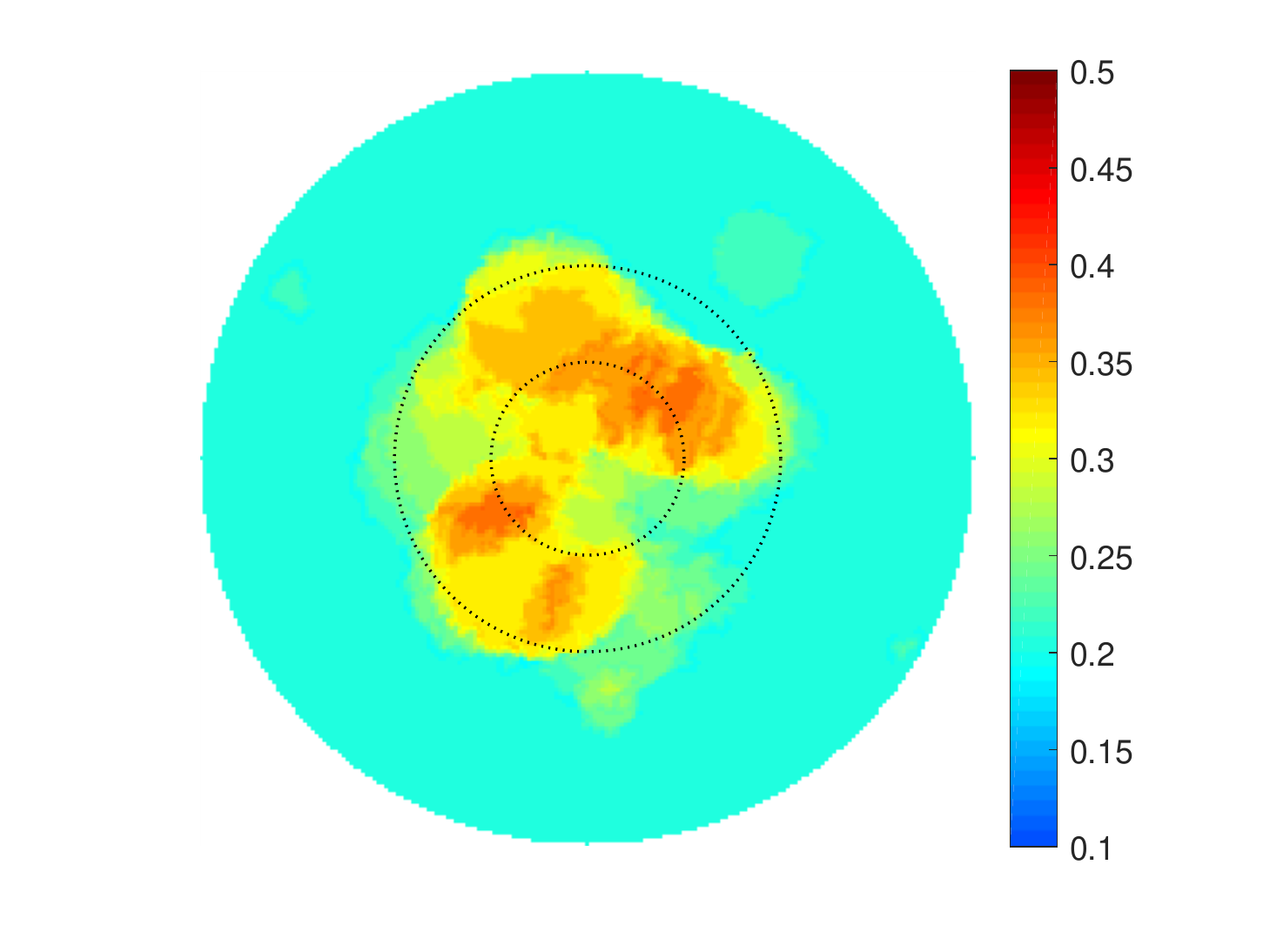}}
  \subfigure[SNOPT: Step~3]{\includegraphics[width=0.33\textwidth]{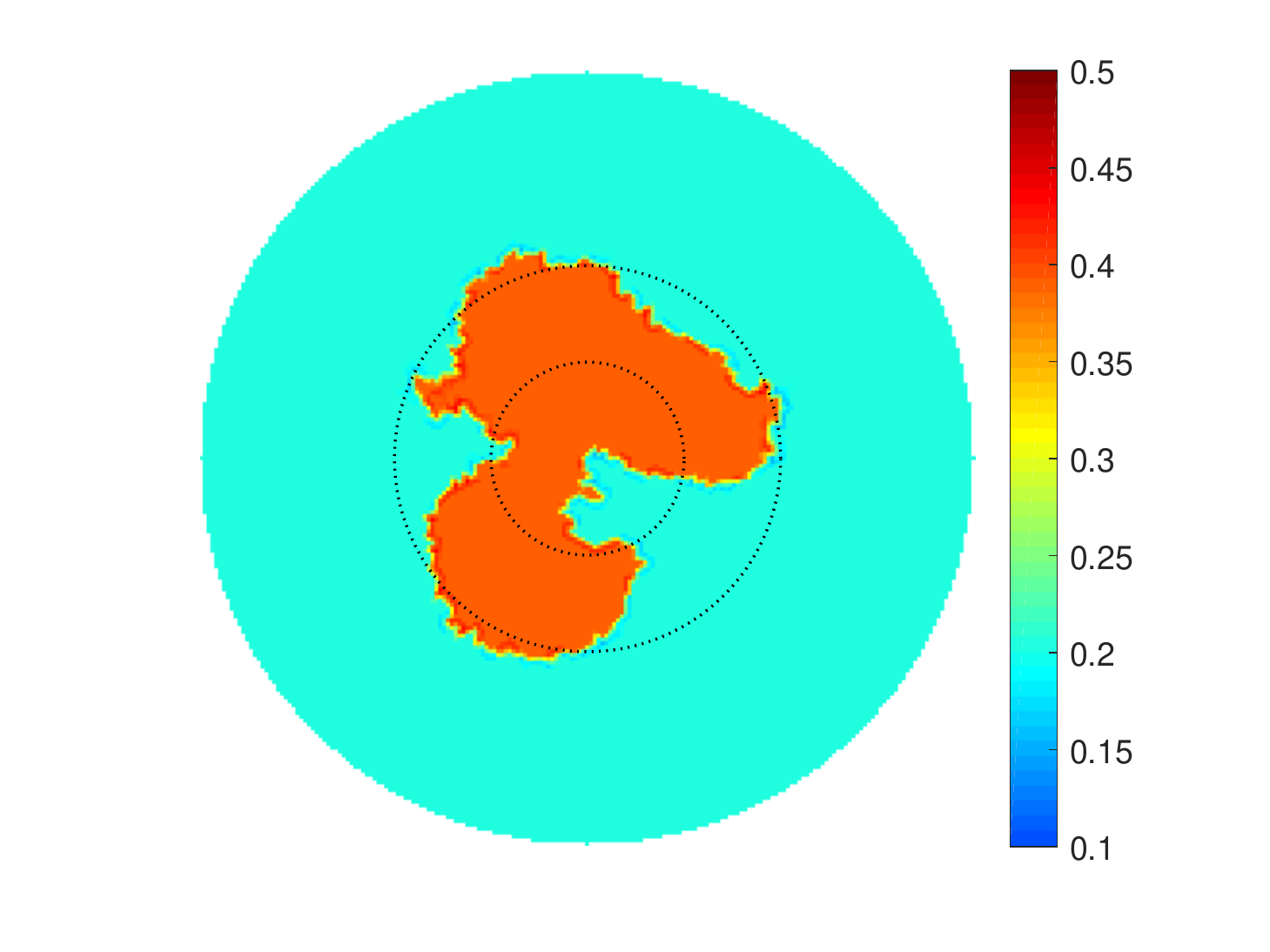}}}
  \mbox{
  \subfigure[solution error]{\includegraphics[width=0.33\textwidth]{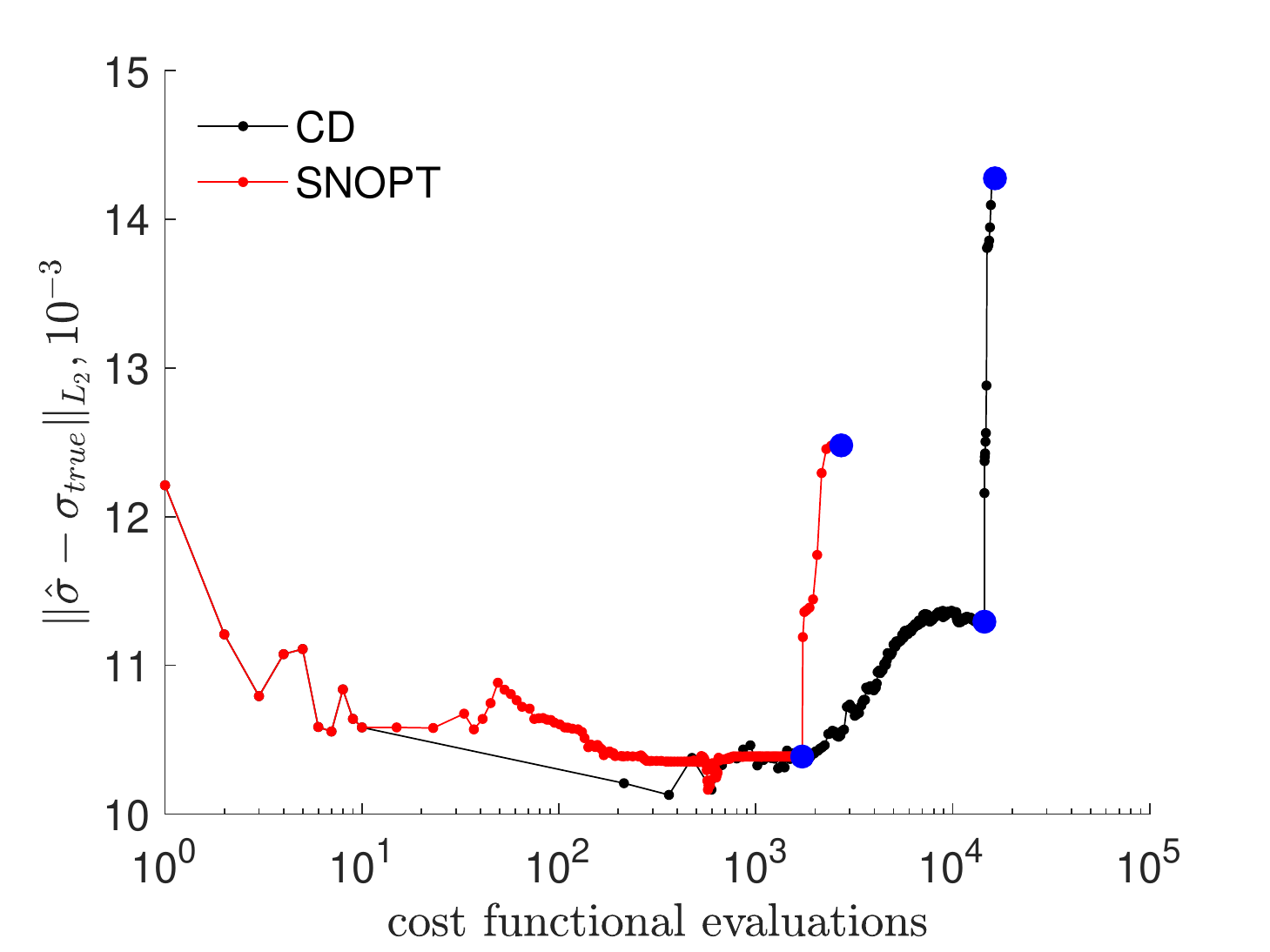}}
  \subfigure[CD: Step~2]{\includegraphics[width=0.33\textwidth]{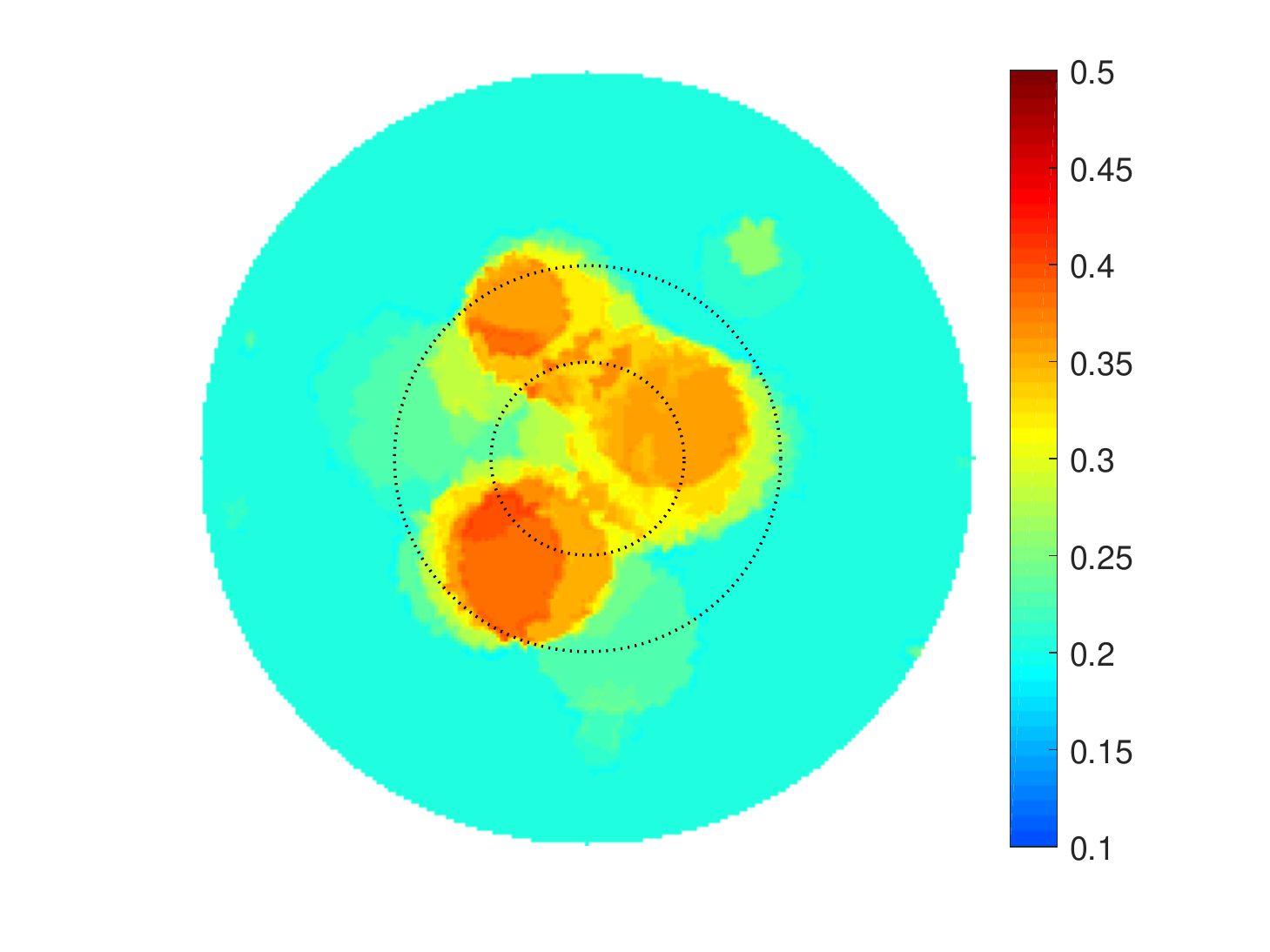}}
  \subfigure[CD: Step~3]{\includegraphics[width=0.33\textwidth]{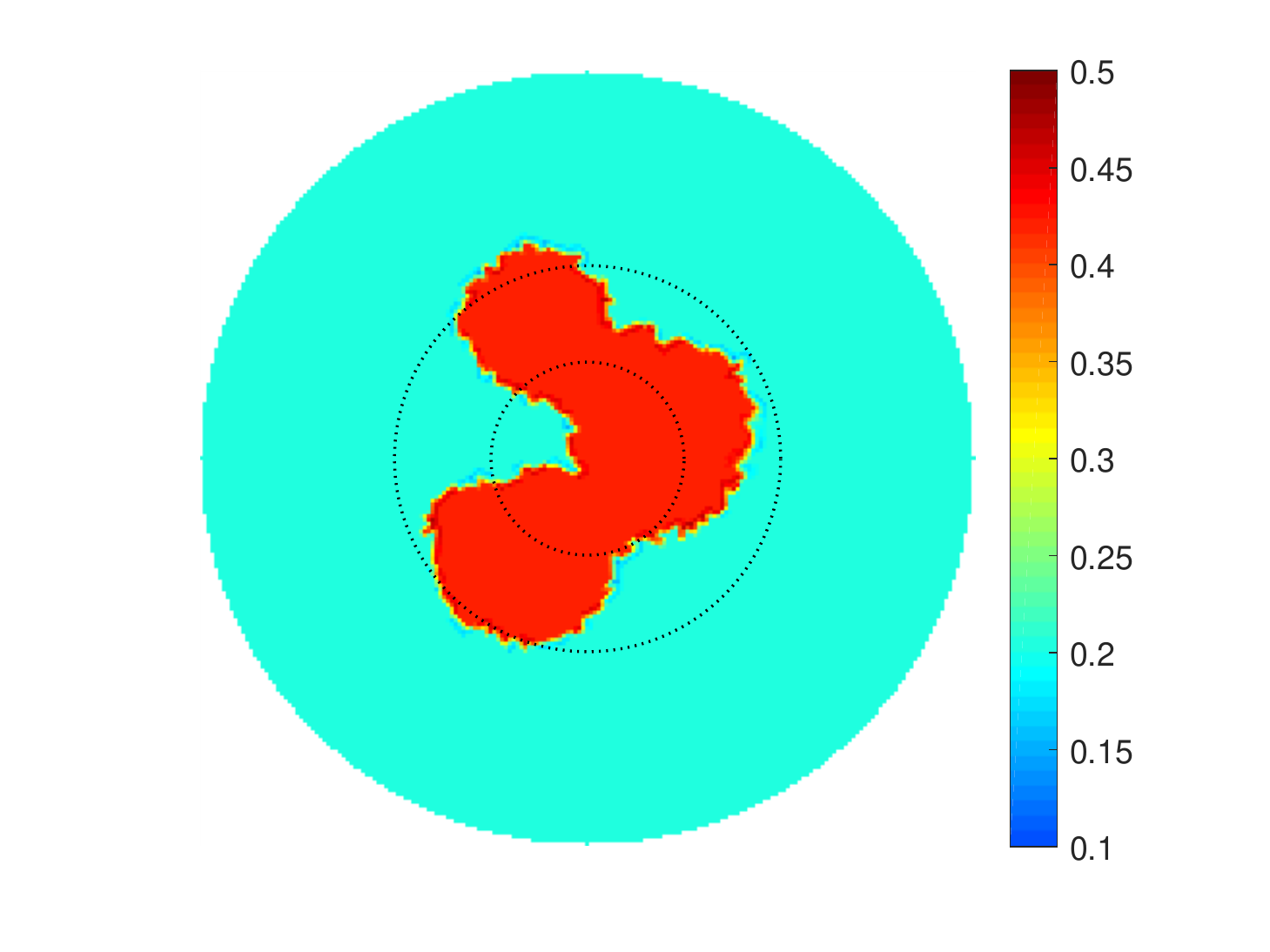}}}
  \end{center}
  \caption{(a)~EIT model~\#2: true electrical conductivity $\sigma_{true}(x)$.
    (b,c,e,f)~Solution images obtained by (b,c)~{\tt SNOPT} and (e,f)~CD after completing
    (b,e)~Step~2 and (c,f)~Step~3. The dashed circles are added to represent the C-shape
    of the cancer-affected region taken from known $\sigma_{true}(x)$ in (a).
    (d)~Solution errors $\| \sigma^k - \sigma_{true} \|_{L_2}$
    as functions of a number of cost functional evaluations evaluated while employing
    (black) CD and (red) {\tt SNOPT} optimizers. Blue dots represent solutions obtained
    after Step~2 and 3 phases are complete.}
  \label{fig:model_15}
\end{figure}

We create our next model (\#2) to check the method's performance when
cancer-affected areas depart from circular shapes, e.g., as the one in
Figure~\ref{fig:model_15}(a), showing an image with a C-shape region. Although
using the same collection $\mC(10,0000)$ of samples containing only circles
of various radii does not assume such shapes, the reconstructions obtained
by both {\tt SNOPT} and CD seem acceptable; see Figures~\ref{fig:model_15}(c)
and \ref{fig:model_15}(f), respectively. Similar to model~\#1, gradient-based
{\tt SNOPT} shows significantly better performance for quality and computational
speed at both Step~2 and 3 optimization phases; refer to
Figure~\ref{fig:model_15}(d) showing the solution error and blue dots
representing solutions obtained after both phases are complete. Images obtained
after Step~2, in Figures~\ref{fig:model_15}(b) and \ref{fig:model_15}(e), confirm
the superior performance added by using gradients improved by the Step~3 tuning
procedure that indicates its suitability for such shapes. As explored in
\cite{ArbicBukshtynov2022,ArbicMS2020}, this performance may be further enhanced,
even in the presence of noise, once we increase the maximum number of circles
$N_{c, max}$ in the samples. It proves the potential of the proposed methodology
in applications with models bearing rather complex geometry.

\begin{figure}[htb!]
  \begin{center}
  \mbox{
  \subfigure[model~\#3]{\includegraphics[width=0.33\textwidth]{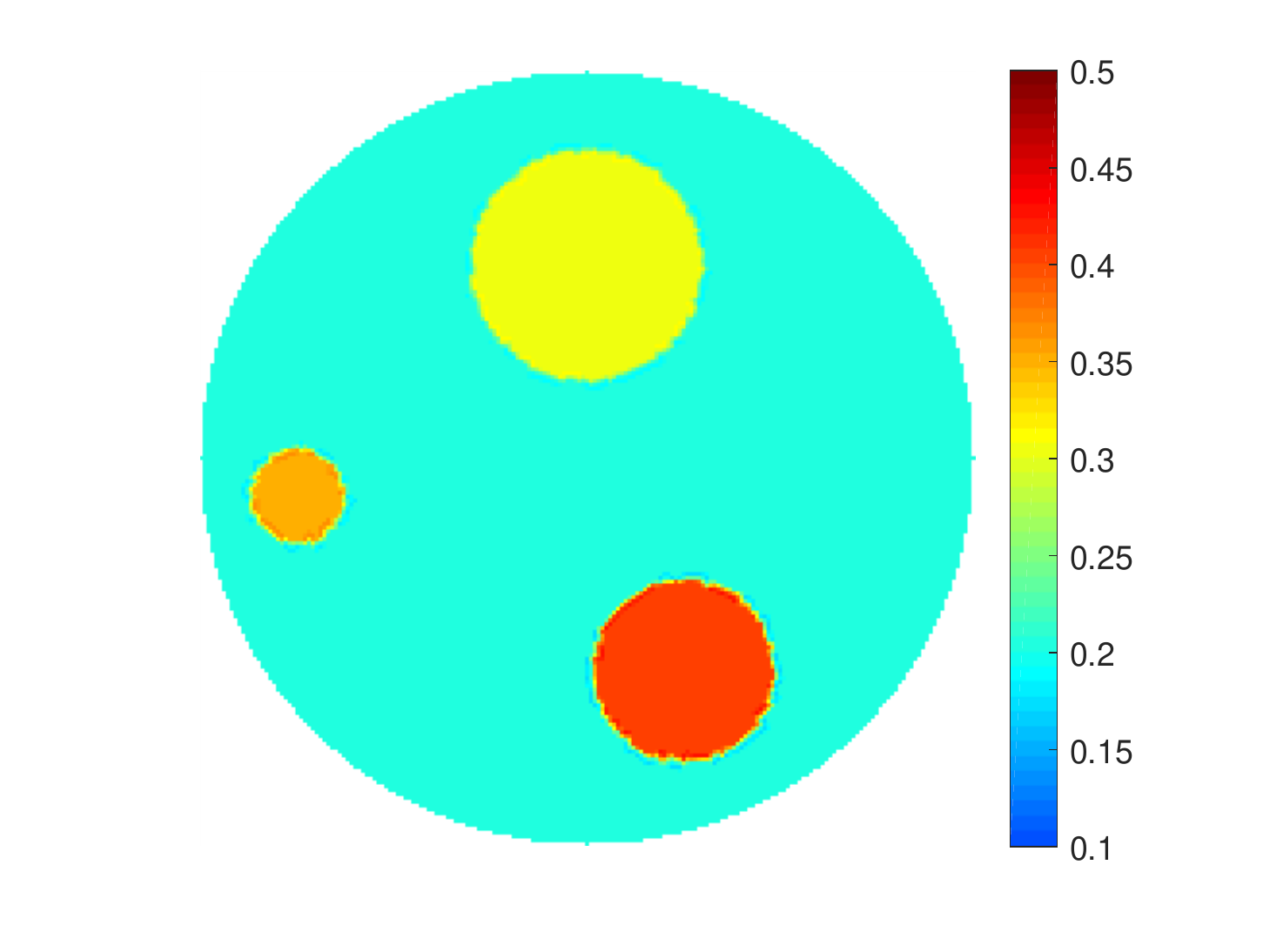}}
  \subfigure[SNOPT: Step~2]{\includegraphics[width=0.33\textwidth]{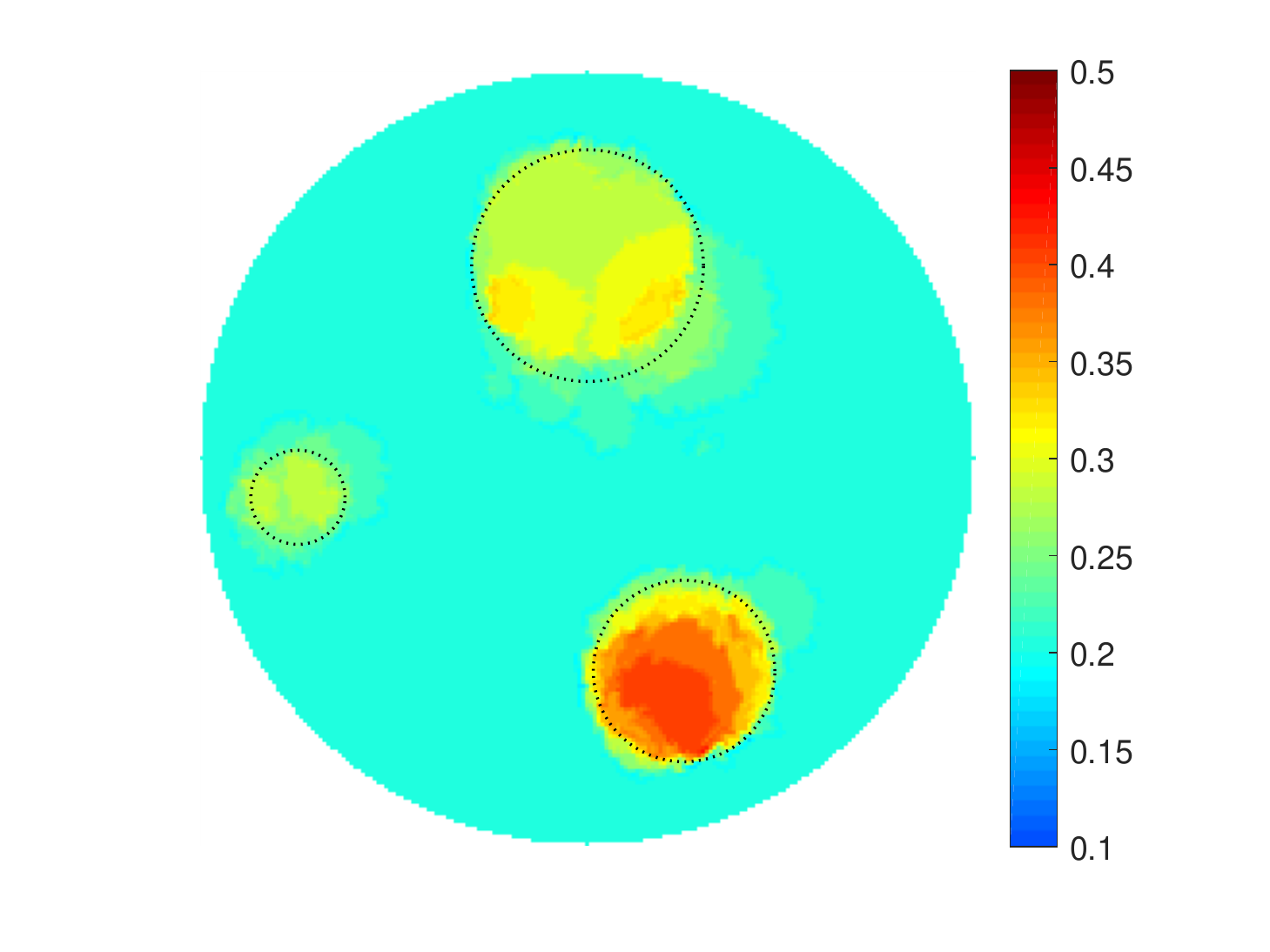}}
  \subfigure[SNOPT: Step~3]{\includegraphics[width=0.33\textwidth]{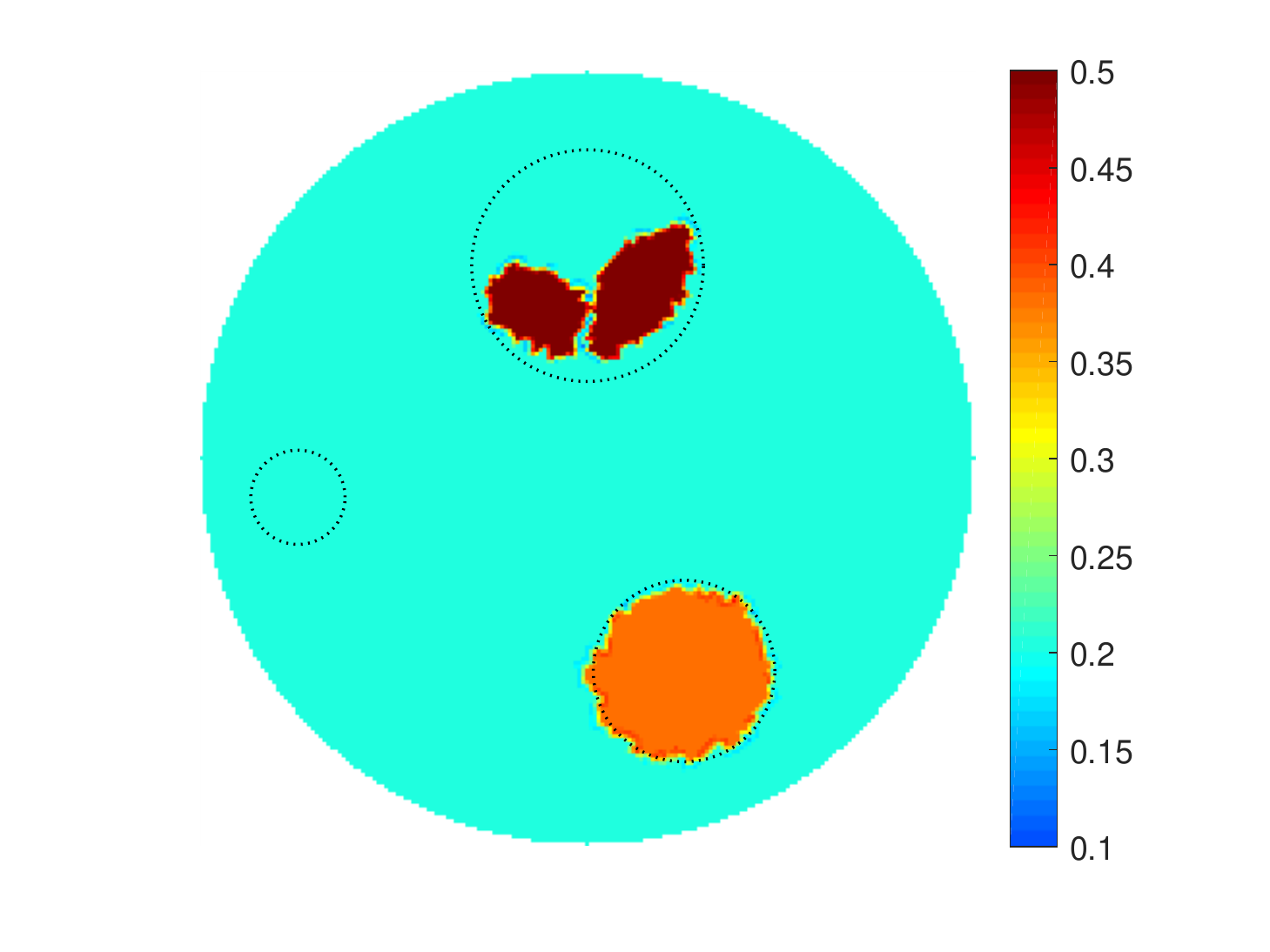}}}
  \mbox{
  \subfigure[solution error]{\includegraphics[width=0.33\textwidth]{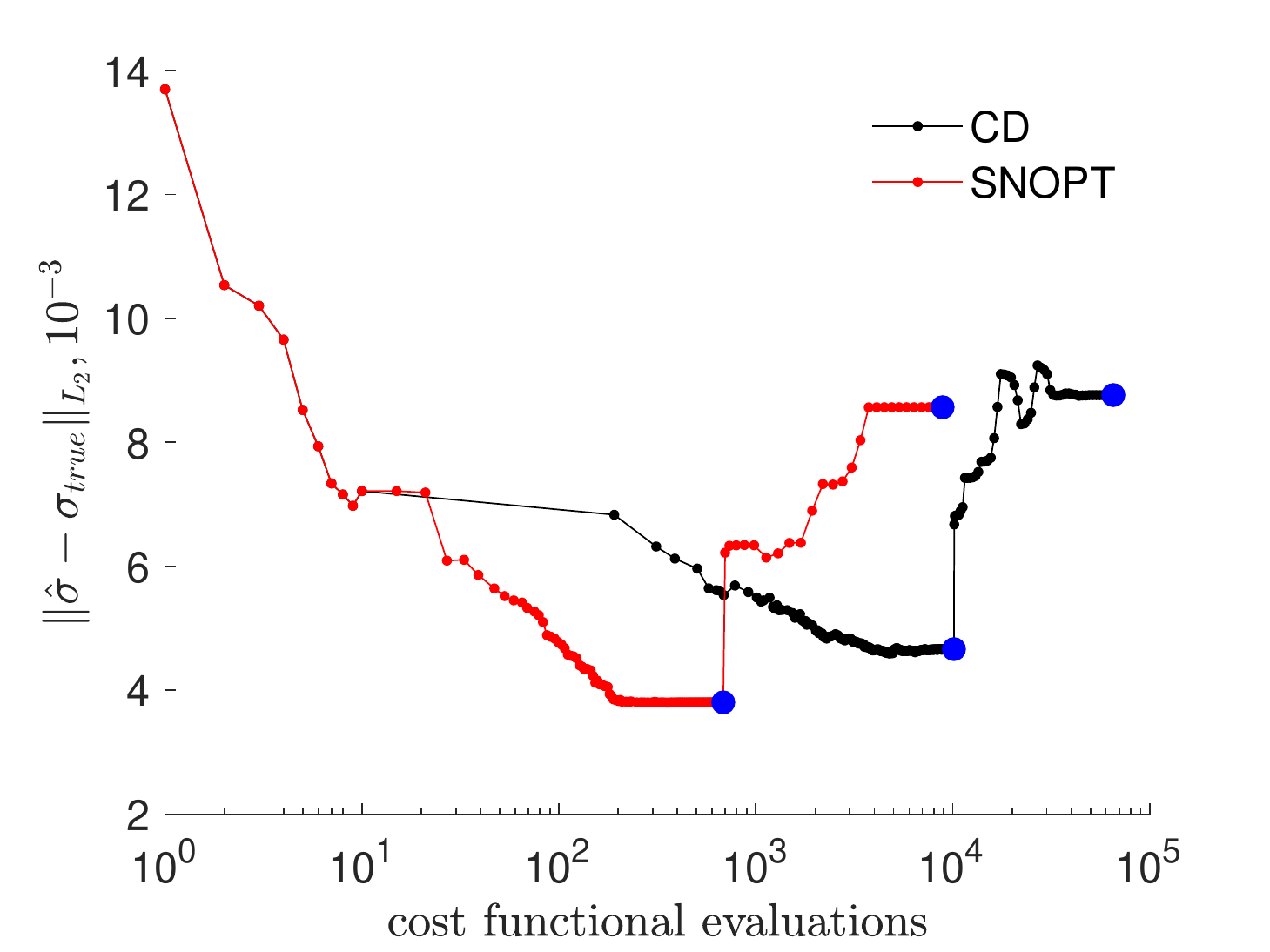}}
  \subfigure[CD: Step~2]{\includegraphics[width=0.33\textwidth]{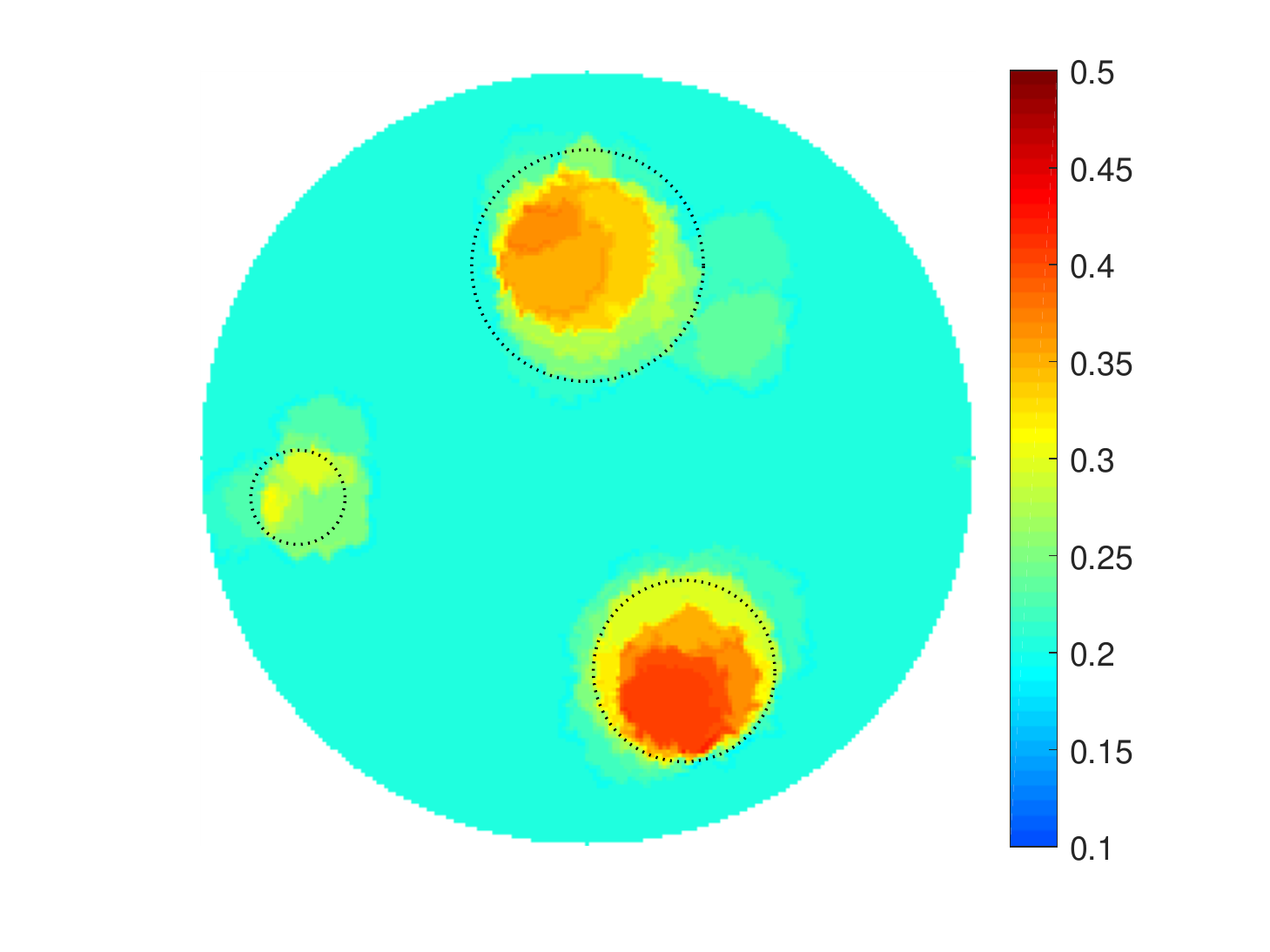}}
  \subfigure[CD: Step~3]{\includegraphics[width=0.33\textwidth]{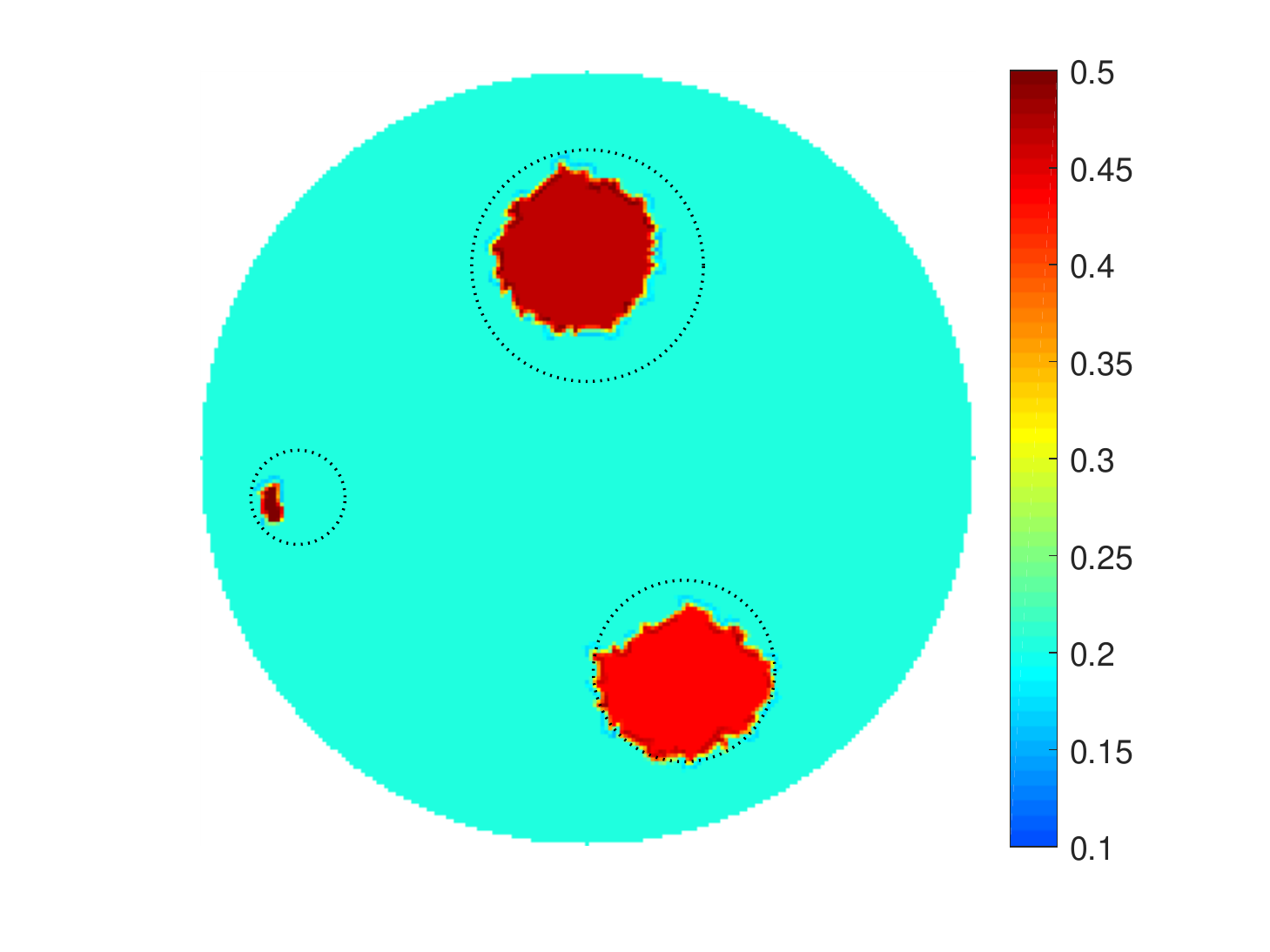}}}
  \end{center}
  \caption{(a)~EIT model~\#3: true electrical conductivity $\sigma_{true}(x)$.
    (b,c,e,f)~Solution images obtained by (b,c)~{\tt SNOPT} and (e,f)~CD after completing
    (b,e)~Step~2 and (c,f)~Step~3. The dashed circles are added to represent the location
    of cancer-affected regions taken from known $\sigma_{true}(x)$ in (a).
    (d)~Solution errors $\| \sigma^k - \sigma_{true} \|_{L_2}$
    as functions of a number of cost functional evaluations evaluated while employing
    (black) CD and (red) {\tt SNOPT} optimizers. Blue dots represent solutions obtained
    after Step~2 and 3 phases are complete.}
  \label{fig:model_18}
\end{figure}

To further experiment with the applicability and performance of various parts in
the proposed computational algorithm, we modified our model~\#1 by changing the
electrical conductivity (high values) inside the cancerous spots while keeping
the same their mutual positioning and sizes. Figure~\ref{fig:model_18}(a) displays
model~\#1 after this modification (model~\#3), where we set $\sigma_c$ to 0.3, 0.4,
and 0.35 for the big, medium-size, and small spots, respectively.

The visual analysis of the fine-scale solutions obtained using gradient-based
({\tt SNOPT}) and derivative-free (CD) searches after completing Step~2 reveals
good results in recovering the positions of all three spots and reconstructing
their shapes for both methods; refer to Figures~\ref{fig:model_18}(b) and
\ref{fig:model_18}(e), respectively. As previously, {\tt SNOPT} supplied by
gradient information performs better for quality (better match for colors,
especially for the two big spots) and computational speed (683 vs.~10,138 cost
functional evaluations to complete Step~2), confirmed by the analysis of the
solution error provided in Figure~\ref{fig:model_18}(d). It concludes with the
applicability of the proposed methodology to differentiate between cancerous
spots with varied conductivities even though the used samples are ``unaware''
of that, making this methodology less dependent on the prior knowledge of
the simulated phenomena. The same analysis applied to the results for binary
tuning during Step~3, Figures~\ref{fig:model_18}(c) and \ref{fig:model_18}(f),
concludes on the lower quality of obtained images for both methods. We explain
it by the presence of noise (0.5\%) in the supplied data -- the known effect
reported earlier -- procedure in Step~2 has less sensitivity to noise
\cite{ArbicBukshtynov2022,ArbicMS2020} than binary tuning in Step~3
\cite{ChunEdwardsBukshtynov2024}. Although the overall outcomes
are rather satisfactory and promising, it leaves space for further development
for the binary tuning phase and better ``communication'' between computational
components while transitioning from Step~2 to Step~3.

\subsection{Applications to Cancer Detection}
\label{sec:frame_real}

\begin{figure}[htb!]
  \begin{center}
  \mbox{
  \subfigure[model~\#4]{\includegraphics[width=0.33\textwidth]{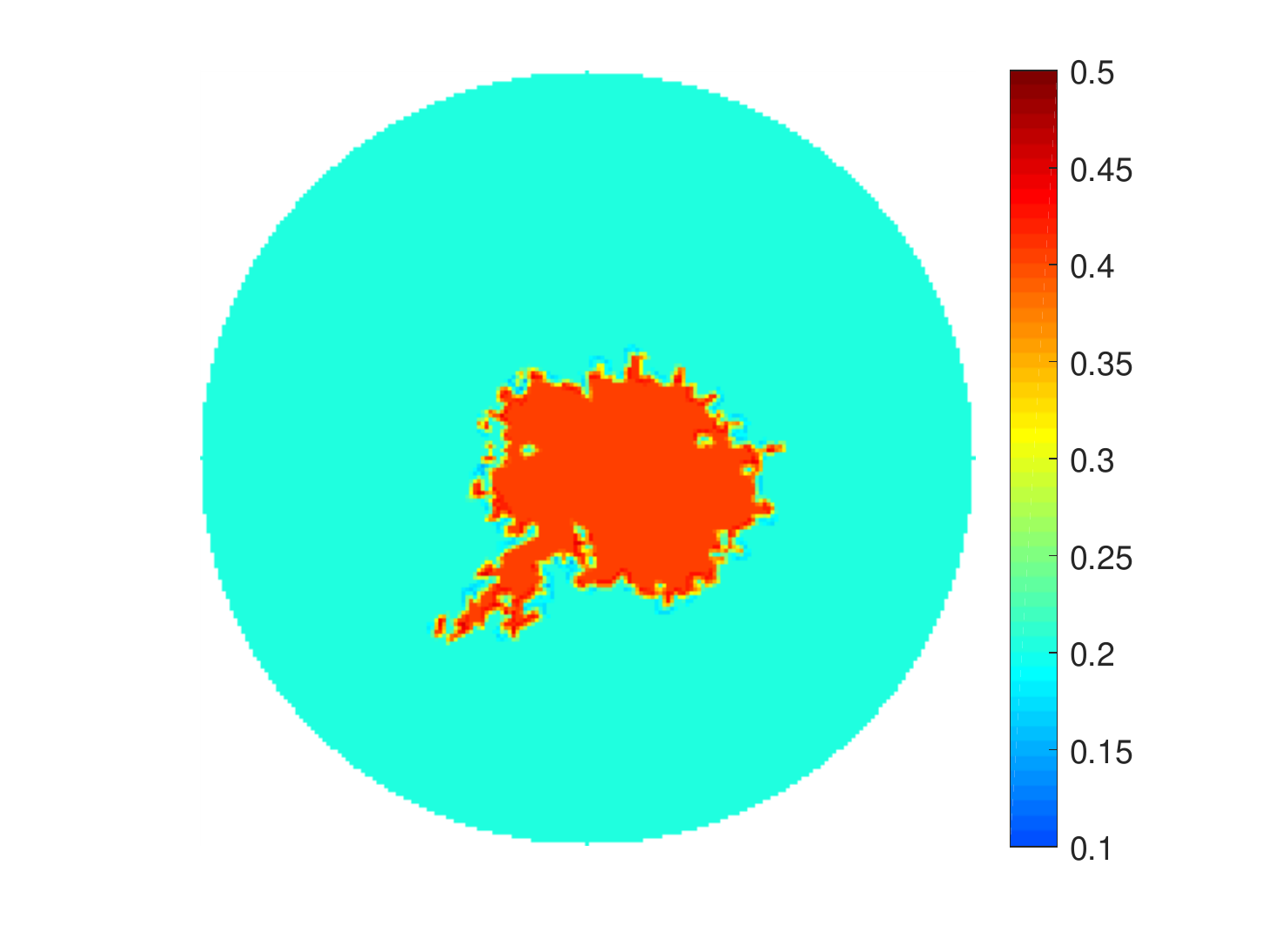}}
  \subfigure[SNOPT: Step~2]{\includegraphics[width=0.33\textwidth]{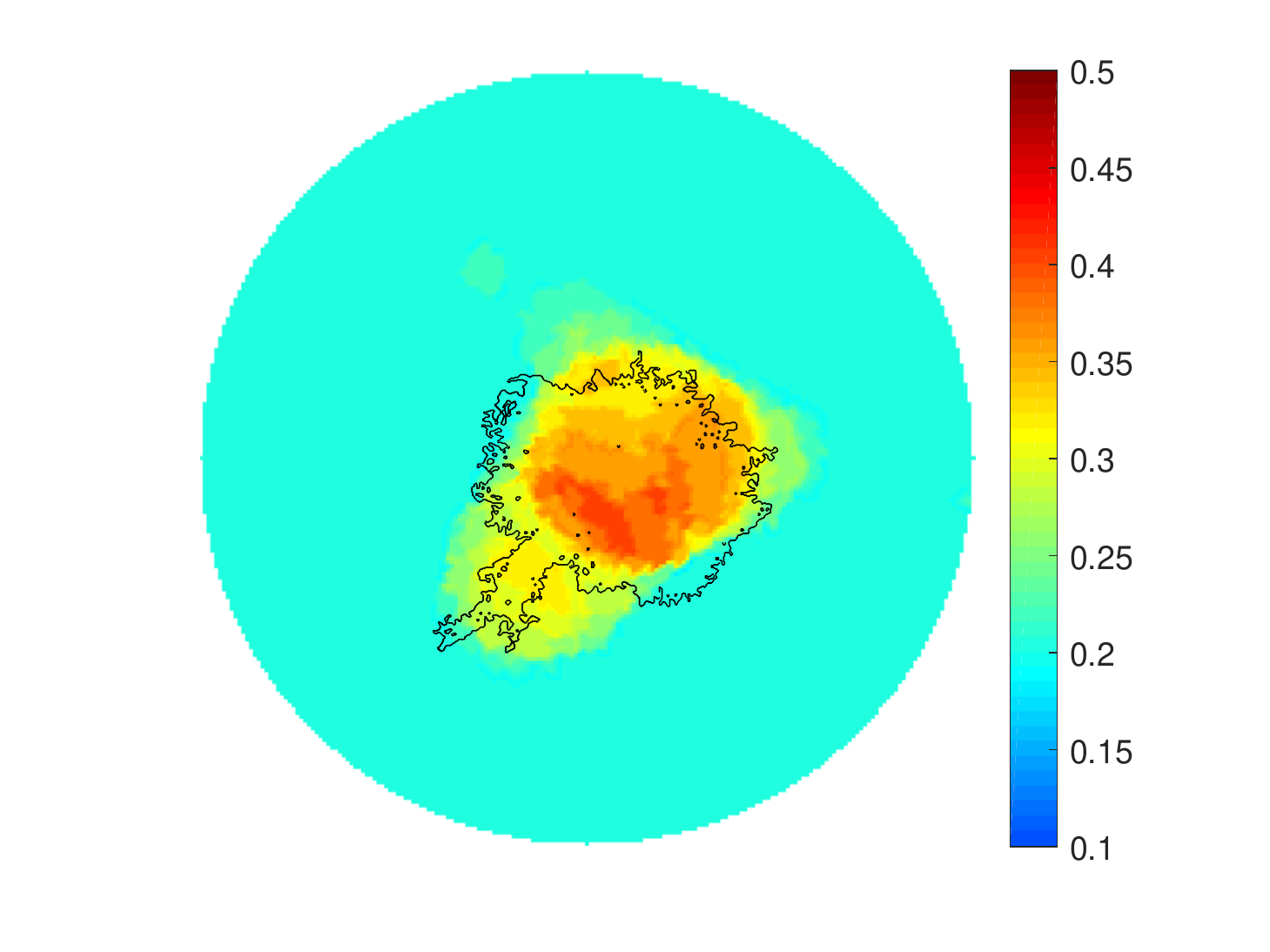}}
  \subfigure[SNOPT: Step~3]{\includegraphics[width=0.33\textwidth]{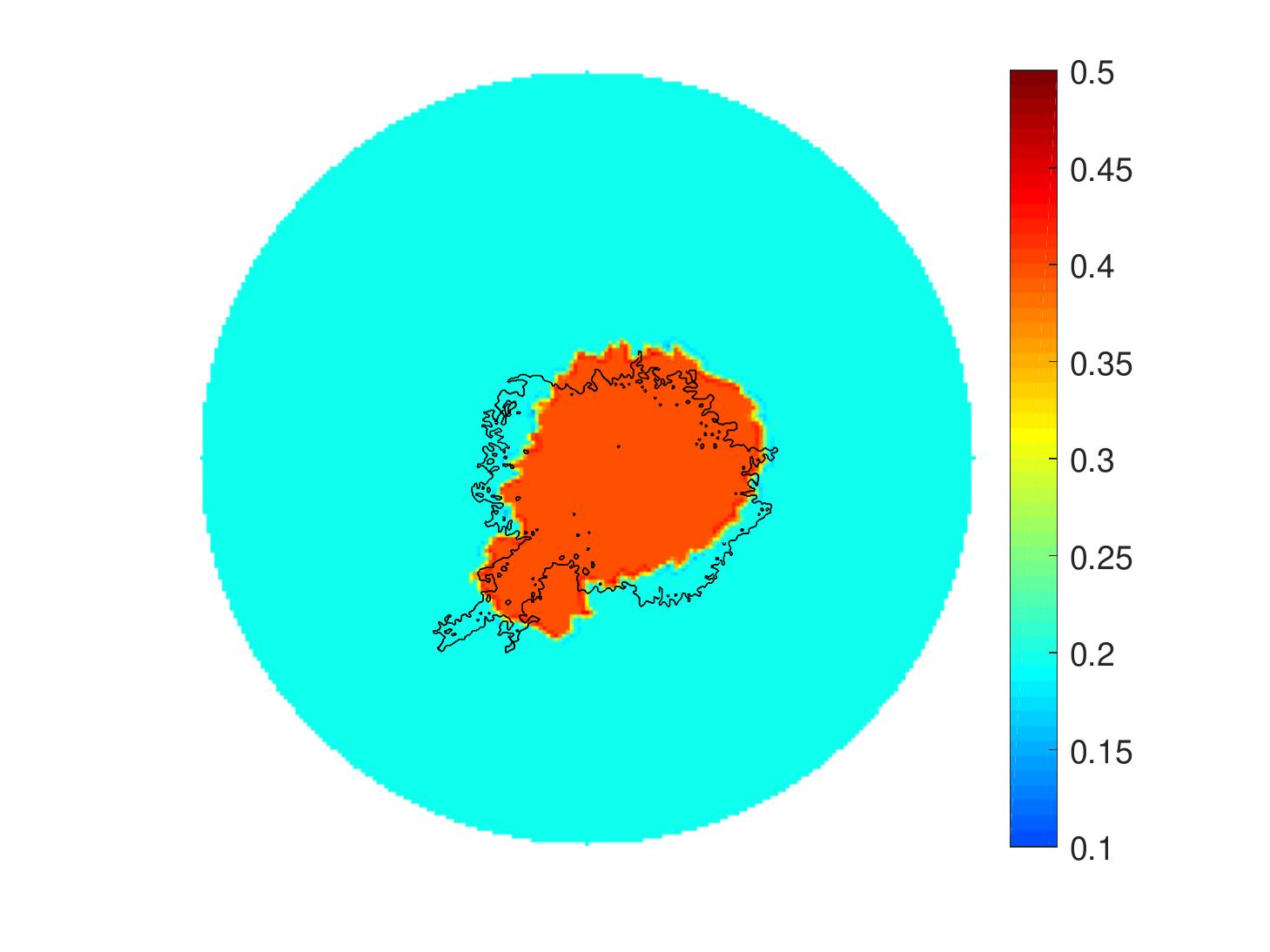}}}
  \mbox{
  \subfigure[solution error]{\includegraphics[width=0.33\textwidth]{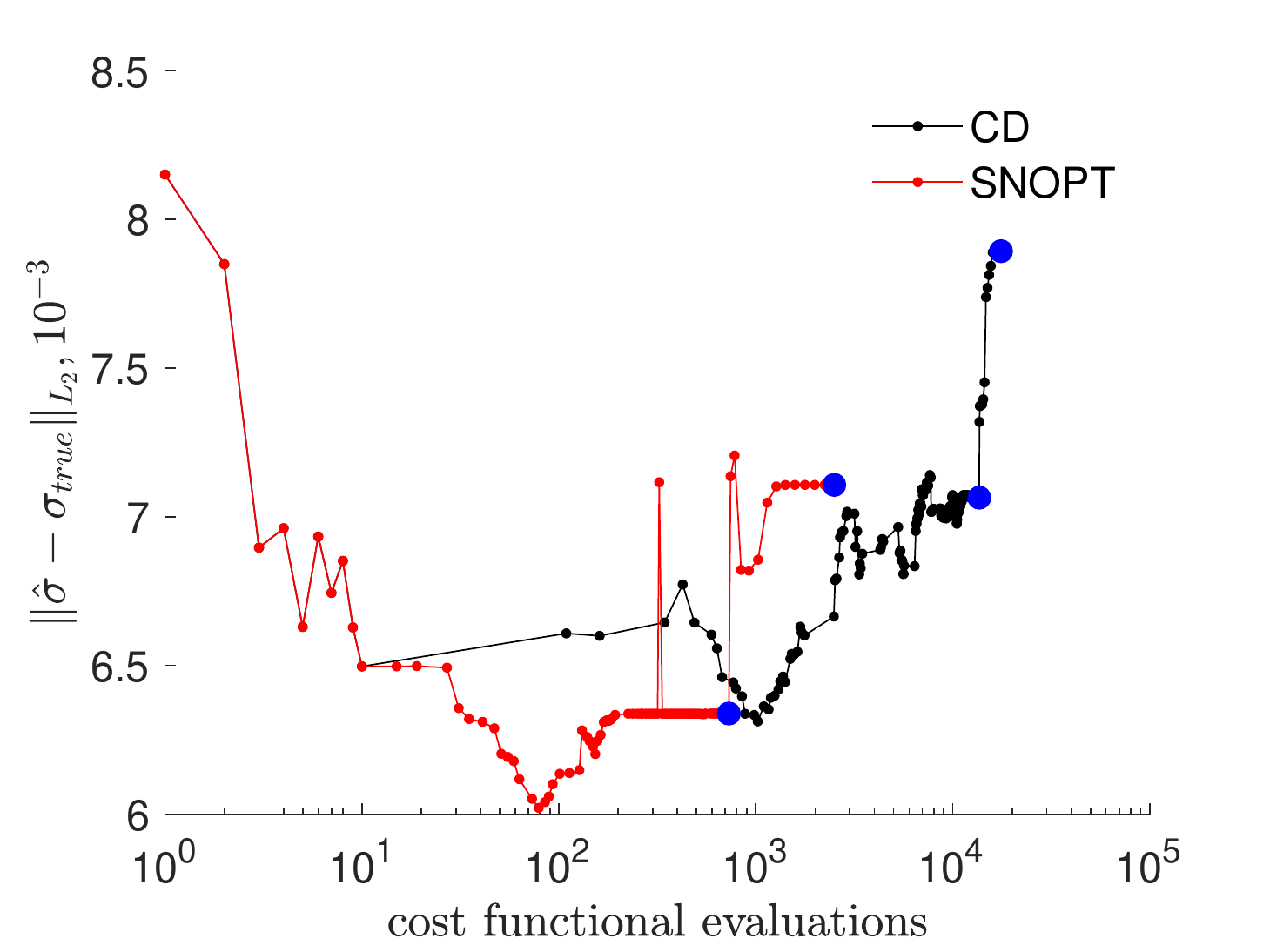}}
  \subfigure[CD: Step~2]{\includegraphics[width=0.33\textwidth]{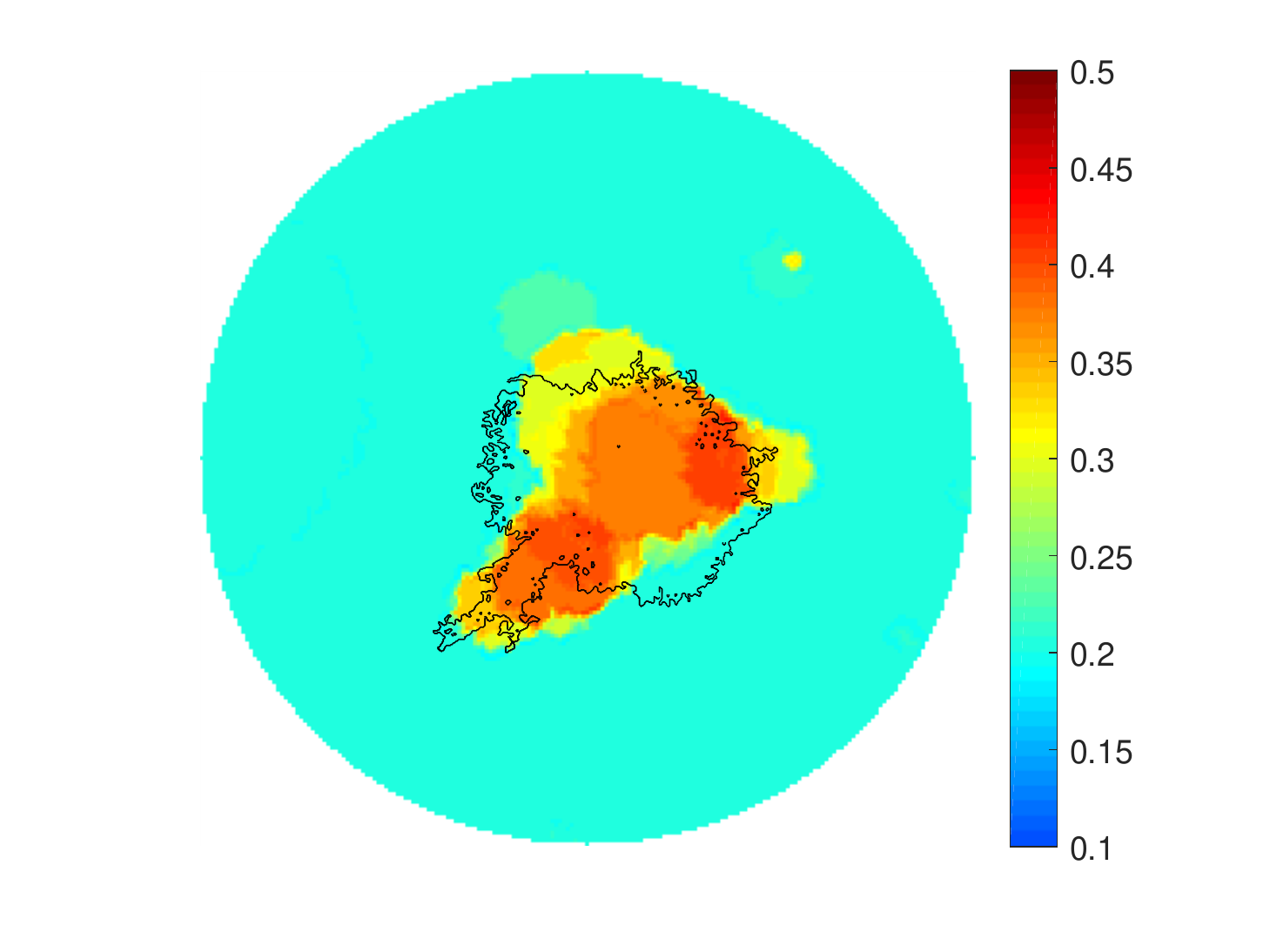}}
  \subfigure[CD: Step~3]{\includegraphics[width=0.33\textwidth]{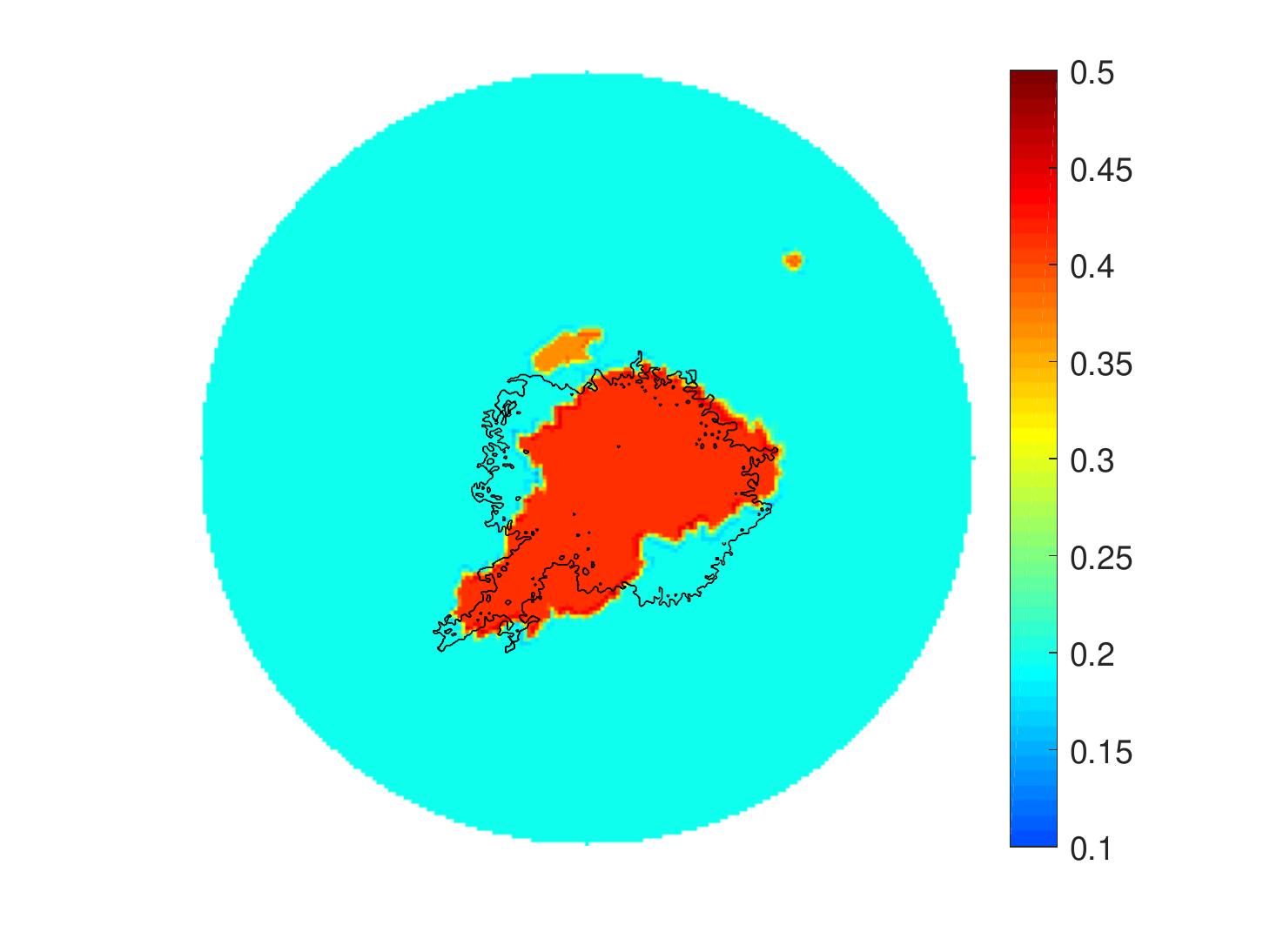}}}
  \end{center}
  \caption{(a)~EIT model~\#4: true electrical conductivity $\sigma_{true}(x)$.
    (b,c,e,f)~Solution images obtained by (b,c)~{\tt SNOPT} and (e,f)~CD after completing
    (b,e)~Step~2 and (c,f)~Step~3. The black dots are added to represent the location of
    the cancer-affected region taken from known $\sigma_{true}(x)$ in (a).
    (d)~Solution errors $\| \sigma^k - \sigma_{true} \|_{L_2}$
    as functions of a number of cost functional evaluations evaluated while employing
    (black) CD and (red) {\tt SNOPT} optimizers. Blue dots represent solutions obtained
    after Step~2 and 3 phases are complete.}
  \label{fig:model_7}
\end{figure}

\begin{figure}[htb!]
  \begin{center}
  \mbox{
  \subfigure[model~\#5]{\includegraphics[width=0.33\textwidth]{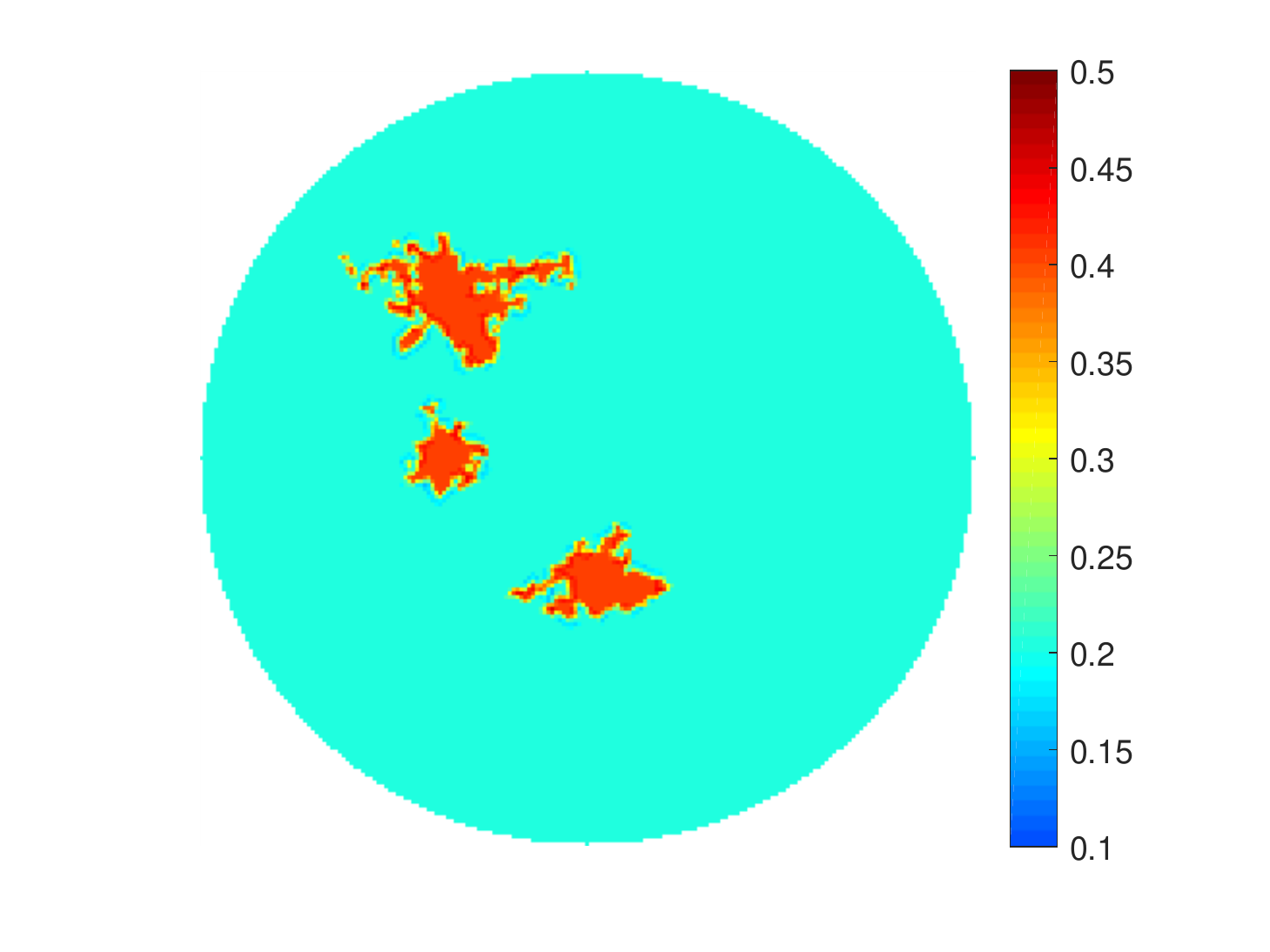}}
  \subfigure[SNOPT: Step~2]{\includegraphics[width=0.33\textwidth]{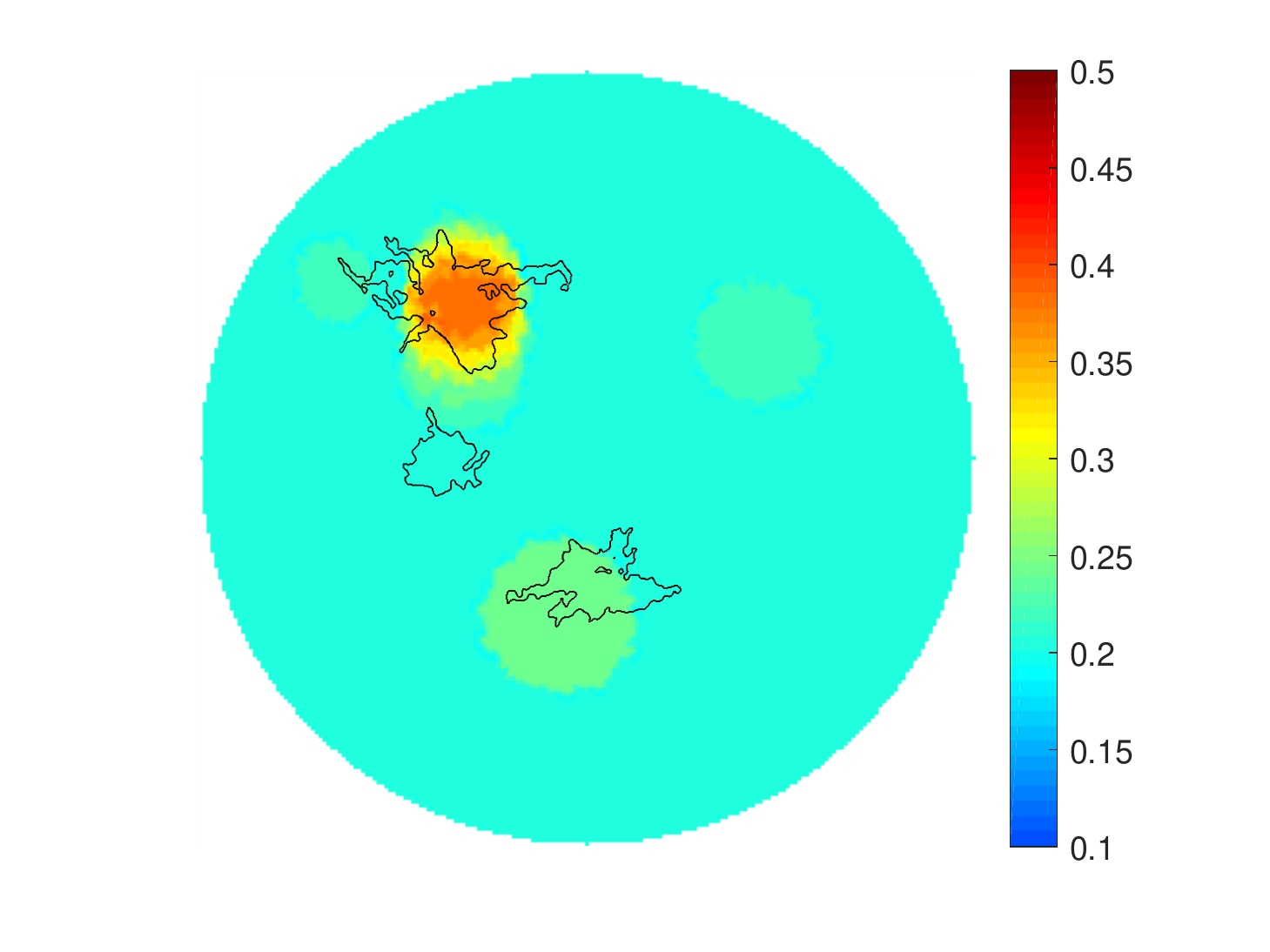}}
  \subfigure[SNOPT: Step~3]{\includegraphics[width=0.33\textwidth]{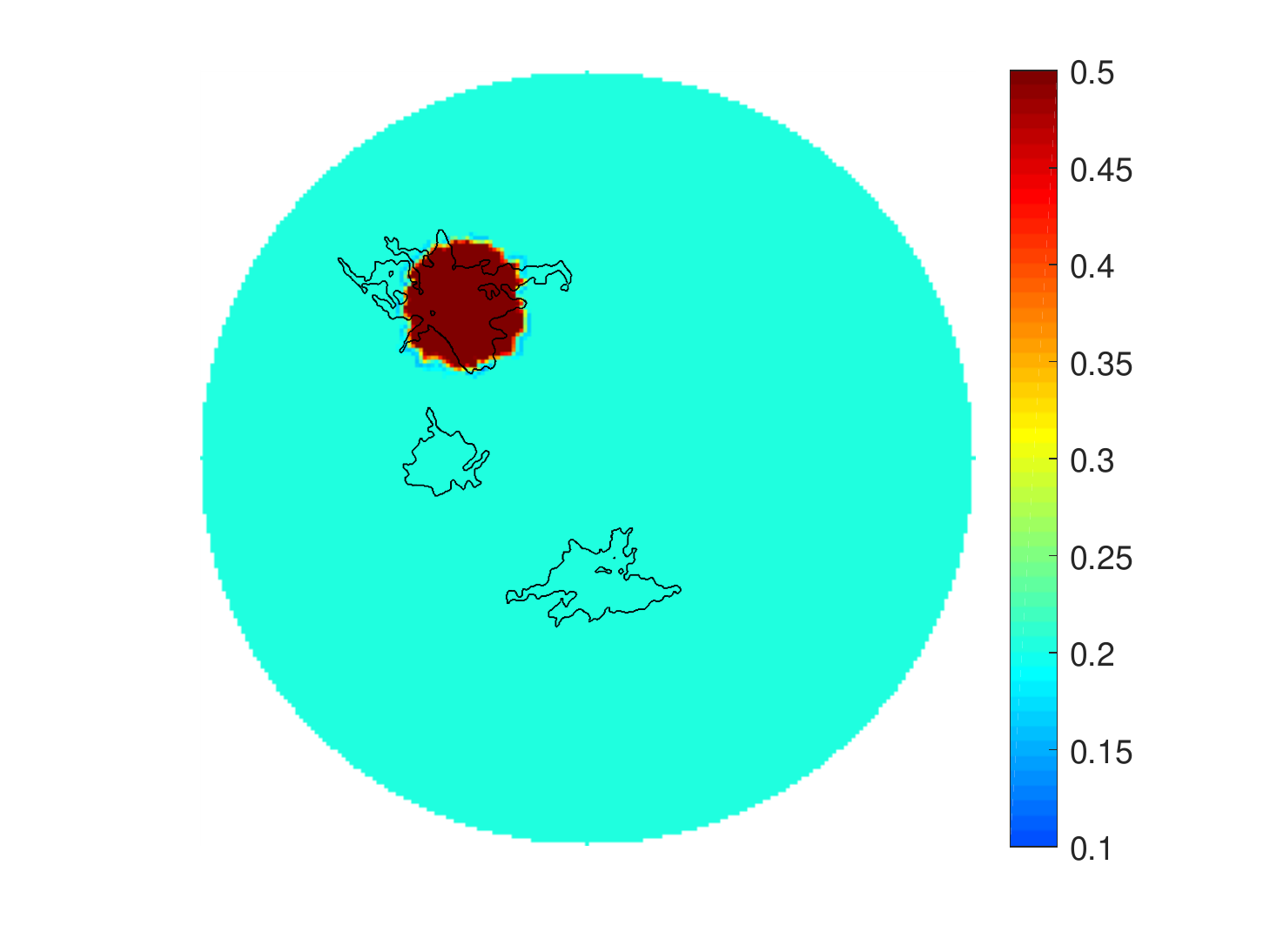}}}
  \mbox{
  \subfigure[solution error]{\includegraphics[width=0.33\textwidth]{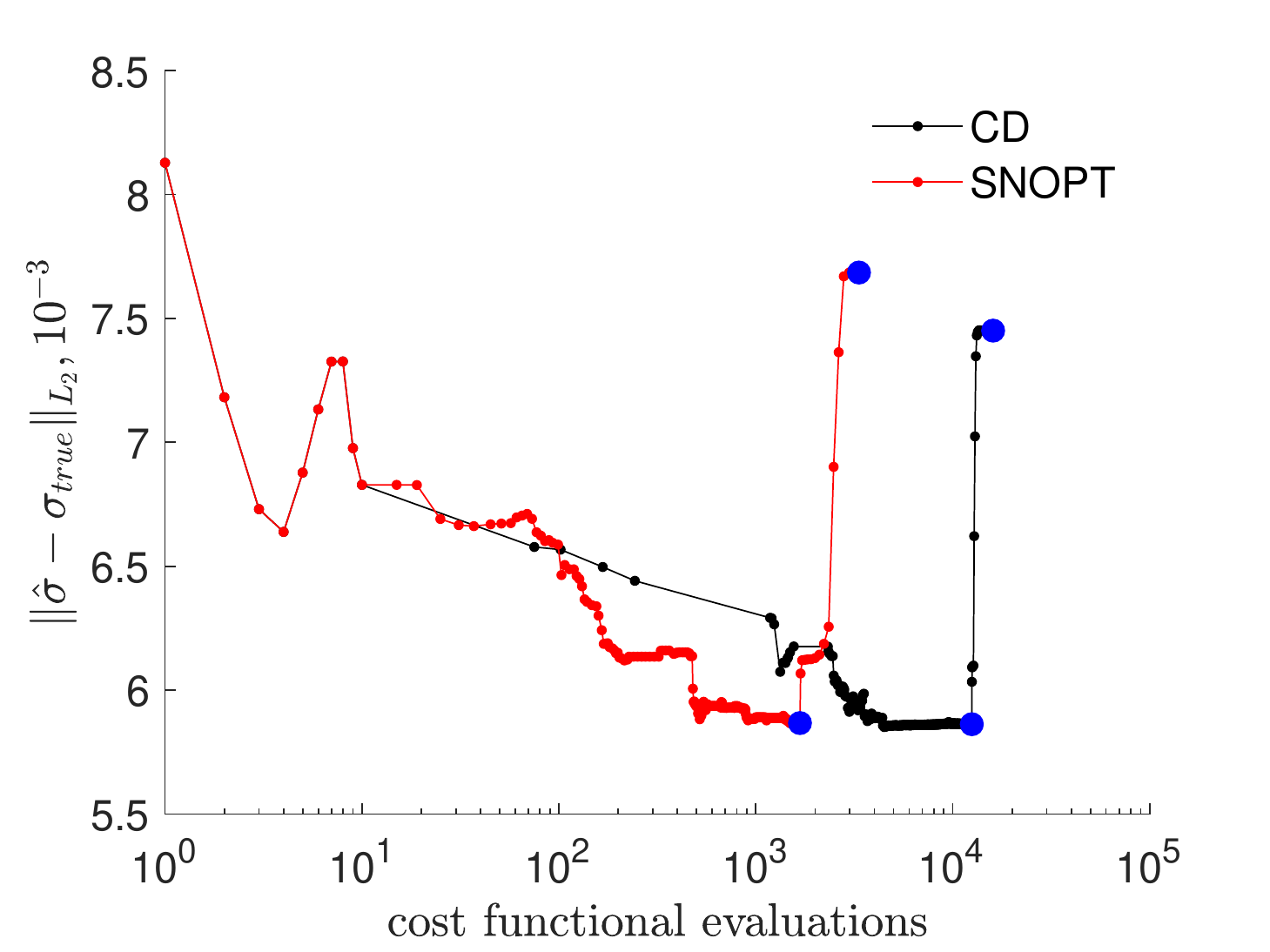}}
  \subfigure[CD: Step~2]{\includegraphics[width=0.33\textwidth]{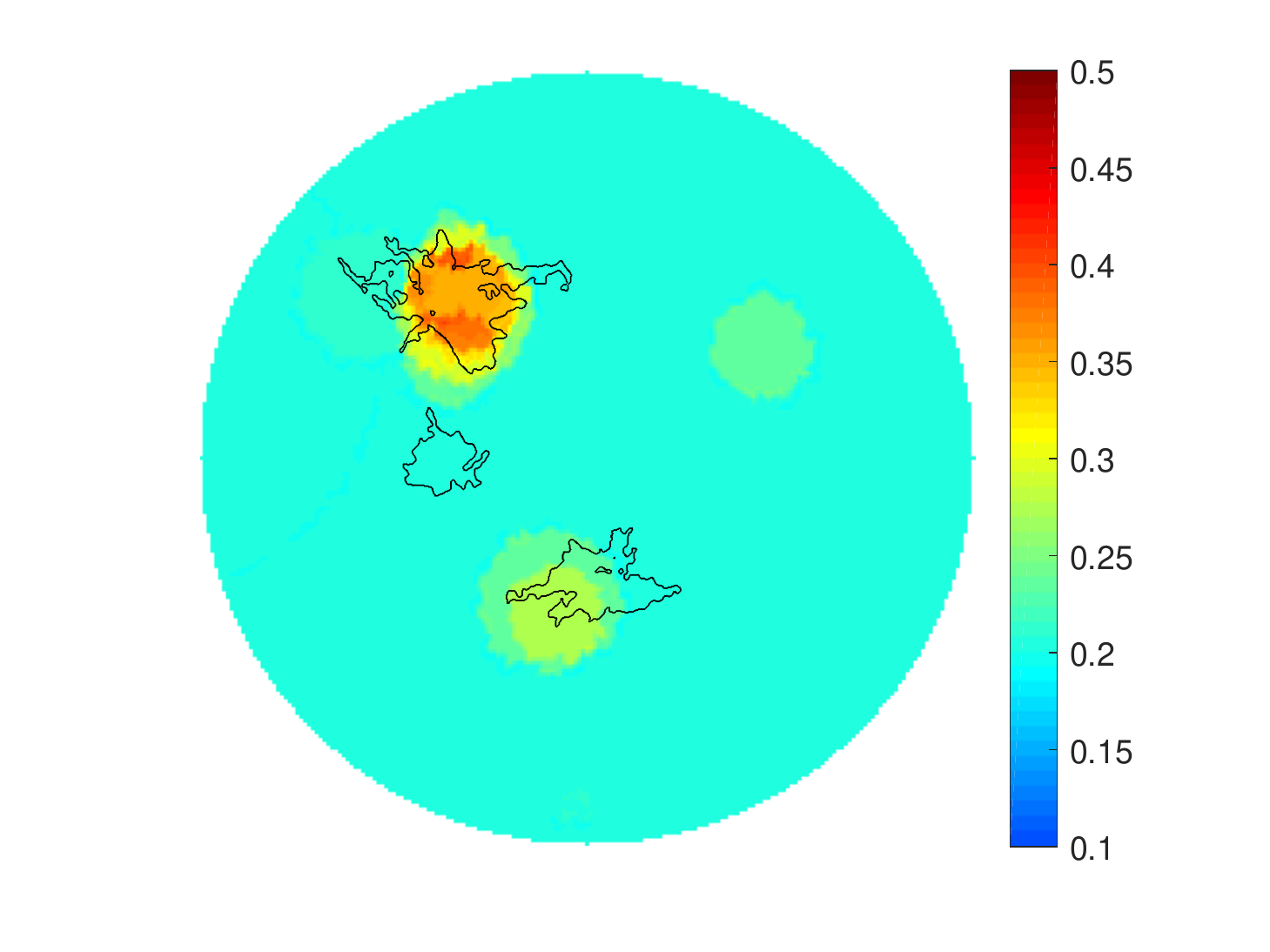}}
  \subfigure[CD: Step~3]{\includegraphics[width=0.33\textwidth]{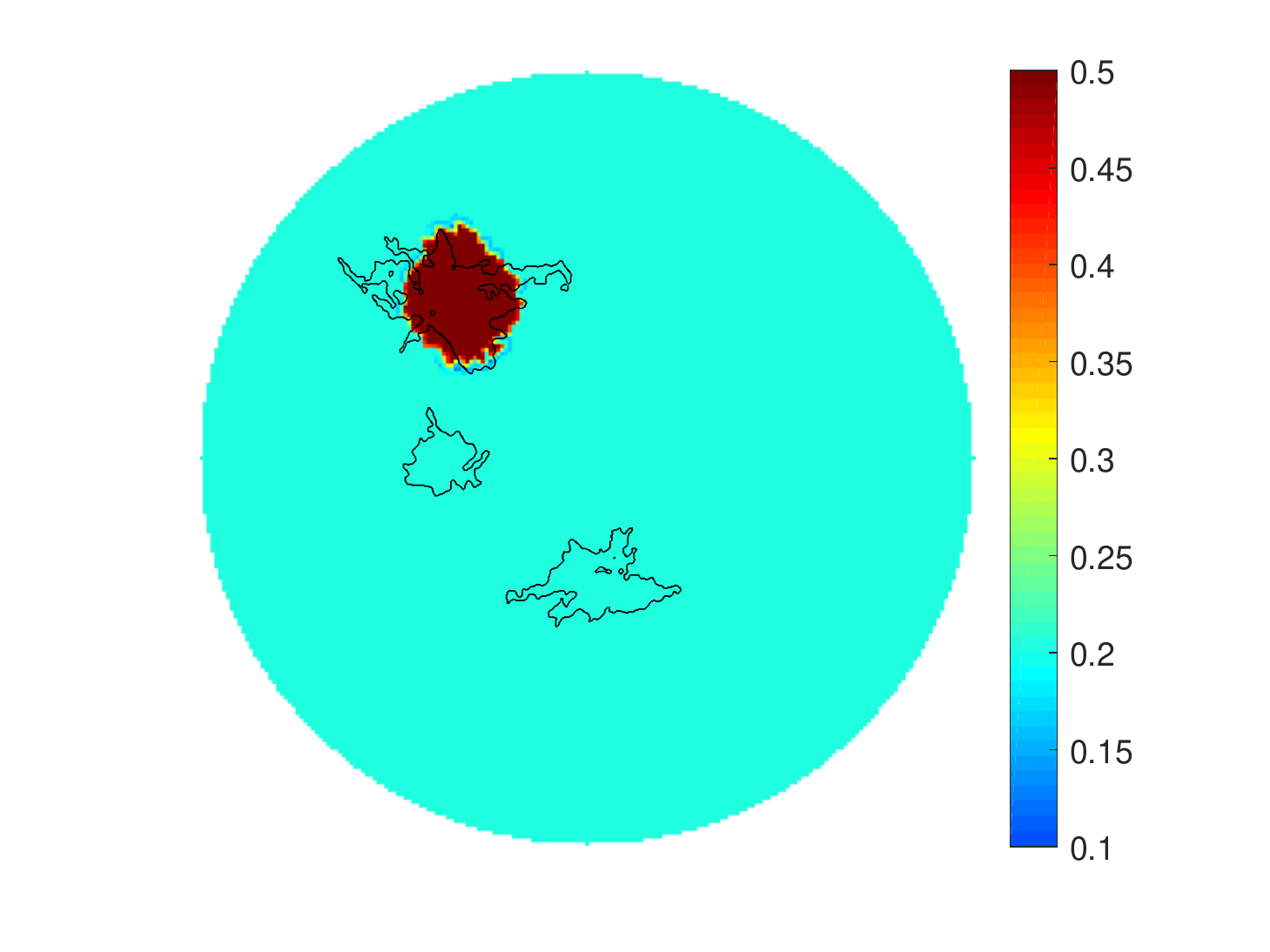}}}
  \end{center}
  \caption{(a)~EIT model~\#5: true electrical conductivity $\sigma_{true}(x)$.
    (b,c,e,f)~Solution images obtained by (b,c)~{\tt SNOPT} and (e,f)~CD after completing
    (b,e)~Step~2 and (c,f)~Step~3. The black dots are added to represent the locations of
    cancer-affected regions taken from known $\sigma_{true}(x)$ in (a).
    (d)~Solution errors $\| \sigma^k - \sigma_{true} \|_{L_2}$
    as functions of a number of cost functional evaluations evaluated while employing
    (black) CD and (red) {\tt SNOPT} optimizers. Blue dots represent solutions obtained
    after Step~2 and 3 phases are complete.}
  \label{fig:model_16}
\end{figure}

In the last part of our numerical experiments with the proposed optimization
framework, our particular interest is in applying it to cases seen in the
medical practice during cancer-related screening procedures. We utilize our
last two models (\#4 and \#5) based on the mammogram~(X-ray) and magnetic
resonance~(MRI) images, respectively, of real breast cancer cases available
in \cite{Bassett2003,Weinstein2009}; refer to \cite{ChunEdwardsBukshtynov2024}
for the details on creating these models by converting the actual images to
their binary versions and obtaining synthetic data in place of the true
measurements. Model~\#4 shows an invasive ductal carcinoma with an irregular
shape and spiculated margins; see Figure~\ref{fig:model_7}(a). Our final model
(\#5) refers to even more complicated cases seen in the medical practice when
multiple regions suspected of being cancerous are present and characterized by different
sizes and nontrivial shapes; see Figure~\ref{fig:model_16}(a). This model shows
multiple (three) spots also identified as invasive ductal carcinoma with irregular
shapes and spiculated margins.

Figures~\ref{fig:model_7}(b) and \ref{fig:model_7}(e) demonstrate the images obtained
while running fine-scale optimization driven by gradient-based {\tt SNOPT} and
derivative-free CD optimizers, respectively. Although the solution to the inverse EIT
problem of model~\#4 is very challenging due to the nontrivial shape of the cancerous
spot at the center, the quality of both images is rather good and further significantly
improved after applying binary tuning in Step~3; refer to Figures~\ref{fig:model_7}(c)
and \ref{fig:model_7}(f). Here, consistent with the previous results, {\tt SNOPT} conducts
the reconstruction better in quality (better match for both color and shape) and
computational speed (995 vs.~14,099 cost functional evaluations to complete both Steps~2
and 3), confirmed by the analysis based on the solution error; see
Figure~\ref{fig:model_7}(d). The results obtained here by gradient-based {\tt SNOPT}
are apparently of much higher quality than those reported in
\cite{ChunEdwardsBukshtynov2024} obtained by gradient- and PCA-based multiscale
optimization paired with the added regularization.

Finally, we refer to the results obtained for the most complicated case by model~\#5.
Figures~\ref{fig:model_16}(b,c) and \ref{fig:model_16}(e,f) present the images produced
by {\tt SNOPT} and CD optimizers, respectively. We admit the same level of accuracy in
reconstructing two spots and the biggest one after Steps~2 and 3, respectively, by both
optimizers. Figure~\ref{fig:model_16}(d), showing the solution error, also confirms more
or less the same performance for both. This model indeed sets the limits for the current
implementation of the proposed computational framework due to the nontrivial shapes of
the cancerous spots and their small sizes. Although the overall performance is fairly
modest, the final results here are of much better quality than reported previously in
\cite{ChunEdwardsBukshtynov2024}, where extra data was added through applied regularization.
It makes the application of our computational framework very promising, having enough
space for future developments in both quality of the reconstructed images related to
very complex medical cases and computational efficacy sought to help implement the
approach for use in the medical practice.

\section{Concluding Remarks}
\label{sec:remarks}

In this work, we developed and validated a highly efficient computational framework
for reconstructing images of near-binary types related to physical properties of
various models used in biomedical applications, e.g., to assist with the recognition
of cancer-affected regions while solving an inverse problem of cancer detection using
the electrical impedance tomography technique. In particular, we explored the
possibility of applying the proposed solution methodology to the IPCD problems to detect
cancerous regions surrounded by healthy tissues, which also have small sizes and
boundaries of irregular shapes that appeared to be vital for early cancer detection and
easy control of the dynamics of cancer development or treatment progress. We also prove
its suitability for models of various complexity seen in diverse applications for
biomedical sciences, physics, geology, chemistry, etc. The efficiency of the new approach
in both computational speed and accuracy is achieved by combining the advantages of recently
developed optimization methods that use sample solutions with customized geometry and
multiscale control space reduction, all paired with gradient-based techniques. A practical
implementation of this approach has an easy-to-follow design tuned by a nominal number of
parameters to govern the entire suite of computational and optimization facilities. We
conclude on the high potential of the proposed computational methodology to minimize
possibilities for false positive and false negative screening and improve the overall
quality of EIT-based procedures.

Despite the superior performance of the proposed framework, there are many ways this
optimization algorithm can be tested and further extended, e.g., by applying advanced
minimization techniques to perform local and global searches, using adaptive schemes
for a smooth transition of the obtained solutions between computational phases,
flexible (efficient) termination criteria, and different metrics to evaluate the quality
of the overall reconstruction. We plan further research on the
measurement structure, e.g., considering 32-electrode schemes in 2D, shifting to 3D models,
and improving sensitivity by optimizing the configuration of available data.
We also expect to have more benefits for saving computational time from applying parallelization
to the entire computational framework, including solving forward EIT problems. Finally, we
are interested in other expansions leading to reliable and accurate results for near-real
cases by using bimodal distributions and fully anisotropic models.

\subsubsection*{Acknowledgements}
We wish to thank the anonymous reviewers for their valuable comments and
suggestions to improve the clarity of the presented approach and the overall
readability of this paper.

\bibliographystyle{spmpsci}
\bibliography{biblio_Bukshtynov,biblio_EIT,biblio_OPT}

\end{document}